\documentstyle[12pt]{amsart}
\let\text=\mbox
\newcommand{\skipthistext}[1]{}

\newcommand{\R}{{\mathbb R}}

\newcommand{\Q}{{\mathbb Q}}

\renewcommand{\pf}{\noindent {\bf Proof} \hspace{2mm}}

              \newcommand{\W}{{\cal W}}
              \newcommand{\U}{{\cal U}}
              \newcommand{\T}{{\cal T}}
              \newcommand{\V}{{\cal V}}

              \renewcommand{\O}{{\cal O}}
              
              \newcommand{\G}{{\cal G}}

\newcommand{\abs}[1]{\centerline{\sc Abstract}
\vspace{3mm}
\centerline{\parbox{150mm}{\small #1}}}
\newcounter{thm}

\def\thmnumber{
\addtocounter{thm}{1}\if\arabic{section}0\relax
\else\arabic{section}.\fi%
\if\arabic{subsection}0\relax
\else\arabic{subsection}.\fi%
\arabic{thm}
\hspace{2mm}}

\def\setthmnumber{\setcounter{thm}{0}}

\newenvironment{thm}{
\vspace{3mm}\par\noindent
{\bf Theorem \thmnumber} \it}
{\vspace{3mm}\par\rm}

\newenvironment{defi}{
\vspace{3mm}\par\noindent
{\bf Definition \thmnumber}\it}
{\vspace{3mm}\par\rm}

\newenvironment{lem}{
\vspace{3mm}\par\noindent
{\bf Lemma \thmnumber} \it}
{\vspace{3mm}\par\rm}

\newenvironment{prop}{
\vspace{3mm}\par\noindent
{\bf Proposition \thmnumber}\it}
{\vspace{3mm}\par\rm}

\newenvironment{re}{
\vspace{3mm}\par\noindent
{\bf Remark \thmnumber}\it}
{\vspace{3mm}\par\rm}

\newcommand{\sectioni}[1]{\setthmnumber\setcounter{subsection}{0}
\section{#1}}

\renewcommand{\subsection}[1]{%
\vspace{7mm}\par\noindent%
\addtocounter{subsection}{1}%
{\sc\arabic{section}.\arabic{subsection} \ {#1}}
\vspace{4mm}\par\noindent\setthmnumber}

\begin{document}
\title{On a notion of maps between orbifolds\\
I. function spaces}
\author{Weimin Chen}
\date{\today\\
\hspace{2mm}{\it Keywords}. Orbifold. Mapping space. 
}
\maketitle

\abs{This is the first of a series of papers which are devoted to a 
comprehensive theory of maps 
between orbifolds. In this paper, we define the maps in the more general
context of orbispaces, and establish several basic results concerning 
the topological structure of the space of such maps. In particular, we
show that the space of such maps of $C^r$ class between smooth orbifolds 
has a natural Banach orbifold structure if the domain of the map is 
compact, generalizing the corresponding result in the manifold case. 
Motivations and applications of the theory come from string theory
and the theory of pseudoholomorphic curves in symplectic orbifolds.
}

\sectioni{Introduction}

The space of differentiable maps is a fundamental object associated to
smooth manifolds. More precisely, given smooth manifolds $M$ and $N$, 
one considers the set $C^r(M,N)$ of $C^r$-maps $f:M\rightarrow N$, ie. 
maps with continuous partial derivatives up to a given order $r$. It is 
well-known that $C^r(M,N)$ has a natural topology, called the Whitney 
topology, which in the case when $M$ is compact 
also gives $C^r(M,N)$ an infinite dimensional manifold (ie. Banach manifold) 
structure. A particularly important case is when $M$ is the circle $S^1$, in
this case the space of maps is called the loop space of $N$. These function 
spaces are fundamental in many branches of mathematics, including algebraic 
and differential topology of manifolds, differential geometry and global 
analysis, as well as mathematical physics. 

The purpose of this paper is to formulate a notion of maps between orbifolds,
which is given in the more general perspective of orbispaces, such that 
orbifolds form a category under such maps. Several basic results concerning 
topological structure of the space of maps were established. In particular, 
we showed that the space of $C^r$-maps between smooth orbifolds $X$, 
$X^\prime$ is a Banach orbifold when $X$ is compact. Based on
these structural results, we developed in the sequel \cite{C1} a basic 
machinery for studying homotopy classes of such maps. 

Considerations in this work originated from construction of Gromov-Witten
invariants of symplectic orbifolds \cite{CR1, CR2} and from our attempt to
understand some of the mathematical implication of the work of Dixon, Harvey,
Vafa and Witten \cite{DHVW} on string theories of orbifolds. While the
treatment in \cite{CR1, CR2} was {\it ad hoc} in nature, the aforementioned
structural theorem on mapping spaces provided a solid foundation for the 
Fredholm 
theory of pseudoholomorphic curves in symplectic orbifolds, which was more 
instrumental in dimension $4$ in the recent, more geometric application in 
\cite{C2, C3}. On the other hand, 
concerning the orbifold string theories it is worth noting that, with the 
notion of maps in this paper, the loop space of an orbifold gives exactly 
the configuration space of strings considered by the physicists in 
\cite{DHVW}, where the orbifolds are global quotients $X=Y/G$. Given the fact 
that spaces of differentiable maps have been involved in many considerations 
in geometry, topology and mathematical physics, it is our hope that this work 
will pave the way for similar considerations, especially those originating 
from mathematical physics, in the orbifold category.

\vspace{3mm}

{\it A brief history of orbifolds.} \hspace{2mm}
The concept of orbifold was first introduced by Satake under the name
``V-manifold'' in his 1957 paper \cite{Sa}, where orbifolds were perceived
as a class of singular spaces which can be modeled locally on a smooth
manifold modulo a finite group. The smooth manifold together with the group
action is called a local uniformizing system. The purpose of Satake was to
demonstrate that the basic differential geometry of smooth manifolds can be
suitably extended to this class of singular spaces. The local structure of
orbifold, ie. being locally modeled as the quotient of a smooth manifold by 
a finite group action, was merely used as sort of generalized smooth structure 
here. This was clearly reflected in the notion of ``V-manifold map'' 
introduced by
Satake, which was roughly speaking a continuous map with smooth liftings to
local uniformizing systems. The basic intuition here is that orbifolds behave
very much like smooth manifolds, as long as only $\Q$ coefficients are
concerned.

The more popular name ``orbifold'' was due to Thurston \cite{Th}, who 
rediscovered this concept in the late 70's as a useful device in studying
geometric structures on $3$-manifolds. Here is the basic idea. The quotient
manifold $Y/\Gamma$ provides a useful, geometric device for studying a free,
proper, and discontinuous action of $\Gamma$ on $Y$. In order to extend this
to the case where the action may not be free, one has to allow singularities 
in the quotient space $Y/\Gamma$. The key issue here, however, is how to 
recover the action from the quotient space when it is not free. To this end, 
Thurston introduced the notion of fundamental group of an orbifold, which 
has the property that when $Y$ is simply-connected, the fundamental group of 
the orbifold $Y/\Gamma$ is always $\Gamma$, regardless the action is free or
not. This is obviously different from the fundamental group of $Y/\Gamma$
regarded as a singular topological space. Thurston's notion of fundamental 
group proves to be particularly useful in a context where the orbifold 
arises as the quotient space of a group action. 

There are two major sources of orbifolds. One is as the quotient space of
either a proper discontinuous action or a compact Lie group action with only
finite isotropy subgroups. The other is as a space or a variety with only 
``mild'' singularities (for instance, as those appeared in birational 
transformation or in 
degeneration of certain Riemannian metrics), where there is not necessarily
a natural global group action involved. Despite the abundance of appearance
in the literature, orbifolds have been mainly involved as a useful technical
device rather than an object of permanent interest, and have been treated 
in a rather practical, {\it ad hoc} manner, keeping in mind the basic 
observations of Satake and Thurston. In particular, there has been no 
comprehensive study or theory on orbifolds (eg. as the one we had on 
manifolds). 

There is a somewhat formal and more contemporary formalism of orbifolds 
(compared with that in Satake \cite{Sa} or Thurston \cite{Th}), using the 
categorical language of \'{e}tale (smooth) topological groupoids (cf. eg. 
\cite{BH}). One of the advantages of this formalism is that canonical, 
functorial constructions on topological groupoids provided additional tools 
associated to the orbifolds (these, of course, work for manifolds also, but 
only become nontrivial when applied to orbifolds). For example, one can
associate an orbifold with a classifying space (which is the classifying space 
of the associated  \'{e}tale topological groupoids, and which is unique up to
homotopy equivalence), and define topological invariants of orbifolds 
through the classifying spaces, cf. eg. \cite{Ha2}. Thurston's orbifold
fundamental group may be recovered as the fundamental group of the classifying
space. This approach is particularly useful in the case of cofunctors,
such as cohomologies, fibre bundles, characteristic classes, etc.

\vspace{3mm}

{\it String theories on orbifolds.}\hspace{2mm}
The suggestion of having a comprehensive study of orbifolds (as the one 
for manifolds) seemed to come first from string theory. In 1985, physicists 
Dixon, Harvey, Vafa and Witten \cite{DHVW} considered string theories on a
compact Calabi-Yau manifold $Y$ equipped with a finite group action of $G$
preserving the Calabi-Yau structure of $Y$ (more generally, $Y$ could be
noncompact with the quotient $X=Y/G$ being a compact orbifold). For the 
purpose of symmetry breaking, it was necessary to not only consider strings
$y(t)$ which satisfies a periodic boundary condition, but also to consider
those satisfying boundary conditions which are periodic up to the action 
of $G$:
$$
y(t+2\pi)=g\cdot y(t), \mbox{  for some } g\in G  \leqno (1.1)
$$
(These more general boundary conditions are called ``twisted boundary
conditions''.) An interesting idea in \cite{DHVW} was to study such string 
theories on $Y$ as a (closed) string theory on the quotient orbifold $X=Y/G$.

One of the advantages of introducing string theories on orbifolds is that
string propagation on an orbifold may be regarded as an arbitrarily good 
approximation to the string propagation on any of the smooth resolutions
of the orbifold. Given such a relation between string theories on a 
Calabi-Yau orbifold $X$ and any of its crepant resolutions $\tilde{X}$
(if there exists one), the physicists made some remarkable prediction 
concerning the Euler number $e(\tilde{X})$ of $\tilde{X}$ in terms of 
the orbifold $X$.

More concretely, in string theories on a smooth manifold, 
the Euler number of the manifold 
may be interpreted as twice of the ``number of generations'' in the 
physical theory. Extending this to orbifolds, one may simply define 
the ``stringy Euler number'' of an orbifold to be twice of the ``number 
of generations''  in the physical theory. With this understood,
the physicists derived, using the path integral method, the following formula 
for the ``stringy Euler number'' of the Calabi-Yau orbifold $X=Y/G$:
$$
e_{string}(X)=\frac{1}{|G|}\sum_{hg=gh}e(Y^{<g,h>}) \leqno (1.2)
$$
where $Y^{<g,h>}$ denotes the common fixed-point set of $g$ and $h$ in $Y$.
By the said relation between string theories on $X$ and any of its crepant
resolutions $\tilde{X}$, the physicists then obtained the following identity
$$
e(\tilde{X})=\frac{1}{|G|}\sum_{hg=gh}e(Y^{<g,h>}). \leqno (1.3)
$$
The right-hand side of $(1.2)$ was later reformulated by Hirzebruch and
H\"{o}fer \cite{HH}, and consequently $(1.3)$ becomes 
$$
e(\tilde{X})=e(X)+\sum_{(g),g\neq 1} e(Y^g/Z(g)), \leqno (1.4)
$$
where $Y^g$ is the fixed-point set of $g$ in $G$, $(g)$ stands for the
conjugacy class of $g$ in $G$ and $Z(g)$ is the centralizer of $g$ in $G$.
The advantage of the right-hand side of $(1.4)$ over that of $(1.2)$ is 
that: (1) it explicitly shows that the ``stringy Euler number'' 
$e_{string}(X)$ of an orbifold $X$ differs in general from the Euler number 
$e(X)$, (2) the right-hand side of $(1.4)$ can be defined for any orbifold, 
which is not necessarily of the form $Y/G$ with $G$ finite.

The physicists's prediction $(1.3)$ (or equivalently $(1.4)$) was soon 
related to the so-called ``McKay correspondence'' in mathematics (cf.
\cite{McK}), and has since stimulated a great deal of interests in this
subject. See \cite{Reid} for a recent survey on McKay correspondence.

Similar predictions may be made on other types of invariants, such as
Gromov-Witten invariants, elliptic cohomology, etc.,  which also have
an alleged string theory interpretation (cf. \cite{V, W, Se}). In order 
to fully explore in this direction of research, it is desirable to develop 
a comprehensive theory of orbifolds, which would provide a framework within 
which these ``stringy invariants'' of orbifolds may be properly interpreted 
and further studied. The theory of maps between 
orbifolds developed in this work was motivated by such a desire. It turns
out that this theory of maps (1) is amenable to techniques of differential
geometry and global analysis on orbifolds as pioneered by Satake, (2)
reflects the distinct topological structure of orbifolds as discovered by
Thurston, and (3) is consistent to the considerations in string theories 
of orbifolds by the physicists in \cite{DHVW}.  

\vspace{3mm}

{\it Gromov-Witten invariants of orbifolds.}\hspace{1mm}
The theory of maps between orbifolds developed in this work originated 
in the construction of Gromov-Witten invariants of symplectic orbifolds
in \cite{CR1, CR2}. 

Recall that Gromov-Witten invariants of a symplectic manifold $(M,\omega)$ 
come as a certain algebraic count of $J$-holomorphic curves in $M$ for some 
fixed, $\omega$-tamed almost complex structure $J$. The actual counting of
$J$-holomorphic curves goes roughly as follows. One introduces a Banach
manifold $B$, which is the space of maps of a certain fixed Sobolev or 
H\"{o}lder type from a Riemann surface into $M$, and a Banach bundle 
$E\rightarrow B$ such that the set of $J$-holomorphic curves in question
is given as the zero set of a Fredholm section $s:B\rightarrow E$. In nice 
situations, one can show that for a generic choice of $J$ the Fredholm section 
$s$ intersects transversely with the zero section of $E$, and consequently, as 
the zero set the moduli space of $J$-holomorphic curves is a smooth manifold 
whose dimension can be computed via the Riemann-Roch theorem. Furthermore,
one can compactify the moduli space using Gromov's compactness theorem, and 
again in nice situations one can show that the compactification of the moduli 
space has only codimension at least $2$ ``boundary'' components, thus has 
a well-defined ``fundamental class'' in the Banach manifold $B$. In these 
nice situations, the Gromov-Witten invariants are defined by evaluating 
certain ``universal'' cohomology classes 
on the Banach manifold $B$ against the fundamental class of the compactified 
moduli space of $J$-holomorphic curves, which turns out to be independent of 
the choice of the almost complex structure $J$. (For details, see eg. 
\cite{McDS}.) In any event, such a setup in the Banach manifold framework is 
the so-called Fredholm theory of pseudoholomorphic curves.

The above construction in the manifold case does not generalize to orbifolds
in any obvious way because there is a lack of corresponding theory of maps
between orbifolds, and thus no Fredholm theory is available. More concretely,
one may regard a $J$-holomorphic curve in a symplectic orbifold as the 
image of a $J$-holomorphic V-manifold map (in the sense of Satake) from an 
orbifold Riemann surface into the symplectic orbifold. The problem is that
the corresponding space of V-manifold maps is not known to have any infinite
dimensional manifold or orbifold structure. On the other hand, there is also a 
related but more conceptual issue due to the lack of a good theory of maps (or
in a more general sense, due to a lack of comprehensive theory of orbifolds),
that is, how to determine whether the constructed Gromov-Witten invariants of 
symplectic orbifolds (no matter how technically the construction was done) 
actually give the corresponding ``stringy invariants'' in the orbifold string 
theory \cite{DHVW}.

The problem of orbifold Gromov-Witten theory was solved in \cite{CR1, CR2}
at the technical level. A crucial step in the construction is to add an
additional piece of data to a V-manifold map, which is an isomorphism class of 
pull-back bundles by the V-manifold map. The secret behind this is that
one needs to specify a deformation type for each $J$-holomorphic curve in
order to ``count'' the $J$-holomorphic curves properly. 
The resulting Gromov-Witten invariants 
turn out to be a certain algebraic count of $J$-holomorphic curves together
with a specified deformation type in the orbifold. Formally, one would have
defined the Gromov-Witten invariants in \cite{CR1, CR2} by working with the 
space of pairs $(f,\xi)$, where $f$ is a V-manifold map and $\xi$ is an 
isomorphism class of pull-back bundles by $f$. However, no Fredholm theory 
based on the space of pairs $(f,\xi)$ was available back then, so the 
construction in \cite{CR1, CR2} was carried out in a rather {\it ad hoc} 
manner. 

The pair $(f,\xi)$ in the preceding paragraph is a prototype of the maps 
defined in this paper. In this regard, the theory of maps developed in
this work has accomplished the following two goals: (1) providing a
necessary foundation for the Fredholm theory of pseudoholomorphic curves
in symplectic orbifolds, and (2) providing a mathematical framework for
the stringy interpretation of the Gromov-Witten invariants constructed 
in \cite{CR1, CR2}. In particular, the construction of Gromov-Witten 
invariants of symplectic orbifolds may be done in the same line as in 
the smooth case, and the one in \cite{CR1, CR2} can be substantially 
simplified. 

In dimension $4$, Gromov's pseudoholomorphic curve theory has had more 
geometric applications (cf. eg. \cite{McDS}), where the manifold structure 
of the corresponding moduli space of 
pseudoholomorphic curves played a fundamental role. (To the contrary,
in the problem of counting pseudoholomorphic curves one only needs a 
well-defined fundamental class of the moduli space.) In order to apply
similar ideas to a situation where quotient singularities are present,
a Fredholm theory in the orbifold context must be in place first.
For a theory of pseudoholomorphic curves in symplectic $4$-orbifolds and 
some related applications, see the recent survey \cite{C4}, and for more 
details, see \cite{C2, C3}.

\vspace{3mm}

{\it Extension of the equivariant category.} \hspace{1mm}
Define the ``equivariant category'' as follows: the objects are pairs $(Y,G)$,
where $Y$ is a smooth manifold equipped with a smooth action of a finite group 
$G$, and the morphisms from $(Y,G)$ to $(Y^\prime, G^\prime)$ are pairs
$({f},\rho)$, where ${f}:Y\rightarrow Y^\prime$ is a differentiable 
map and $\rho:G\rightarrow G^\prime$ is a homomorphism, such that ${f}$
is $\rho$-equivariant, ie. ${f}(g\cdot y)=\rho(g)\cdot {f}(y)$
for any $y\in Y$ and $g\in G$. Then the category of orbifolds, where the
objects are smooth orbifolds and the morphisms are differentiable maps 
between orbifolds defined in this paper, may be regarded as an extension of
the equivariant category in the following sense.

Given any two objects $(Y,G)$ and $(Y^\prime, G^\prime)$ of the equivariant 
category, let $X=Y/G$ and $X^\prime=Y^\prime/G^\prime$ be the corresponding
orbifolds. Let $[(Y,G);(Y^\prime, G^\prime)]$ denote the set of morphisms 
$({f},\rho):(Y,G)\rightarrow (Y^\prime, G^\prime)$, and let $[X;X^\prime]$
denote the set of morphisms $\Phi:X\rightarrow X^\prime$. Then each 
$({f},\rho)\in [(Y,G);(Y^\prime, G^\prime)]$ defines a differentiable map 
$\Phi\in [X;X^\prime]$: $({f},\rho)\mapsto \Phi$,
which induces an identification as a subset (cf. Lemma 3.1.2)
$$
[(Y,G);(Y^\prime, G^\prime)]/G^\prime\subset [X;X^\prime], \leqno (1.5)
$$
where the action of $G^\prime$ on $[(Y,G);(Y^\prime, G^\prime)]$ is given by
$$
g^\prime\cdot ({f},\rho)=(g^\prime\circ {f}, Ad(g^\prime)\circ\rho),\;
\forall g^\prime\in G^\prime, ({f},\rho)\in [(Y,G);(Y^\prime, G^\prime)].
$$
Moreover, each $\Phi:X\rightarrow X^\prime$ induces a homomorphism between 
Thurston's orbifold fundamental groups 
$\Phi_\ast:\pi_1(X)\rightarrow\pi_1(X^\prime)$, such that $\Phi$ is defined by
a pair $({f},\rho)$ if and only if $\Phi_\ast$ maps the subgroup $\pi_1(Y)$
of $\pi_1(X)$ into the subgroup $\pi_1(Y^\prime)$ of $\pi_1(X^\prime)$
(cf. Lemma 2.4.1 in \cite{C1}). Finally, for any $r>0$, let 
$[(Y,G);(Y^\prime, G^\prime)]^r$
be the set of pairs $({f},\rho)$ where ${f}$ is a $C^r$-map, and
let $[X;X^\prime]^r$ be the space of $C^r$-maps from $X$ to $X^\prime$.
Then when $Y$ is compact, $[(Y,G);(Y^\prime, G^\prime)]^r/G^\prime$ is 
naturally a Banach orbifold, and since $X=Y/G$ is also compact, 
$[X;X^\prime]^r$ is a Banach orbifold by Theorem 1.4 of this paper. 
In this case, $(1.5)$ gives 
$$
[(Y,G);(Y^\prime, G^\prime)]^r/G^\prime \subset [X;X^\prime]^r \leqno (1.6)
$$
as an open and closed Banach suborbifold. 

The above observation suggests that when developing a theory of orbifolds
(based on the theory of maps in this work),
one should look at the corresponding theory in the equivariant category 
first, and then seek a proper generalization to the orbifold category.

An important technique in the equivariant category is given by the various
localization theorems. It would be interesting to see if such a technique 
can be extended to the orbifold category (or more generally the orbispace
category). For any orbifold $X$, one can write $X=Y/G$ for a smooth manifold
$Y$ and a compact Lie group $G$, where $Y$ is the bundle of orthonormal frames
on $X$ and $G=O(n)$ (here $n=\dim X$). Technically, one may reduce a 
``localization problem'' for $X$ to a localization problem for $(Y,G)$, for
example, as done in the derivation of the index theorem over orbifolds in
\cite{Ka1}. However, such an approach is not natural and an intrinsic
localization principle would be more desirable. 

Not every orbifold can be written as a global quotient $Y/G$. For a general
orbifold where $G$ is missing, what would be an appropriate replacement for
the ``structure group'' $G$ ? The answer seems to lie in the notion of 
``complex of groups'' introduced by Haefliger \cite{Ha3}. More concretely,
given any orbifold $X$ (or more generally an orbispace), one can associate
a complex of groups $G(\{U_i\})$ to a given cover $\{U_i\}$ of $X$, where
each $U_i$ is a uniformized open set, such that for any refinement $\{V_k\}$
of $\{U_i\}$, there is a canonically defined homomorphism from $G(\{V_k\})$
to $G(\{U_i\})$. In this way, each orbifold is associated with a direct limit
of complexes of groups, such that a map between orbifolds induces 
a ``morphism''
between the corresponding direct limits of complexes of groups. When $X=Y/G$
is a global quotient and the cover $\{U_i\}$ is taken to be $\{X\}$, the 
complex of groups $G(\{U_i\})$ reduces to the group $G$. See Remark 2.1.2 (7)
for more details. 

We would like to point out that Haefliger \cite{Ha4} has developed the
homological algebra aspects of complexes of groups, generalizing most of 
the notions having been developed for groups, for instance as in the book 
of Brown \cite{Br}. 

\vspace{3mm}

{\it Summary of main results.}\hspace{1mm}
We shall develop the theory of maps between orbifolds in the more general
context of orbispaces and \'{e}tale topological groupoids. The following 
definition of orbispaces is taken from Haefliger \cite{Ha3}. See \S 2.1 
for a review on \'{e}tale topological groupoids.

\begin{defi}
An orbispace is a topological space $X$ equipped with an {\em(}equivalence 
class of {\em)} \'{e}tale topological groupoid  $\Gamma$, such that {\em (1)} 
$X$ is covered by a set of open sets $\{U_i\}$, where 
for each $U_i$ there is a space $\widehat{U_i}$ acted on by a discrete group 
$G_i$ such that $U_i=\widehat{U_i}/G_i$, {\em (2)} the space of units of 
$\Gamma$ is the disjoint union $\bigsqcup_i\widehat{U_i}$ such that the
restriction of $\Gamma$ to each $\widehat{U_i}$ is the product
groupoid $G_i\times\widehat{U_i}$, and {\em (3)} $X$ is the space of
$\Gamma$-orbits.
\end{defi}

The orbispaces studied in this work are required to further satisfy
certain technical conditions, which are given in terms of the defining 
\'{e}tale topological groupoid $\Gamma$.

\vspace{2mm}

\noindent{\bf Technical Assumptions:}
{\em
\begin{itemize}
\item [{(C1)}] The \'{e}tale topological groupoid $\Gamma$ is
locally connected.
\item [{(C2)}] Denote by $\alpha,\omega$ the maps sending each
morphism in $\Gamma$ to its right and left units. For any
$U_i,U_j$, set $\Gamma(U_i,U_j)=\{\gamma\in\Gamma\mid\alpha(\gamma)
\in \widehat{U_i},\omega(\gamma)\in \widehat{U_j}\}$. Then the
restriction of $\alpha$ {\em(}resp. $\omega${\em)} to any connected 
component of $\Gamma(U_i,U_j)$ is a homeomorphism onto a connected 
component of the image of $\alpha|_{\Gamma(U_i,U_j)}$ {\em (}resp.
$\omega|_{\Gamma(U_i,U_j)}${\em )} in $\widehat{U_i}$ {\em(}resp. 
$\widehat{U_j}${\em)}.
\end{itemize}
}

We remark that the condition (C2) is not preserved under the usual 
equivalence of \'{e}tale topological groupoids. For example, consider
Thurston's ``teardrop'' orbifold which is a $2$-sphere with
one orbifold point of order $n$ (cf. \cite{Th}). One can 
easily construct a defining \'{e}tale topological groupoid $\Gamma$ 
which does not satisfy (C2) (eg. when the space of units of $\Gamma$,
$\bigsqcup_i\widehat{U_i}$, contains only two connected components). 
On the other hand, we will show that 
smooth orbifolds satisfy (C1) and (C2) after taking an appropriate 
refinement of any given cover of local uniformizing systems 
(cf. Proposition 2.1.3). 

\vspace{3mm}

With the preceding understood, the main results in this paper
are listed below.

\begin{thm}
Consider the set of orbispaces which satisfy {\em (C1), (C2)}.
There is a naturally defined notion of maps under which the said
set of orbispaces forms a category.
\end{thm}

From now on, without explicitly mentioning to the contrary, all 
orbispaces are assumed to satisfy (C1), (C2).

\vspace{3mm}

The following structural theorem of mapping space will also serve
as a technical foundation for the subsequent development. (In this
regard, the key technical lemma is Lemma 3.2.2.)

\begin{thm}
Let $X$ be a paracompact, locally compact and Hausdorff
orbispace\footnote{see the beginning of \S 3.2 for the precise definition
of this condition.}, $X^\prime$ be any orbispace. The set of maps
from $X$ to $X^\prime$ is naturally an orbispace {\em(}as defined in 
Definition 1.1{\em)} under a canonical \'{e}tale topological groupoid.
\end{thm}

Specializing in the case of smooth orbifolds, one may consider the
space of maps of $C^r$ class, namely, those maps which can be
represented locally by $C^r$-maps between local uniformizing
systems. A natural topology can be given to the space of $C^r$-maps 
between smooth orbifolds, which generalizes the Whitney
topology in the case of smooth manifolds (cf. \cite{Hir}). The
following theorem generalizes the corresponding basic results on
smooth manifolds.

\begin{thm}
Let $X$, $X^\prime$ be any smooth orbifolds where $X$ is compact.
\begin{itemize}
\item [{(1)}] The space of $C^r$-maps from $X$ to $X^\prime$ is
naturally a smooth Banach orbifold. In particular, it is Hausdorff
and second countable.
\item [{(2)}] For any $l\geq r$, the set of $C^l$-maps is a dense
subset of the space of $C^r$-maps.
\end{itemize}
\end{thm}

Suppose $X$, $X^\prime$ are complete Riemannian orbifolds and $X$ is
compact. Then each $C^r$-map from $X$ to $X^\prime$ has a natural $C^r$-norm. 
(The corresponding topology is called the $C^r$-topology.)
The following result extends the classical
Arzela-Ascoli theorem to the orbifold setting.

\begin{thm}
Let $X$, $X^\prime$ be any complete Riemannian orbifolds where
$X$ is compact. For any sequence of $C^r$-maps from $X$ to
$X^\prime$ which have bounded $C^r$-norms, there is a subsequence
which converges to a $C^{r-1}$-map in the $C^{r-1}$-topology.
\end{thm}

\begin{re}{\em
(1) A map between orbispaces defined in this paper is a certain equivalence 
class of homomorphisms between the corresponding defining \'{e}tale 
topological groupoids. In this sense,
there are several related notions of maps in the
literature, cf. Haefliger \cite{Ha1}, Hilsum-Skandalis \cite{HS}, and Pronk
\cite{Pr}. (See also Moerdijk \cite{M}). Our definition is most closely
related to that of an Haefliger structure in \cite{Ha1}. More precisely,
let $T$ be a locally connected topological space which is trivially regarded 
as an orbispace, $X$ be an orbispace with a defining  \'{e}tale topological 
groupoid $\Gamma$. Then a map from $T$ to $X$ defined in this paper may be 
canonically identified with a $\Gamma$-structure on $T$ defined in \cite{Ha1}. 

\vspace{1.5mm}

(2) The space of maps as the set of certain equivalence classes of groupoid
homomorphisms is naturally the orbit space of a certain ``tautological''
groupoid. It is not obvious, however, that one can define a nature topology 
on the groupoid which makes it into an \'{e}tale topological groupoid. Our 
main result Theorem 1.3 asserts that this can be done assuming the technical
conditions (C1) and (C2). It would be interesting to know whether these
conditions can be removed. 

\vspace{1.5mm}

(3) It proved to be important and much more convenient
to allow non-effective actions in the local uniformizing systems of an
orbifold. This will be assumed throughout the work unless it is explicitly
mentioned to the contrary. Theorem 1.4 and Theorem 1.5 remain valid for
orbifolds in this more general sense.

\vspace{1.5mm}

(4) The Arzela-Ascoli theorem is fundamental in proofs of various compactness 
theorems. In \cite{CR1} the orbifold version of Gromov's compactness theorem
for pseudoholomorphic curves was proved by an {\it ad hoc} method, relying
essentially on the unique continuity property of pseudoholomorphic curves
(so were many other arguments in \cite{CR1}), thus requiring the involved 
almost complex structures be of $C^\infty$ class. Theorem 1.5 allows us to 
bring the proof of orbifold Gromov Compactness Theorem to the usual line 
of arguments.
}

\hfill $\Box$

\end{re}

{\it A brief history of this work.}\hspace{1mm}
Some of the ideas and results in this work were first written down in an
article under the title ``A homotopy theory of orbispaces'' \cite{C} in
January 2001, which were also presented in the talk \cite{C0} 
(compare also \cite{LU}) in the Madison conference in May 2001. 
The main point of this preliminary version is that
the category of orbispaces introduced in \cite{C} may provide a mathematical
framework for considerations originating from the orbifold string theory 
\cite{DHVW}. More concretely, the following observations were made in \cite{C}.
\begin{itemize}
\item [{(1)}] Let $X=Y/G$ be an orbifold which is a global quotient. Then the
loop space of $X$, ie., the space of maps from $S^1$ into $X$, can be 
identified with the space $P(Y,G)/G$, where
$$
P(Y,G)\equiv\{(\gamma,g)|\gamma:\R\rightarrow Y, g\in G \mbox{ such that }
\gamma(t+2\pi)=g\cdot\gamma(t)\} \leqno (1.7)
$$
and $G$ acts on $P(Y,G)$ by $h\cdot (\gamma,g)=(h\circ\gamma, Ad(h)(g))$
(cf. Lemma 3.5.1 in \cite{C}).
\item [{(2)}] For any orbifold $X$, the space of ``constant'' loops in $X$, 
ie., the fixed-point set of the canonical $S^1$-action on the loop space, 
can be naturally identified with 
$$
\widetilde{X}\equiv \{(p,(g)_{G_p})|p\in X, g\in G_p\} \leqno (1.8)
$$ 
(here $G_p$ stands for the isotropy group at $p$ and $(g)_{G_p}$ the 
conjugacy class of $g$ in $G_p$), which is a disjoint union of orbifolds 
of various dimensions, containing $X$ as an open and closed suborbifold
of top dimension (cf. Proposition 3.5.3 in \cite{C}).
\item [{(3)}] The loop space of an orbifold has a natural infinite 
dimensional orbifold structure (cf. Theorem 3.5.5 in \cite{C}). (This was
also independently observed in \cite{GH}.)
\end{itemize}
Based on these observations, we made the following speculations.
\begin{itemize}
\item [{(1)}] Notice that the space $P(Y,G)$ in $(1.7)$ is exactly the
configuration space of strings satisfying the so-called ``twisted boundary 
conditions'' (cf. $(1.1)$) considered by the physicists in \cite{DHVW}. Thus 
observation (1) suggests that the loop space of an orbifold in \cite{C} 
may serve as the configuration space of strings for string theories of 
orbifolds in general.
\item [{(2)}] The space $\widetilde{X}$ in $(1.8)$ (first introduced by
Kawasaki in \cite{Ka}) played a key role in the orbifold 
Gromov-Witten theory \cite{CR1, CR2}. More concretely, the quantum cohomology 
of an orbifold $X$ is given as an abelian group by the rational cohomology 
of $\widetilde{X}$ (rather than that of $X$) with degrees properly shifted. 
On the other hand, by the stringy interpretation of quantum cohomology in 
Vafa \cite{V}, the
quantum cohomology ring is a deformation of the rational cohomology ring 
of the space of constant loops. Thus observation (2) suggests that the 
Gromov-Witten invariants in \cite{CR1, CR2} are indeed the 
``stringy invariants'' of orbifolds in \cite{DHVW}.  
\item [{(3)}] In Witten \cite{W}, elliptic genus was interpreted as the 
$S^1$-equivariant index of the Dirac operator on the loop space. By 
observation (3), one may attempt to develop a theory of orbifold elliptic 
genera and orbifold elliptic cohomology by extending the constructions 
in \cite{W} to the loop space of an orbifold. (Later we learned that 
orbifold elliptic genera had already been studied, cf. \cite{Liu}.) 
We also remarked in \cite{C} that one may also attempt to extend the work
of Chas and Sullivan on ``string topology'' \cite{CS} to the 
orbifold category, by working with the loop space of orbifold in \cite{C}. 
\end{itemize}

The bulk of the current version of this work was carried out in the academic
year of 2001-2002, during which time the author was visiting the Institute 
for Advanced Study in Princeton. The work was completed in its current form
in 2003. There were several substantial improvements over the preliminary 
version \cite{C}. For example, the structural theorems were not available 
in \cite{C} except for the case of loop space. The homotopy groups
defined in \cite{C} were only part of those defined in \cite{C1}, which are 
equivalent to the homotopy groups of the corresponding classifying spaces. 
The CW-complex theory in \cite{C1} was not yet developed. 

In January 2001, the author had a conversation with Dennis Sullivan about
the subject discussed in \cite{C}. Sullivan pointed out to the author that
the existing methods for groupoids or alike work best in the case of
cofunctors (eg. cohomology theories). In this regard, the direction taken in
the present study is to develop a formalism that works as well in the case of 
functors. 

\vspace{3mm}

The rest of this paper is organized as follows.

Section 2 is concerned with a basic foundation for the category of
orbispaces introduced in this work. In \S 2.1, we present
the set of orbispaces which satisfy (C1), (C2) in terms of local charts.
This formalism is the most natural one in which (C1), (C2) are
expressed. We end \S 2.1 with an elementary proof that smooth
orbifolds satisfy (C1), (C2) after taking an appropriate
refinement of any given cover of local uniformizing systems.
In \S 2.2, we discuss the definition of the maps studied in this
work, and give the proof of Theorem 1.2.

Section 3 is devoted to a basic structural study on mapping spaces.
Theorem 1.3, Theorem 1.4 and Theorem 1.5 are proved in \S 3.2, \S 3.3 
and \S 3.4 respectively. In \S 3.1, three preliminary lemmas are presented.

\vspace{6mm}

\centerline{\bf Acknowledgments}

\vspace{3mm}

I wish to thank Yongbin Ruan for introducing me to the work of Dixon,
Harvey, Vafa and Witten on string theories of orbifolds, and for the
joint work on orbifold Gromov-Witten theory, without which the work in this 
paper and its sequel would not have become possible. I am also very grateful
to Andr\'{e} Haefliger for several very useful email correspondences,
and to Dennis Sullivan for some enlightening conversations. I also wish
to thank the Institute for Advanced Study for her warm hospitality and
excellent academic atmosphere, where the bulk of this research in its 
current form was done.
This work has been partially supported by NSF grant DMS-9971454, and by 
the financial support from the Institute for Advanced Study through NSF 
grant DMS-9729992. During the final preparation, it is partially supported 
by NSF grant DMS-0304956.

\sectioni{A category of orbispaces}
\subsection{Groupoid versus local chart}

We begin by recalling briefly the basic definitions regarding
groupoids. See e.g. Haefliger \cite{Ha3} or Bridson-Haefliger
\cite{BH} for details.

A groupoid $\Gamma$ is a small category whose morphisms are all
invertible. The set of objects of $\Gamma$ is naturally identified
with the set of units $U$. There are mappings $\alpha,\omega:\Gamma
\rightarrow U$ sending each morphism $\gamma\in\Gamma$ to its
initial object (also called right unit) and its terminal object
(also called left unit) respectively. For any $x\in U$, the set
$\Gamma_x=\{\gamma\in\Gamma\mid\alpha(\gamma)=\omega(\gamma)=x\}$ is
naturally a group, called the isotropy group of $x$. The set
$\Gamma\cdot x=\{y\in U\mid\exists\gamma\in\Gamma \mbox{ s.t. }
\alpha(\gamma)=x,\omega(\gamma)=y\}$ is called the $\Gamma$-orbit
of $x$. The set of $\Gamma$-orbits is denoted by $\Gamma \backslash
U$.

A topological groupoid is a groupoid $\Gamma$ where $\Gamma$ is
also a topological space, such that with the induced topology on
the set of units $U$, the mappings $\alpha,\omega:\Gamma\rightarrow
U$ as well as the mappings of taking composition and inverse are
continuous. The set of $\Gamma$-orbits is naturally a topological
space with the quotient topology. A topological groupoid is called
\'{e}tale if $\alpha,\omega:\Gamma\rightarrow U$ are local
homeomorphisms.

Given any pair $(Y,G)$, where $G$ is a topological group acting
continuously on a topological space $Y$ from the left, one may
canonically put a topological groupoid structure on $G\times Y$
as follows. Define composition and inverse by $(h,g\cdot y)\circ
(g,y)=(hg,y)$ and $(g,y)^{-1}=(g^{-1},g\cdot y)$. The space of
units is naturally identified with the space $Y$, and $\alpha,\omega$
are given by $\alpha(g,y)=y$ and $\omega(g,y)=g\cdot y$. The
groupoid $G\times Y$ is \'{e}tale if and only if $G$ is a discrete
group.

A homomorphism of topological groupoids from $\Gamma$ to
$\Gamma^\prime$ is a continuous map $\Phi:\Gamma\rightarrow
\Gamma^\prime$, which commutes with the mappings $\alpha,\omega$ and
the mappings of taking composition and inverse, and which induces a
continuous map between the spaces of units. When both $\Gamma$
and $\Gamma^\prime$ are \'{e}tale, $\Phi$ is called an equivalence
if it is a local homeomorphism and induces a bijection between
$\Gamma\backslash U$ and $\Gamma^\prime\backslash U^\prime$ and an
isomorphism $\Gamma_x\cong\Gamma^\prime_{\Phi(x)}$ for any $x\in U$.

With the preceding understood, we return to the conditions (C1), (C2)
in the introduction. First, by (C1) the space of units
$U=\bigsqcup_i \widehat{U_i}$ is locally connected, which implies
that the topological space $X=\Gamma\backslash U$ is also locally
connected. Second, the \'{e}tale topological groupoid $\Gamma$ is
only determined up to equivalence. Because the condition (C2) may not
be preserved under equivalence of \'{e}tale topological groupoids, it
should be understood throughout this work that when an orbispace is said
to satisfy (C1), (C2), it is meant that some preferred \'{e}tale
topological groupoid in the equivalence class has been chosen, with
respect to which (C1), (C2) are satisfied.

Orbispaces satisfying (C1), (C2) may be equivalently formulated
in a way in terms of local charts instead of groupoids, which is more
concrete and geometric. In fact, it is in this formalism that conditions
(C1), (C2) are most naturally expressed and sufficiently exploited.
We give details of this formulation in the following proposition.

\begin{prop}
Let $X$ be an orbispace under an \'{e}tale topological groupoid
$\Gamma$ which satisfies {\em (C1), (C2)}. Then $\Gamma$ may be
given by the following set of data:
\begin{itemize}
\item An atlas of local charts
$\{(\widehat{U_i},G_{U_i},\pi_{U_i})\}$, where each $U_i$ is a
connected open subset of $X$, $\widehat{U_i}$ is a locally
connected space with continuous left action by a discrete group
$G_{U_i}$, and $\pi_{U_i}:\widehat{U_i}\rightarrow X$ is a
continuous map inducing a homeomorphism $\widehat{U_i}/G_{U_i}
\cong U_i$. Note that the atlas of local charts
$\{(\widehat{U_i},G_{U_i},\pi_{U_i})\}$ is naturally identified
with the set $\U=\{U_i\}$, which forms an open cover of $X$.
\item A collection of discrete sets
$$
\T=\{T(U_i,U_j)\mid U_i,U_j\in\U, \mbox{ s.t. } U_i\cap
U_j\neq\emptyset\}
$$
which satisfies the following conditions:
\begin{itemize}
\item $T(U_i,U_i)=G_{U_i}$, and for any $i\neq j$, $T(U_i,U_j)=
\bigsqcup_{u\in I_{ij}} T_{W_u}(U_i,U_j)$ where $\{W_u\mid u\in I_{ij}\}$
is the set of connected components of $U_i\cap U_j$.
\item Each $\xi\in T_{W_u}(U_i,U_j)$ is assigned with a
homeomorphism $\phi_\xi$, whose domain and range are connected
components of the inverse image of $W_u$ in $\widehat{U_i}$ and
$\widehat{U_j}$ respectively, such that $\pi_{U_i}=\pi_{U_j}\circ\phi_\xi$.
(For each $\xi\in G_{U_i}$, we define $\phi_\xi$ to be the
self-homeomorphism of $\widehat{U_i}$ induced by the action of $\xi$.)
\item For any $\xi\in T(U_i,U_j)$, $\eta\in T(U_j,U_k)$, and any $x\in
\mbox{Domain }(\phi_\xi)$ such that $\phi_\xi(x)\in\mbox{Domain }
(\phi_\eta)$, there exists an $\eta\circ\xi(x)\in T(U_i,U_k)$ such
that $\phi_{\eta\circ\xi(x)}(x)=\phi_\eta(\phi_\xi(x))$. Moreover,
$x\mapsto \eta\circ\xi(x)$ is locally constant, and the composition
$(\xi,\eta)\mapsto \eta\circ\xi(x)$ is associative and coincides with
the group multiplication in $G_{U_i}$ when restricted to $T(U_i,U_i)=G_{U_i}$.
\item Every $\xi\in T(U_i,U_j)$ has an inverse $\xi^{-1}\in T(U_j,U_i)$,
such that $\mbox{Domain }(\phi_\xi)=\mbox{Range }(\phi_{\xi^{-1}})$,
$\mbox{Domain }(\phi_{\xi^{-1}})=\mbox{Range }(\phi_\xi)$, and
$\xi^{-1}\circ\xi(x)=1$ $\forall x\in \mbox{Domain }(\phi_\xi)$,
$\xi\circ\xi^{-1}(x)=1$ $\forall x\in\mbox{Domain }(\phi_{\xi^{-1}})$.
\end{itemize}
\end{itemize}
\end{prop}

\pf Assume that $\Gamma$ is given, satisfying (C1), (C2),
and that the space of units of $\Gamma$ is the disjoint union
$\bigsqcup_i\widehat{U_i}$ with the restriction to each $\widehat{U_i}$
being the product groupoid $G_i\times\widehat{U_i}$. We obtain the
atlas of local charts by setting $G_{U_i}=G_i$, and $\pi_{U_i}:\widehat{U_i}
\rightarrow X$ the orbit map. (Without loss of generality, we may
always assume that $U_i=\widehat{U_i}/G_i$ is connected.)

As for the collection of discrete sets $\T=\{T(U_i,U_j)\}$, we define
$T(U_i,U_i)=G_{U_i}$, and for any $i\neq j$, we define $T(U_i,U_j)$ to
be the set of connected components of $\Gamma(U_i,U_j)=\{\gamma\in\Gamma
\mid \alpha(\gamma)\in \widehat{U_i},\omega(\gamma)\in \widehat{U_j}\}$.
For any $\xi\in T(U_i,U_j)$, $i\neq j$, by (C2), there exists an inverse
of $\alpha$ defined from a connected component of the image of
$\alpha|_{\Gamma(U_i,U_j)}$ in $\widehat{U_i}$ to the connected
component of $\Gamma(U_i,U_j)$ which is named by $\xi$. We define
$\phi_\xi=\omega\circ s_\xi$ where $s_\xi$ is the said inverse of $\alpha$.
We define $\eta\circ\xi(x)$ to be the connected component which
contains the morphism $s_\eta(\phi_\xi(x))\circ s_\xi(x)$, and
define the inverse $\xi^{-1}$ to be the connected component which
contains the morphism $s_\xi(x)^{-1}$ for some (and hence for all)
$x\in\mbox{Domain }(\phi_\xi)$. It is easy to verify that these
objects satisfy the claimed conditions. We leave the details to
the reader.

Conversely, assume that the said set of data is given, we may recover
the groupoid $\Gamma$ as follows. Define $\Gamma$ to be the subset
of $\bigsqcup_{i,j} T(U_i,U_j)\times\widehat{U_i}$ with induced
topology which consists of pairs $(\xi,x)$ such that $x\in\mbox{Domain
}(\phi_\xi)$. The space of units is the disjoint union $\bigsqcup_i
\widehat{U_i}$ via the identification $x\mapsto (1,x)$. The groupoid
structure of $\Gamma$ is given as follows: $\alpha(\xi,x)=x$,
$\omega(\xi,x)=\phi_\xi(x)$, $(\eta,y)\circ (\xi,x)=(\eta\circ\xi(x),x)$
where $y=\phi_\xi(x)$, and $(\xi,x)^{-1}=(\xi^{-1},\phi_\xi(x))$. It is
an easy exercise to verify that $\Gamma$ is an \'{e}tale topological
groupoid which satisfies (C1), (C2). We leave the details to
the reader.

\hfill $\Box$

Some remarks are in order to further explain a few points.

\begin{re}{\em

\vspace{1.5mm}

(1) One may always require that each local chart $\widehat{U_i}$
is also connected by replacing it with one of its connected component,
which amounts to replace $\Gamma$ with an equivalent groupoid such
that (C1), (C2) are still satisfied.

\vspace{1.5mm}

(2) Since $x\mapsto \eta\circ\xi(x)$ is locally constant, $\eta\circ\xi(x)$
depends only on the connected component of
$\phi_\xi^{-1}(\mbox{Domain }(\phi_\eta))$ in which $x$ lies. More
generally, for any sequence $\xi_{k,k+1}\in T(U_k,U_{k+1})$, $1\leq k\leq
n-1$, such that
$$
\phi_{\xi_{n-1,n}}\circ\cdots\circ\phi_{\xi_{23}}\circ\phi_{\xi_{12}}
\leqno (2.1.1)
$$
is defined, then $\xi_{n-1,n}\circ\cdots\circ\xi_{23}\circ\xi_{12}(x)$
depends only on the connected component of the domain of $(2.1.1)$
in which $x$ lies. We introduce
$$
\Lambda(\xi_{12},\xi_{23},\cdots,\xi_{n-1,n})=\{ {\bf a}|{\bf
a}\mbox{ is a connected component of the domain of } (2.1.1)\},
\leqno (2.1.2)
$$
and for any ${\bf a}\in\Lambda(\xi_{12},\xi_{23},\cdots,\xi_{n-1,n})$,
set
$$
\xi_{n-1,n}\circ\cdots\circ\xi_{23}\circ\xi_{12}({\bf a})
=\xi_{n-1,n}\circ\cdots\circ\xi_{23}\circ\xi_{12}(x), \forall x\in
{\bf a}.
$$
When there is no ambiguity, we simply write
$\xi_{n-1,n}\circ\cdots\circ\xi_{23}\circ\xi_{12}$ for
$\xi_{n-1,n}\circ\cdots\circ\xi_{23}\circ\xi_{12}({\bf a})$.

\vspace{1.5mm}

(3) For any $\xi\in T(U_i,U_j)$, there exists an isomorphism of groups
$\lambda_\xi$ defined by
$$
\lambda_\xi(g)=\xi\circ g\circ\xi^{-1}(x), \leqno (2.1.3)
$$
where $x$ is any point in $\mbox{Domain }(\phi_{\xi^{-1}})$ and $g$ is in
the subgroup of $G_{U_i}$ fixing the set $\mbox{Domain }(\phi_\xi)$. The
range of $\lambda_\xi$ is the subgroup of $G_{U_j}$ fixing the set
$\mbox{Range }(\phi_\xi)$. In particular, for any $\xi\in T(U_i,U_i)=G_{U_i}$,
$\lambda_\xi=Ad(\xi):G_{U_i}\rightarrow G_{U_i}$. We observe that for any
$\xi$, $\phi_\xi$ is $\lambda_\xi$-equivariant, and the pairs
$(\phi_\xi,\lambda_\xi)$ satisfy
$$
(\phi_{\eta\circ\xi({\bf a})},\lambda_{\eta\circ\xi({\bf a})})=
(\phi_\eta,\lambda_\eta)\circ (\phi_\xi,\lambda_\xi)
$$
when restricted to $({\bf a},G_{\bf a})$ for any ${\bf
a}\in\Lambda(\xi,\eta)$, where $G_{\bf a}$ is the subgroup
of $\mbox{Domain }(\lambda_\xi)$ which fixes ${\bf a}$. Hence each
element $\xi$ is associated with a `transformation of local charts'
$(\phi_\xi,\lambda_\xi)$. Roughly speaking, (C2) dictates that
each `transformation of local charts' be defined over a domain that is
maximally large as possible.  In a certain sense, \'{e}tale topological
groupoids satisfying (C1), (C2) are the closest generalization of
the product groupoids, or in other words, orbispaces satisfying (C1),
(C2) are the closest generalization of the orbispaces of global quotients.

\vspace{1.5mm}

(4) We may maximize the atlas of local charts by adding all the
connected open subsets of each $U_i\in\U$ to $\U=\{U_i\}$.  This
amounts to change $\Gamma$ to an equivalent groupoid (still
satisfying (C1), (C2)), which is done as follows. Let $V$ be any
connected open subset of some $U_{i_0}\in\U$, and $\widehat{V}$ be
the inverse image of $V$ in $\widehat{U_{i_0}}$. We define a
groupoid $\Gamma_1=\Gamma\bigsqcup\Gamma|_{\widehat{V}}\bigsqcup
\Gamma^\prime$ where $\Gamma^\prime=\{\gamma\in\Gamma\mid
\mbox{ either } \alpha(\gamma)\in\widehat{V}\mbox{ or }
\omega(\gamma)\in\widehat{V}\}$. The space of units of
$\Gamma_1$ is the disjoint union
$\widehat{V}\bigsqcup(\bigsqcup_i\widehat{U_i})$, hence the atlas
of local charts for $\Gamma_1$ is $\{V\}\bigcup\U$. Moreover, the
inclusion $\Gamma\subset\Gamma_1$ induces an equivalence of
groupoids, and $\Gamma_1$ satisfies (C1), (C2). One may iterate
this process to maximize $\U$ by Zorn's Lemma.

\vspace{1.5mm}

(5) We introduce the following equivalence relation: two orbispace
structures on $X$ is said to be directly equivalent if the atlas
of local charts of one of them is contained in that of the other,
and is said to be equivalent if they are related by a finite chain
of directly equivalent orbispace structures.

\vspace{1.5mm}

(6) Product groupoids $G\times Y$, where $Y$ is locally connected,
satisfy (C1), (C2) trivially. Thus the global quotient spaces are the
most basic examples of orbispaces considered in this work. Other
known examples of orbispaces satisfying (C1), (C2) include smooth
orbifolds (which is shown next), and orbihedra of Haefliger in
\cite{Ha3}. On the other hand, this class of orbispaces is closed
under taking subspaces in the following sense. Let $Y$ be a subset
of $X$ such that the inverse image of $Y$ in
$\bigsqcup_i\widehat{U_i}$, $\widehat{Y}$, is locally connected.
Then $Y$ becomes an orbispace under the groupoid $\Gamma|_{\widehat{Y}}$,
which also satisfies (C1), (C2).

\vspace{1.5mm}

(7) We shall explain here that, with the technical assumptions (C1) and (C2),
each orbispace may be associated with a direct limit of complexes of groups 
defined in Haefliger \cite{Ha3}. We hope that this direct limit would play
the role of `structure group' (ie. the role of 
$G$ in a global quotient $X=Y/G$) for a general orbispace. 

We first recall the definition of complex of groups. Let $K$ be a simplicial 
cell complex. We set $V(K)$ for the set of barycenters of cells of $K$, and 
$E(K)$ for the set of edges of the barycentric subdivision of $K$. Each edge
$a\in E(K)$ is natually oriented, ie., if the initial point $i(a)$ of $a$ 
is the barycenter of a cell $\sigma$ and the terminal point $t(a)$ is the
barycenter of $\tau$, then $\dim\tau<\dim\sigma$. Two edges $a,b\in
E(K)$ are said to be composable if $i(a)=t(b)$ and the composition 
$c=ab$ is the edge with $i(c)=i(b)$ and $t(c)=t(a)$ such that $a, b$ and $c$ 
form the boundary of a $2$-simplex in the barycentric subdivision of $K$. 

With the above notations understood, a complex of groups 
$G(K)=(K,G_\sigma,\psi_a,g_{a,b})$ on $K$ is given by the following data:
\begin{itemize}
\item [{(1)}] a group $G_\sigma$ for each cell $\sigma\in V(K)$,
\item [{(2)}] an injective homomorphism $\psi_a:G_{i(a)}\rightarrow G_{t(a)}$
for each edge $a\in E(K)$,
\item [{(3)}] an element $g_{a,b}\in G_{t(a)}$ for each pair $a,b\in E(K)$ of 
composable edges such that
$$
Ad(g_{a,b})\circ\psi_{ab}=\psi_a\circ\psi_b, \leqno (2.1.4a)
$$
and the set of elements $\{g_{a,b}\}$ satisfies the cocycle condition
$$
\psi_a(g_{b,c})g_{a,bc}=g_{a,b}g_{ab,c} \leqno (2.1.4b)
$$
for any triple $a,b,c\in E(K)$ of composable edges.
\end{itemize}
 Homomorphisms of complexes of groups are defined as follows.
Let $G(K)=(K,G_\sigma,\psi_a,g_{a,b})$, 
$G(K^\prime)=(K^\prime,G^\prime_{\sigma^\prime},\psi^\prime_{a^\prime},
g^\prime_{a^\prime,b^\prime})$ be complexes of groups on $K, K^\prime$
respectively, and let $f:K\rightarrow K^\prime$ be a simplicial map. 
A homomorphism $\phi=(f,\phi_\sigma,g_a^\prime): 
G(K)\rightarrow G(K^\prime)$ over $f$ is given by the following data:
\begin{itemize}
\item [{(1)}] a homomorphism $\phi_\sigma: G_\sigma\rightarrow 
G^\prime_{f(\sigma)}$ for each cell $\sigma\in V(K)$,
\item [{(2)}] an element $g_{a}^\prime\in G^\prime_{f(t(a))}$ for each edge 
$a\in E(K)$ such that
$$
Ad(g^\prime_a)\circ\psi^\prime_{f(a)}\circ\phi_{i(a)}=\phi_{t(a)}\circ\psi_a, 
\leqno (2.1.5a)
$$
and for each pair $a,b\in E(K)$ of composable edges, 
$$
\phi_{t(a)}(g_{a,b})g^\prime_{ab}=g^\prime_a\psi^\prime_{f(a)}(g^\prime_b)
g^\prime_{f(a),f(b)}.   \leqno (2.1.5b)
$$
(Note: when $f(a)$ is not a well-defined edge in $E(K^\prime)$, ie. when
$f(i(a))=f(t(a))$, we set $\psi^\prime_{f(a)}=id:G^\prime_{f(i(a))}\rightarrow
G^\prime_{f(t(a))}$, and moreover, we set $g^\prime_{f(a),f(b)}=1$ if either 
$f(a)$ or $f(b)$ is not well-defined.)
\end{itemize}
Finally, two complexes of groups $(K,G_\sigma,\psi_a,g_{a,b})$, 
$(K,G_{\sigma},\psi^\prime_{a}, g^\prime_{a,b})$ on $K$ are said to
differ by a coboundary if there is a homomorphism $\phi=(id,id,g_a)$
over $id:K\rightarrow K$ between them.

Now let $X$ be an orbispace with atlas of local charts $\U$. For any cover
$\{U_i\}$ of $X$ where $U_i\in \U$, we define a complex of groups $G(\{U_i\})$
as follows. The corresponding simplicial cell complex $K(\{U_i\})$ consists
of $n$-cells $(\{U_{i_0},\cdots, U_{i_n}\},j)$, $n\geq 0$, where 
$U_{i_0},\cdots, U_{i_n}$ are distinct elements in $\{U_i\}$ such that
$U_{i_0}\cap\cdots\cap U_{i_n}\neq\emptyset$, and $j$ is a connected 
component of $U_{i_0}\cap\cdots\cap U_{i_n}$. The faces of 
$(\{U_{i_0},\cdots, U_{i_n}\},j)$ are obtained by removing one of $U_{i_k}$'s
and by replacing $j$ with the corresponding connected component containing $j$.

To obtain the remaining data, we rely on the following fact due to the
assumptions (C1) and (C2): For each cell 
$\sigma=(\{U_{i_0},\cdots, U_{i_n}\},j)$, $j$ is a local chart of $X$, which
we denote by $U_\sigma$. With this understood, we define the group $G_\sigma$ 
in  $G(\{U_i\})$ associated to $\sigma$ by $G_\sigma\equiv G_{U_\sigma}$. 
Now observe that an edge $a$ with $i(a)=\sigma$ corresponds to removing one 
or several $U_{i_k}$'s in $\{U_{i_0},\cdots, U_{i_n}\}$, and if we
let $j_a$ be the corresponding connected component for the cell $t(a)$,
we have $j\subset j_a$. We define $\psi_a$ in $G(\{U_i\})$ by fixing a 
$\xi\in T(j,j_a)$ and set $\psi_a\equiv \lambda_\xi$ (see (2.1.3) for the 
definition of $\lambda_\xi$). Then it follows easily that there are unique 
elements $g_{a,b}$ such that $(2.1.4a)$ and $(2.1.4b)$ are satisfied, thus 
giving the construction of $G(\{U_i\})$. If we choose a different 
$\xi\in T(j,j_a)$, the corresponding complex of groups differs by a 
coboundary. 

We leave to the reader to verify that for any refinement $\{V_k\}$ of
$\{U_i\}$, there is a canonically defined homomorphism from $G(\{V_k\})$
to $G(\{U_i\})$, so that in this way, the orbispace $X$ is associated with
a direct limit of complexes of groups. Moreover, a map between two orbispaces 
(to be defined in the next subsection) induces a ``morphism'' between the 
corresponding direct limits of complexes of 
groups. Finally, we observe that if $X=Y/G$ and $X$ is connected, and 
$\{U_i\}=\{X\}$, then the simplicial cell complex $K(\{U_i\})$ consists of 
a single vertex, and the complex of groups $G(\{U_i\})$ reduces to the group 
$G$.
}
\end{re}

\hfill $\Box$

The remainder of this subsection is occupied with an elementary
proof that after taking a suitable refinement of any given cover of
local uniformizing systems, smooth orbifolds satisfy (C1), (C2).
We remark that passing to a refinement is necessary here,
as one can easily construct an example not satisfying (C2)
using Thurston's `teardrop' orbifold (cf. \cite{Th}).

Recall the definition of orbifolds due to Satake \cite{Sa}.
Suppose $X$ is a Hausdorff, paracompact space. An orbifold
structure on $X$ is given by an open cover $\{U_i\}$ satisfying the
following conditions: (1) Each $U_i$ is given with a uniformizing system
$(\widehat{U_i},G_i,\pi_i)$, where $\widehat{U_i}$ is a connected open
subset of $\R^n$, $G_i$ is a finite group of self-diffeomorphisms
of $\widehat{U_i}$, and $\pi_i:\widehat{U_i}\rightarrow U_i$ is a
continuous map inducing a homeomorphism $\widehat{U_i}/G_i\cong U_i$.
(2) For any pair $U_i,U_j$ such that $U_j\subset U_i$, there is an
associated set of injections $\{\phi\}$, where
$\phi:\widehat{U_j}\rightarrow \widehat{U_i}$ is an
open embedding such that $\pi_i\circ\phi=\pi_j$.
The group $G_i$ acts on the set of injections freely and transitively
via post-compositions. Moreover, for any $U_k\subset U_j\subset U_i$,
the composition of an injection associated to $U_k\subset U_j$
with an injection associated to $U_j\subset U_i$ is an injection
associated to $U_k\subset U_i$. (3) For any point $p\in X$ such that
$p\in U_i\cap U_j$, there is a $U_k$ satisfying $p\in U_k$ and
$U_k\subset U_i\cap U_j$. In Satake's original definition, it was
assumed that for each $i$ the fixed-point set of any element of $G_i$
in $\widehat{U_i}$ is of codimension at least $2$. This assumption
implies the properties about injections listed in Condition (2) above.
Later in \cite{Ka}, Kawasaki removed Satake's original assumption,
and instead, he imposed these properties as part of the axioms in the
definition of orbifolds.

The orbifolds of Satake may also be described equivalently using
\'{e}tale topological groupoids. More concretely, let $U=\bigsqcup_i
\widehat{U_i}$ be the disjoint union. Then the set of open embeddings
$\{\phi\}$ generates a pseudogroup of local diffeomorphisms of $U$.
The associated \'{e}tale topological groupoid $\Gamma$, which is the
space of germs given with the \'{e}tale topology, defines an orbispace
structure on $X$ in the sense of Definition 1.1.
See e.g. \cite{BH} for more details.

\begin{prop}
Given any orbifold $(X,\{U_i\})$ as defined by Satake, with the
associated \'{e}tale topological groupoid denoted by $\Gamma$,
there is a refinement $\{V_\alpha\}$ of $\{U_i\}$, such that the
restriction of $\Gamma$ to the inverse image of $\bigsqcup_\alpha
V_\alpha$ in $\bigsqcup_i\widehat{U_i}$ is an equivalent groupoid
with respect to which {\em (C1), (C2)} are satisfied.
\end{prop}

\pf
We need to assume an auxiliary Riemannian metric on the
orbifold $X$. Recall that a Riemannian metric on $X$ is a collection
of metrics $\{g_i\}$ where $g_i$ is a $G_i$-equivariant Riemannian
metric on $\widehat{U_i}$ with respect to which each injection is
isometric. Such a metric on $X$ exists because $X$ is paracompact.
Given a Riemannian metric on $X$, a path $\gamma:I\rightarrow X$ is
called a geodesic if for any $t\in I$, the restriction of
$\gamma$ in a neighborhood of $t$ can be lifted to a geodesic
$\tilde{\gamma}_t$ in $\widehat{U_i}$ for some $i$. We observe two
facts about geodesics: (1) the length of a geodesic
$\gamma:I\rightarrow X$ is well-defined, which is
$\ell(\gamma)=\int_I |\frac{d\tilde{\gamma}_t}{dt}(t)|dt$, and
(2) if $\gamma(I)\subset U_i$ for some $i$, then there is a geodesic
$\tilde{\gamma}:I\rightarrow \widehat{U_i}$ such that
$\gamma=\pi_i\circ \tilde{\gamma}$, and moreover, $\ell(\gamma)=
\ell(\tilde{\gamma})=\int_I|\frac{d\tilde{\gamma}}{dt}|dt$.

Given any point $p\in X$, suppose $p\in U_i$ for some $i$ and
let $\hat{p}\in\pi_i^{-1}(p)\subset \widehat{U}_i$. There is
a family of neighborhoods $U_p(r)$ of $p$ in $X$,
parameterized by $r$ where $0<r\leq r_p$ for some $r_p$ depending only on
$p$ and the fixed Riemannian metric. Each $U_p(r)$ has a natural
uniformizing system $(B_p(r),G_p,\pi_p)$, where $B_p(r)$ is the open
ball of radius $r$ in $(T\widehat{U_i})_{\hat{p}}$ centered at the
origin, $G_p$ is the stabilizer of $G_i$ at $\hat{p}$, and
$\pi_p:(T\widehat{U_i})_{\hat{p}}\rightarrow X$ is the map sending
each line $tv$, where $v\in (T\widehat{U_i})_{\hat{p}}$,
to a geodesic $\gamma_v$ in $X$ such that
$\gamma_v(0)=p$ and locally at $0$, $\gamma_v$ can be lifted to a
geodesic $\tilde{\gamma}_v$ in $\widehat{U_i}$ with
$\tilde{\gamma}_v(0)=\hat{p}$ and $\frac{d\tilde{\gamma}_v}{dt}(0)=v$.
These neighborhoods $U_p(r)$ have the following nice property:
for any $U_i\in \{U_i\}$, if $U_p(r)\subset U_i$, then the
uniformizing system $(B_p(r),G_p,\pi_p)$ is isomorphic to an induced
one of $U_p(r)$ from $(\widehat{U_i},G_i,\pi_i)$. Each $B_p(r)$
inherits a $G_p$-equivariant metric from the auxiliary metric on $X$.
We may further assume that all $\overline{B_p(r)}$, $r\leq r_p$, are
geodesically convex, meaning that any two points in $\overline{B_p(r)}$
are joined by a unique geodesic in $\overline{B_p(r)}$ such that the distance
from a point on the geodesic to the center of $\overline{B_p(r)}$ is a
strict convex function.

Let $\W=\{U_p(r)|p\in X, r\leq\frac{r_p}{5}\}$. We will show next that
for any $U_{p_1}(r_1),\cdots, U_{p_n}(r_n)\in\W$ with $r_{p_1}\leq r_{p_k}$,
$2\leq k\leq n$, and $U_{p_1}(r_1)\cap\cdots\cap U_{p_n}(r_n)\neq\emptyset$,
the subset $\pi_{p_1}^{-1}(\overline{U_{p_1}(r_1)}\cap\cdots\cap
\overline{U_{p_n}(r_n)})$ of $\overline{B_{p_1}(r_1)}$ has the property
that every two points in $\pi_{p_1}^{-1}(\overline{U_{p_1}(r_1)}\cap\cdots
\cap\overline{U_{p_n}(r_n)})$ are joined by a unique geodesic in
$\pi_{p_1}^{-1}(\overline{U_{p_1}(r_1)}\cap\cdots\cap\overline{U_{p_n}(r_n)})$.
The proof goes as follows. Suppose $x_1,x_2$ are any two points in
$\pi_{p_1}^{-1}
(\overline{U_{p_1}(r_1)}\cap\cdots\cap\overline{U_{p_n}(r_n)})$. Let
$\tilde{\gamma}_1$ be the unique geodesic in $\overline{B_{p_1}(r_1)}$
connecting $x_1,x_2$. We shall prove that $\tilde{\gamma}_1$ lies in the
subset $\pi_{p_1}^{-1}(\overline{U_{p_1}(r_1)}\cap\cdots\cap
\overline{U_{p_n}(r_n)})$. It suffices to show that $\gamma=\pi_{p_1}
(\tilde{\gamma}_1)$ lies in $\overline{U_{p_k}(r_k)}$ for any $2\leq k\leq n$.
To this end, recall that we have assumed that $r_{p_1}\leq r_{p_k}$ for any
$2\leq k\leq n$. It follows that $\ell(\gamma)\leq 2r_1\leq\frac{2r_{p_1}}{5}
\leq \frac{2r_{p_k}}{5}$, so that the distance between $p_k$ and any point on
$\gamma$ is less than or equal to
$r_k+\frac{2r_{p_k}}{5}\leq\frac{3r_{p_k}}{5}$,
which implies that $\gamma$ lies in $U_{p_k}(r_{p_k})$.
Hence there is a geodesic
lifting $\tilde{\gamma}_k$ of $\gamma$ into $B_{p_k}(r_{p_k})$,
with the end points
of $\tilde{\gamma}_k$ lying in $\overline{B_{p_k}(r_k)}$. Now recall that
all $\overline{B_{p_k}(r)}$, $r\leq r_{p_k}$, are geodesically convex. Hence
$\tilde{\gamma}_k$ must lie in $\overline{B_{p_k}(r_k)}$, and therefore
$\gamma=\pi_{p_k}(\tilde{\gamma}_k)$ must lie in $\overline{U_{p_k}(r_k)}$.

Note that each $\pi_{p_1}^{-1}({U_{p_1}(r_1)}\cap\cdots\cap {U_{p_n}(r_n)})$
is a convex, topological ball. In particular, ${U_{p_1}(r_1)}\cap\cdots
\cap {U_{p_n}(r_n)}$ is connected.

Now we take a refinement of $\{U_i\}$, denoted by $\V=\{V_\alpha\}$, where
each $V_\alpha\in\W$. Let $V=\bigsqcup_{U_p(r)\in\V} B_p(r)$ be the disjoint
union. Then the \'{e}tale topological groupoid $\Gamma$ induces an
equivalent \'{e}tale topological groupoid $\Gamma^\prime$ having $V$
as the space of units. ($\Gamma^\prime$ is simply the restriction of
$\Gamma$ by thinking $V$ as a subset of $\bigsqcup_i\widehat{U_i}$,
the space of units of $\Gamma$.) We shall prove that $\Gamma^\prime$
satisfies (C1), (C2).

First, (C1) is trivial. To verify (C2), assume $U_{p_1}(r_1),U_{p_2}(r_2)\in
\V$ have nonempty intersection, and without loss of generality, assume
$r_{p_1}\leq r_{p_2}$. Set $\Gamma(U_{p_1}(r_1),
U_{p_2}(r_2))=\{\gamma\in\Gamma^\prime\mid\alpha(\gamma)\in B_{p_1}(r_1),
\omega(\gamma)\in B_{p_2}(r_2)\}$. Then the image of
$\alpha|_{\Gamma(U_{p_1}(r_1),U_{p_2}(r_2))}$ in $B_{p_1}(r_1)$ is
$\pi_{p_1}^{-1}(U_{p_1}(r_1)\cap U_{p_2}(r_2))$, which is a convex ball
by the nature of the set $\W$. Since $\alpha|_{\Gamma(U_{p_1}(r_1),
U_{p_2}(r_2))}$ is a covering map, it must be a homeomorphism onto
its image when restricted to each connected component of
$\Gamma(U_{p_1}(r_1),U_{p_2}(r_2))$. In particular, each connected component
of $\Gamma(U_{p_1}(r_1),U_{p_2}(r_2))$ is homeomorphic to an open ball.
On the other hand, $\omega|_{\Gamma(U_{p_1}(r_1),U_{p_2}(r_2))}$ is a covering
map onto $\pi_{p_2}^{-1}(U_{p_1}(r_1)\cap U_{p_2}(r_2))$. The group of deck
transformations for any connected component is isomorphic to a subgroup $K$
of $H$, where $H$ is the subgroup of $G_{p_1}$ which leaves the subset
$\pi_{p_1}^{-1}(U_{p_1}(r_1)\cap U_{p_2}(r_2))$ of $B_{p_1}(r_1)$ invariant.
Furthermore, $K$ acts on $\pi_{p_1}^{-1}(U_{p_1}(r_1)\cap U_{p_2}(r_2))$
freely. But this is impossible unless $K$ is trivial, because the action
of $K$ is easily seen to extend over to the closure of
$\pi_{p_1}^{-1}(U_{p_1}(r_1)\cap U_{p_2}(r_2))$, which is a closed ball.
By Brouwer's fixed point theorem, it can not be free unless $K$ is trivial.
Thus $\omega|_{\Gamma(U_{p_1}(r_1),U_{p_2}(r_2))}$ must be a homeomorphism
onto its image also when restricted to each connected component of
$\Gamma(U_{p_1}(r_1),U_{p_2}(r_2))$. This verifies (C2).

\hfill $\Box$

\subsection{Equivalence class of groupoid homomorphisms}

Let $X$ be an orbispace, defined by an \'{e}tale topological groupoid
$\Gamma$. Denote by $\U=\{U_i\}$ the atlas of local charts on $X$.
Then any subset $\{U_\alpha\}$ of $\U$ is associated with an \'{e}tale
topological groupoid $\Gamma\{U_\alpha\}$, which is the restriction of
$\Gamma$ to $\bigsqcup_\alpha\widehat{U_\alpha}$. More explicitly,
$\Gamma\{U_\alpha\}$ is the subset of
$\bigsqcup_{\alpha,\beta} T(U_\alpha,U_\beta)\times \widehat{U_\alpha}$
which consists of pairs $(\xi,x)$ such that $x\in\mbox{Domain }(\phi_\xi)$.

\begin{lem}
Let $X$, $X^\prime$ be orbispaces whose atlas of local charts are
denoted by $\U=\{U_i\}$, $\U^\prime=\{U_{i^\prime}^\prime\}$
respectively, and $\{U_\alpha\}\subset\U$,
$\{U^\prime_{\alpha^\prime}\}\subset\U^\prime$ be any subsets
where each $\widehat{U_\alpha}$ is connected. Then any homomorphism
from $\Gamma\{U_\alpha\}$ to $\Gamma\{U_{\alpha^\prime}^\prime\}$
may be written as $(\{f_\alpha\},\{\rho_{\beta\alpha}\})$, where
there is a mapping $U_\alpha\mapsto U_\alpha^\prime\in
\{U^\prime_{\alpha^\prime}\}$, with respect to which
$f_\alpha:\widehat{U_\alpha}\rightarrow\widehat{U_\alpha^\prime}$,
$\rho_{\beta\alpha}:T(U_\alpha,U_\beta)\rightarrow
T(U_\alpha^\prime,U_\beta^\prime)$, such that
\begin{itemize}
\item [{(a)}] $\phi_{\rho_{\beta\alpha}(\xi)}\circ f_\alpha(x)
=f_\beta\circ\phi_\xi(x)$ for any $\alpha,\beta$, where $\xi\in
T(U_\alpha,U_\beta)$, $x\in\mbox{Domain }(\phi_\xi)$.
\item [{(b)}] $\rho_{\gamma\alpha}(\eta\circ\xi(x))=\rho_{\gamma\beta}(\eta)
\circ\rho_{\beta\alpha}(\xi)(f_\alpha(x))$ for any
$\alpha,\beta,\gamma$, where $\xi\in T(U_\alpha,U_\beta)$,
$\eta\in T(U_\beta,U_\gamma)$, and $x\in\phi_\xi^{-1}(\mbox{Domain }
(\phi_\eta))$.
Note that $f_\alpha(x)\in\phi_{\rho_{\beta\alpha}(\xi)}^{-1}(\mbox{Domain }
(\phi_{\rho_{\gamma\beta}(\eta)})$, which follows from the
equations in {\em (a)} and the assumption that
$x\in\phi_\xi^{-1}(\mbox{Domain }(\phi_\eta))$.
\end{itemize}
\end{lem}

\pf
Any homomorphism from $\Gamma\{U_\alpha\}$ to
$\Gamma\{U_{\alpha^\prime}^\prime\}$ will induce a continuous map
between the corresponding spaces of units, $\bigsqcup_\alpha
\widehat{U_\alpha}$ and $\bigsqcup_{\alpha^\prime}\widehat
{U_{\alpha^\prime}^\prime}$. Since each $\widehat{U_\alpha}$ is
connected, this map determines a mapping $U_\alpha\mapsto
U_\alpha^\prime\in\{U^\prime_{\alpha^\prime}\}$ by the rule that
$\widehat{U_\alpha}$ is being mapped into $\widehat{U_\alpha^\prime}$.
We define $f_\alpha:\widehat{U_\alpha}\rightarrow\widehat{U_\alpha^\prime}$
to be the restriction of this map to $\widehat{U_\alpha}$.

Since a groupoid homomorphism commutes with the mapping $\omega$, it
follows that the image of any $(\xi,x)\in
T(U_\alpha,U_\beta)\times\widehat{U_\alpha}$ under the
homomorphism is $(\xi^\prime,f_\alpha(x))$ for some $\xi^\prime\in
T(U_\alpha^\prime,U_\beta^\prime)$. We define $\rho_{\beta\alpha}:
T(U_\alpha,U_\beta)\rightarrow T(U_\alpha^\prime,U_\beta^\prime)$
by setting $\rho_{\beta\alpha}(\xi)=\xi^\prime$. Then the
equations in (a) are a consequence of the fact that the
homomorphism commutes with the mapping $\omega$. The equations in
(b) are simply another way of saying that the homomorphism
commutes with composition in the groupoid.

Thus we have shown that any homomorphism from $\Gamma\{U_\alpha\}$
to $\Gamma\{U_{\alpha^\prime}^\prime\}$ gives rise to a
$(\{f_\alpha\},\{\rho_{\beta\alpha}\})$ with the claimed
properties, such that it may be written $(\xi,x)\mapsto
(\rho_{\beta\alpha}(\xi),f_\alpha(x))$ for any $(\xi,x)\in
T(U_\alpha,U_\beta)\times\widehat{U_\alpha}$. On the other hand,
any such a $(\{f_\alpha\},\{\rho_{\beta\alpha}\})$ with the claimed
properties defines a homomorphism from $\Gamma\{U_\alpha\}$ to
$\Gamma\{U_{\alpha^\prime}^\prime\}$ in this way, even without the
assumption that each $\widehat{U_\alpha}$ is connected.

\hfill $\Box$

We follow up with a few remarks.

\begin{re}
{\em

\vspace{1.5mm}

(1) In this paper, we shall only consider homomorphisms from
$\Gamma\{U_\alpha\}$ to $\Gamma\{U_{\alpha^\prime}^\prime\}$
which are in the form $(\{f_\alpha\},\{\rho_{\beta\alpha}\})$.
By Lemma 2.2.1, these will cover all the homomorphisms from
$\Gamma\{U_\alpha\}$ to $\Gamma\{U_{\alpha^\prime}^\prime\}$
if each local chart $\widehat{U_\alpha}$ is connected.

\vspace{1.5mm}

(2) For any $\xi_{12},\xi_{23},\cdots,\xi_{n-1,n}$, we denote by
$\xi_{12}^\prime,\xi_{23}^\prime,\cdots,\xi_{n-1,n}^\prime$ their
images under the mappings $\{\rho_{\beta\alpha}\}$. Then the equations
in Lemma 2.2.1 (a) imply that if $x$ lies in the domain of
$$
\phi_{\xi_{n-1,n}}\circ\cdots\circ\phi_{\xi_{23}}\circ\phi_{\xi_{12}},
$$
then $f_\alpha(x)$ lies in the domain of
$$
\phi_{\xi_{n-1,n}^\prime}\circ\cdots\circ\phi_{\xi_{23}^\prime}
\circ\phi_{\xi_{12}^\prime},
$$
where $\alpha$ is any index such that $\xi_{12}\in T(U_\alpha,U_\beta)$
for some $\beta$. This induces a mapping
$$
\underline{\{f_\alpha\}}:\Lambda(\xi_{12},\xi_{23},\cdots,\xi_{n-1,n})
\rightarrow \Lambda(\xi_{12}^\prime,\xi_{23}^\prime,\cdots,
\xi_{n-1,n}^\prime), \leqno (2.2.1)
$$
which is defined by the rule that if $x\in {\bf a}$, then $f_\alpha(x)\in
\underline{\{f_\alpha\}}({\bf a})$. With this understood, the equations
in Lemma 2.2.1 (b) may be equivalently written as
$$
\rho_{\gamma\alpha}(\eta\circ\xi({\bf a}))=\rho_{\gamma\beta}(\eta)
\circ\rho_{\beta\alpha}(\xi)(\underline{\{f_\alpha\}}({\bf a})),
\hspace{2mm} \forall {\bf a}\in\Lambda(\xi,\eta). \leqno (2.2.2)
$$

\vspace{1.5mm}

(3) Note that the equations in Lemma 2.2.1 (b) imply that
$\rho_{\alpha\alpha}:T(U_\alpha,U_\alpha)\rightarrow T(U_\alpha^\prime,
U_\alpha^\prime)$ is in fact a homomorphism from $G_{U_\alpha}$ to
$G_{U_\alpha^\prime}$ for each $\alpha$, which will be denoted by
$\rho_\alpha:G_{U_\alpha}\rightarrow G_{U_\alpha^\prime}$ throughout.
Moreover, each $f_\alpha$ is $\rho_\alpha$-equivariant, and we have
$$
(\phi_{\rho_{\beta\alpha}(\xi)},\lambda_{\rho_{\beta\alpha}(\xi)})
\circ (f_\alpha,\rho_\alpha)=(f_\beta,\rho_\beta)\circ (\phi_\xi,
\lambda_\xi) \leqno (2.2.3)
$$
over $(\mbox{Domain }(\phi_\xi),\mbox{Domain }(\lambda_\xi))$.
(See $(2.1.3)$ for the definition of $\lambda_\xi$.)
In other words, the set of local `equivariant' maps $\{(f_\alpha,
\rho_\alpha)\}$ is compatible with respect to the
`transformations of local charts'.
}
\end{re}

\hfill $\Box$

A map of orbispaces to be defined in this paper will be a certain
equivalence class of groupoid homomorphisms of the form described in
Lemma 2.2.1. First of all, it proves to be convenient to introduce

\begin{defi}
A homomorphism $(\{f_\alpha^{(2)}\},\{\rho_{\beta\alpha}^{(2)}\}):
\Gamma\{U_\alpha\}\rightarrow\Gamma\{U_{\alpha^\prime}^\prime\}$ is
said to be conjugate to $(\{f_\alpha^{(1)}\},
\{\rho_{\beta\alpha}^{(1)}\})$ via $\{g_\alpha\}$, where each $g_\alpha$
is an element of $G_{U_\alpha^\prime}$, if
$$
f_\alpha^{(2)}=g_\alpha\circ f_\alpha^{(1)} \mbox{ and }
\rho_{\beta\alpha}^{(2)}(\xi)=g_\beta\circ\rho_{\beta\alpha}^{(1)}(\xi)
\circ g_\alpha^{-1}
$$
are satisfied for all $\alpha,\beta$ and all $\xi\in T(U_\alpha,U_\beta)$.
(Note that $\#\Lambda(g_\alpha^{-1},\rho_{\beta\alpha}^{(1)}(\xi),g_\beta)
=1$ so that the notation $g_\beta\circ\rho_{\beta\alpha}^{(1)}(\xi)
\circ g_\alpha^{-1}$ has no ambiguity.)
\end{defi}

For the rest of this paper, we shall make the following harmless
assumption on the domain orbispace $X$: in the atlas of local charts
$\U=\{U_i\}$, each $\widehat{U_i}$ is connected.

\vspace{1.5mm}

Now suppose $(\{f_\alpha\},\{\rho_{\beta\alpha}\}):\Gamma\{U_\alpha\}
\rightarrow\Gamma\{U_{\alpha^\prime}^\prime\}$ is a homomorphism where
$\{U_\alpha\}$ is a cover of $X$. Denote by $f:X\rightarrow X^\prime$
the induced map between the underlying spaces. Let $\{U_a\}\subset
\U$, $\{U^\prime_{a^\prime}\}\subset\U^\prime$ be any subsets, with a
correspondence $U_a\mapsto U_a^\prime\in\{U^\prime_{a^\prime}\}$ which
satisfies $f(U_a)\subset U_a^\prime$. Moreover, $\{U_a\}$ is a refinement
of $\{U_\alpha\}$ and is also a cover of $X$.

\begin{lem}
The homomorphism $(\{f_\alpha\},\{\rho_{\beta\alpha}\})$ canonically
induces a family of homomorphisms $(\{f_a\},\{\rho_{ba}\}):
\Gamma\{U_a\}\rightarrow\Gamma\{U_{a^\prime}^\prime\}$, which are
mutually conjugate in the sense of Definition 2.2.3.
\end{lem}

\pf
Since $\{U_a\}$ is a refinement of $\{U_\alpha\}$, there is a mapping
$\theta:a\mapsto\alpha=\theta(a)$ such that $U_a\subset U_{\theta(a)}$.
We fix a choice $\theta$ of such mappings to start with. By the assumption
that $U_a\subset U_{\theta(a)}$ and that each $\widehat{U_a}$ is connected,
we see that $T(U_a,U_{\theta(a)})$ is nonempty and $\mbox{Domain }(\phi_\xi)
=\widehat{U_a}$ for each $\xi\in T(U_a,U_{\theta(a)})$. On the other hand,
from the assumption that $f(U_a)\subset U_a^\prime$, we deduce that
$f(U_a)\subset U_a^\prime\cap U_{\theta(a)}^\prime$, and hence $T(U_a^\prime,
U_{\theta(a)}^\prime)\neq\emptyset$. For each $a$, we pick a
$\xi_a\in T(U_a,U_{\theta(a)})$, and then fix a choice $\xi_a^\prime$
of elements $\xi^\prime\in T(U_a^\prime,U_{\theta(a)}^\prime)$ satisfying
$$
f_{\theta(a)}\circ\phi_{\xi_a}(\widehat{U}_a)\subset
\mbox{Range }(\phi_{\xi^\prime}). \leqno (2.2.4)
$$

Now with a fixed choice of data
$(\theta,\{\xi_a\},\{\xi_a^\prime\})$, we define maps $\{f_a\}$, $f_a:
\widehat{U_a}\rightarrow\widehat{U_a^\prime}$, by
$$
f_a=(\phi_{\xi_a^\prime})^{-1}\circ
f_{\theta(a)}\circ\phi_{\xi_a},
\leqno (2.2.5)
$$
and for each $(a,b)$ with $U_a\cap U_b\neq\emptyset$, define a mapping
$\rho_{ba}:T(U_a,U_b)\rightarrow T(U_a^\prime,U_b^\prime)$ as follows.
For any $\eta\in T(U_a,U_b)$, we set $\theta(\eta)=\xi_b\circ\eta\circ\xi_a^{-1}
\in T(U_{\theta(a)},U_{\theta(b)})$ (note that $\#\Lambda(\xi_a^{-1},\eta,\xi_b)
=1$), and define
$$
\rho_{ba}(\eta)=(\xi_b^\prime)^{-1}\circ\rho_{\theta(b)\theta(a)}(\theta(\eta))
\circ\xi_a^\prime(x), \hspace{2mm} \forall x\in f_a(\mbox{Domain
}(\phi_\eta)). \leqno (2.2.6)
$$
One can easily check that $\phi_{\rho_{ba}(\eta)}\circ f_a(x)
=f_b\circ\phi_\eta(x)$ for any $\eta\in T(U_a,U_b)$, $x\in\mbox{Domain }
(\phi_\eta)$. For $(\{f_a\},\{\rho_{ba}\})$ to be a homomorphism from
$\Gamma\{U_a\}$ to $\Gamma\{U_{a^\prime}^\prime\}$, it remains to show that
$$
\rho_{ca}(\zeta\circ\eta({\bf a}))=\rho_{cb}(\zeta)\circ\rho_{ba}(\eta)
(\underline{\{f_a\}}({\bf a})) \leqno (2.2.7)
$$
for any $\eta\in T(U_a,U_b)$, $\zeta\in T(U_b,U_c)$ and
${\bf a}\in\Lambda(\eta,\zeta)$.

Observe that $\Lambda(\eta,\zeta)$ may be identified with
$\Lambda(\xi_a^{-1},\eta,\zeta,\xi_c)
=\Lambda(\xi_a^{-1},\eta,\xi_b,\xi_b^{-1},\zeta,\xi_c)$ via ${\bf a}
\mapsto \phi_{\xi_a}({\bf a})$. We deduce that
\begin{eqnarray*}
\theta(\zeta\circ\eta({\bf a})) & = &
\xi_c\circ\zeta\circ\eta({\bf a})\circ\xi_a^{-1} \\
& = & \xi_c\circ\zeta\circ\eta\circ\xi_a^{-1}(\phi_{\xi_a}({\bf a})) \\
& = & \xi_c\circ\zeta\circ\xi_b^{-1}\circ\xi_b\circ
\eta\circ\xi_a^{-1}(\phi_{\xi_a}({\bf a})) \\
& = & \theta(\zeta)\circ\theta(\eta)(\theta({\bf a}))
\end{eqnarray*}
for any ${\bf a}\in\Lambda(\eta,\zeta)$, where $\theta({\bf a})$ stands
for the unique element of $\Lambda(\theta(\eta),\theta(\zeta))$ which
contains $\phi_{\xi_a}({\bf a})$. Now we fix a $z\in {\bf a}$, and by
$(2.2.6)$, we have
\begin{eqnarray*}
&   &  \rho_{ca}(\zeta\circ\eta({\bf a}))\\
& = & (\xi_c^\prime)^{-1}\circ
\rho_{\theta(c)\theta(a)}(\theta(\zeta\circ\eta({\bf a})))\circ\xi_a^\prime
(f_a(z))\\
& = & (\xi_c^\prime)^{-1}\circ\rho_{\theta(c)\theta(a)}(\theta(\zeta)
\circ\theta(\eta)(\theta({\bf a})))\circ\xi_a^\prime
(f_a(z))\\
& = & (\xi_c^\prime)^{-1}\circ\rho_{\theta(c)\theta(b)}(\theta(\zeta))
\circ\rho_{\theta(b)\theta(a)}(\theta(\eta))(\underline{\{f_\alpha\}}
(\theta({\bf a})))\circ\xi_a^\prime
(f_a(z))\\
& = & (\xi_c^\prime)^{-1}\circ\rho_{\theta(c)\theta(b)}(\theta(\zeta))\circ
\rho_{\theta(b)\theta(a)}(\theta(\eta))\circ \xi_a^\prime
(f_a(z))\\
& = & (\xi_c^\prime)^{-1}\circ\rho_{\theta(c)\theta(b)}(\theta(\zeta))
\circ\xi_b^\prime\circ (\xi_b^\prime)^{-1}\circ\rho_{\theta(b)\theta(a)}
(\theta(\eta))\circ \xi_a^\prime (f_a(z))\\
& = & ((\xi_c^\prime)^{-1}\circ\rho_{\theta(c)\theta(b)}(\theta(\zeta))
\circ\xi_b^\prime(f_b(\phi_\eta(z))))\circ ((\xi_b^\prime)^{-1}\circ
\rho_{\theta(b)\theta(a)}
(\theta(\eta))\circ\xi_a^\prime(f_a(z)))(f_a(z))\\
& = & \rho_{cb}(\zeta)\circ\rho_{ba}(\eta)(\underline{\{f_a\}}({\bf a})).
\end{eqnarray*}
We thus verified $(2.2.7)$, and we conclude that $(\{f_a\},\{\rho_{ba}\}):
\Gamma\{U_a\}\rightarrow\Gamma\{U_{a^\prime}^\prime\}$ is a
homomorphism.

It remains to examine the dependence of $(\{f_a\},\{\rho_{ba}\})$ on
the choice of $(\theta,\{\xi_a\},\{\xi_a^\prime\})$. Suppose
$\bar{\theta}:a\mapsto\alpha$ is another mapping of indices
such that $U_a\subset U_{\bar{\theta}(a)}$. For any choice of
$\bar{\xi}_a\in T(U_a,U_{\bar{\theta}(a)})$, there is a
$\bar{\xi}_a^\prime\in T(U_a^\prime,U_{\bar{\theta}(a)}^\prime)$
defined by
$$
\bar{\xi}_a^\prime=\rho_{\bar{\theta}(a)\theta(a)}(\bar{\xi}_a\circ
\xi_a^{-1})\circ\xi_a^\prime ({\bf a}_a), \leqno (2.2.8)
$$
where ${\bf a}_a\in\Lambda(\xi_a^\prime,\rho_{\bar{\theta}(a)\theta(a)}
(\bar{\xi}_a\circ\xi_a^{-1}))$ is the element containing the image
of $f_a$, $f_a(\widehat{U_a})$. Then $\bar{\xi}_a^\prime$ satisfies
$(2.2.4)$ with respect to $\bar{\xi}_a$, i.e.,
$f_{\bar{\theta}(a)}\circ\phi_{\bar{\xi}_a}(\widehat{U_a})\subset
\mbox{Range }(\phi_{\bar{\xi}_a^\prime})$. We define
$$
\bar{f}_a=(\phi_{\bar{\xi}_a^\prime})^{-1}\circ f_{\bar{\theta}(a)}
\circ\phi_{\bar{\xi}_a}
$$
and
$$
\bar{\rho}_{ba}(\eta)
=(\bar{\xi}_b^\prime)^{-1}\circ\rho_{\bar{\theta}(b)\bar{\theta}(a)}
(\bar{\theta}(\eta))\circ\bar{\xi}_a^\prime(x),
\; \forall x\in\bar{f}_a(\mbox{Domain }(\phi_\eta)), \eta\in
T(U_a,U_b),
$$
where $\bar{\theta}(\eta)=\bar{\xi}_b\circ\eta\circ\bar{\xi}_a^{-1}$.
We shall verify that $\bar{f_a}=f_a$, $\bar{\rho}_{ba}=\rho_{ba}$.

Observe that $\mbox{Im }(f_{\bar{\theta}(a)}\circ \phi_{\bar{\xi}_a})
\subset\phi_{\bar{\xi}_a^\prime}({\bf a}_a)$, thus we have
\begin{eqnarray*}
\bar{f}_a & = & (\phi_{\bar{\xi}_a^\prime})^{-1}\circ f_{\bar{\theta}(a)}
\circ\phi_{\bar{\xi}_a} \\
          & = & (\phi_{\rho_{\bar{\theta}(a)\theta(a)}(\bar{\xi}_a\circ
\xi_a^{-1})\circ\xi_a^\prime ({\bf a}_a)})^{-1}\circ f_{\bar{\theta}(a)}
\circ\phi_{\bar{\xi}_a}\\
          & = & (\phi_{\xi_a^\prime})^{-1}\circ\phi_{\rho_{\theta(a)
          \bar{\theta}(a)}
(\xi_a\circ\bar{\xi}_a^{-1})}\circ f_{\bar{\theta}(a)}
\circ\phi_{\bar{\xi}_a} \\
          & = & (\phi_{\xi_a^\prime})^{-1}\circ f_{\theta(a)}\circ\phi_{\xi_a}
          =f_a.
\end{eqnarray*}
The verification of $\bar{\rho}_{ba}=\rho_{ba}$ is similar, hence we
leave it to the reader.

Thus the dependence of $(\{f_a\},\{\rho_{ba}\})$ on the choices of
$(\theta,\{\xi_a\},\{\xi_a^\prime\})$ boils down to the dependence on
the choices of $\{\xi_a^\prime\}$, with constrains in $(2.2.4)$, for some
fixed choice of $\{\xi_a\}$. It is easy to see that different choices
of $\{\xi_a^\prime\}$ differ by a pre-composition by an element
$g_a\in G_{U_a^\prime}$ for each index $a$. The resulting
homomorphisms are conjugate via $\{g_a\}$. Hence different choices of
$(\theta,\{\xi_a\},\{\xi_a^\prime\})$ will result in mutually conjugate
homomorphisms $(\{f_a\},\{\rho_{ba}\})$.

\hfill $\Box$

We remark that the family of homomorphisms $(\{f_a\},\{\rho_{ba}\})$
all induce the same map $f:X\rightarrow X^\prime$ between the
underlying spaces. On the other hand, the process $(\{f_\alpha\},
\{\rho_{\beta\alpha}\})\Rightarrow (\{f_a\},\{\rho_{ba}\})$ is
obviously transitive.

\begin{lem}
The following is indeed an equivalence relation on the set of all
groupoid homomorphisms: two homomorphisms are equivalent if
they induce a common family of mutually conjugate homomorphisms
in the sense of the preceding lemma.
\end{lem}

\pf
The only nontrivial part is transitivity, namely, for any
homomorphisms $\tau_1,\tau_2$ and $\tau_3$, if $\tau_1$ is
equivalent to $\tau_2$ and $\tau_2$ is equivalent to $\tau_3$,
then $\tau_1$ is equivalent to $\tau_3$. The proof goes as
follows. There are homomorphisms $\sigma_{12},\sigma_{23}$ which
are induced by $\tau_1, \tau_2$ and $\tau_2, \tau_3$ respectively.
By the transitivity of the process of inducing homomorphisms as
described in the preceding lemma, if there is a homomorphism $\kappa$
induced by both $\sigma_{12}$ and $\sigma_{23}$, then $\kappa$
must also be induced by both $\tau_1$ and $\tau_3$, which implies
that $\tau_1$ is equivalent to $\tau_3$.

Thus the problem boils down to show that $\sigma_{12}, \sigma_{23}$
induce a common family of mutually conjugate homomorphisms. To be
more explicit, let us assume $\sigma_{12}:\Gamma\{U_a\}\rightarrow
\Gamma\{U_{a^\prime}^\prime\}$ and $\sigma_{23}:\Gamma\{U_\alpha\}
\rightarrow\Gamma\{U_{\alpha^\prime}^\prime\}$. We then pick a
cover of $X$, denoted by $\{U_s\}$, where each $U_s$ is a
connected component of $U_a\cap U_\alpha$ for some indices $a$,
$\alpha$. Clearly $\{U_s\}$ is a refinement of both $\{U_a\}$ and
$\{U_\alpha\}$. On the other hand, there exists a
$\{U_{s^\prime}^\prime\}$, where each $U_{s^\prime}^\prime$ is a
connected component of $U_{a^\prime}^\prime\cap U_{\alpha^\prime}
^\prime$ for some indices $a^\prime$, $\alpha^\prime$, such that
each $U_s$ may be assigned with a $U_s^\prime\in\{U_{s^\prime}^\prime\}$
so that $f(U_s)\subset U_s^\prime$. (Here $f:X\rightarrow X^\prime$
is the common map of underlying spaces induced by
$\sigma_{12},\sigma_{23}$.) Now by the preceding lemma, there are
homomorphisms $\kappa_{12}, \kappa_{23}:\Gamma\{U_s\}\rightarrow
\Gamma\{U_{s^\prime}^\prime\}$ induced by $\sigma_{12}, \sigma_{23}$
respectively. Again by the transitivity of the process of inducing
homomorphisms as described in the preceding lemma, $\kappa_{12},
\kappa_{23}$ are also induced by $\tau_2$, hence must be conjugate
to each other. Thus we have shown that $\sigma_{12}, \sigma_{23}$
induce a common family of mutually conjugate homomorphisms.

\hfill $\Box$

\begin{defi}
A map of orbispaces from $X$ to $X^\prime$ is an equivalence class of
homomorphisms $(\{f_\alpha\},\{\rho_{\beta\alpha}\}):\Gamma\{U_\alpha\}
\rightarrow\Gamma\{U_{\alpha^\prime}^\prime\}$ in the sense of
Lemma 2.2.5, where $\{U_\alpha\}$ is a cover of $X$ by local charts.
\end{defi}

Now we conclude this subsection with the proof of Theorem 1.2.

\begin{thm}
With the notion of maps in Definition 2.2.6, the set of orbispaces
(satisfying {\em (C1), (C2)}) forms a category.
\end{thm}

\pf
It suffices to show that if $\Phi:X\rightarrow Y$, $\Psi:Y\rightarrow
Z$, then their composition $\Psi\circ\Phi:X\rightarrow Z$ is
well-defined and associative.

First of all, we fix a homomorphism $\tau=(\{g_a\},\{\eta_{ba}\}):
\Gamma\{V_a\}\rightarrow\Gamma\{W_{a^\prime}\}$ whose equivalence
class is the map $\Psi:Y\rightarrow Z$, and we shall prove that the
composition of $\Phi$ with $\tau$ is well-defined, which will be
denoted by $\tau\circ\Phi:X\rightarrow Z$.

Let $\sigma=(\{f_\alpha\},\{\rho_{\beta\alpha}\}):\Gamma\{U_\alpha\}
\rightarrow\Gamma\{V_{\alpha^\prime}\}$ be any homomorphism which
represents the map $\Phi:X\rightarrow Y$. Denote by $f:X\rightarrow Y$
the induced map between the underlying spaces. We consider the set
$\{U_i\}$ of all connected components of $U_\alpha\cap f^{-1}(V_a)$
for all $\alpha,a$. There is a mapping $\theta$ between the indices
$\{i\}$ and $\{\alpha\}$, $\theta:i\mapsto\alpha=\theta(i)$, such that
$U_i$ is a connected component of $U_{\theta(i)}\cap f^{-1}(V_a)$ for
some index $a$. The mapping $\theta$ defines the cover $\{U_i\}$ as a
refinement of $\{U_\alpha\}$. On the other hand, we fix a correspondence
$i\mapsto a=\hat{\theta}(i)$ where the indices $i$ and $\hat{\theta}(i)$
satisfy the condition that $U_i$ is a connected component of
$U_\alpha\cap f^{-1}(V_{\hat{\theta}(i)})$ for some index $\alpha$. We
set $\{V_{i^\prime}\}=\{V_{\hat{\theta}(i)}\}$, and assign $V_i=
V_{\hat{\theta}(i)}$ to $U_i$ for each index $i$, which verifies the
condition $f(U_i)\subset V_i$. Now by Lemma 2.2.4, $\sigma$ induces a
family of mutually conjugate homomorphisms from $\Gamma\{U_i\}$ to
$\Gamma\{V_{i^\prime}\}$. Let $\sigma^\prime=(\{f_i\},\{\rho_{ji}\})$
be one of the homomorphisms. We define the composition $\tau\circ\Phi:
X\rightarrow Z$ to be the equivalence class of the composition of
$\sigma^\prime$ with $\tau$, $\tau\circ\sigma^\prime=(\{h_i\},
\{\delta_{ji}\}):\Gamma\{U_i\}\rightarrow\Gamma\{W_{i^\prime}\}$,
where each $U_i$ is assigned with $W_i=W_{\hat{\theta}(i)}$,
$h_i=g_{\hat{\theta}(i)}\circ f_i$, and $\delta_{ji}=
\eta_{\hat{\theta}(j)\hat{\theta}(i)}\circ\rho_{ji}$.

The map $\tau\circ\Phi$ does not depend on the various choices made in
the construction, hence it is well-defined. First, it is independent of
the choice of $\sigma^\prime$, because a different one is conjugate to
$\sigma^\prime$ which results in a conjugate composition with $\tau$.
In particular, the choice of the mapping $\theta$ is irrelevant here.
Second, let's examine the dependence on the choice of $V_i$ made via
the mapping $i\mapsto a=\hat{\theta}(i)$. Suppose we have two choices
$\hat{\theta}_1$ and $\hat{\theta}_2$, which give rise to $\{V_{i,1}\}$
and $\{V_{i,2}\}$ by $V_{i,1}=V_{\hat{\theta}_1(i)}$ and
$V_{i,2}=V_{\hat{\theta}_2(i)}$. Let $\sigma_l^\prime=
(\{f_{i,l}\},\{\rho_{ji,l}\}):\Gamma\{U_i\}\rightarrow
\Gamma\{V_{i^\prime,l}\}$, $l=1,2$, be a choice of the corresponding
homomorphisms induced by $\sigma$, which is defined by $(2.2.5)$,
$(2.2.6)$ for some choices of $\xi_{i,l}=\xi_i\in T(U_i,U_{\theta(i)})$,
and $\xi_{i,l}^\prime\in T(V_{i,l},V_{\theta(i)})$ satisfying
$(2.2.4)$. Observe that there is an ${\bf a}\in\Lambda(\xi_{i,1}^\prime,
(\xi_{i,2}^\prime)^{-1})$ containing
$\phi_{\xi_{i,1}^\prime}^{-1}(\mbox{Im }(f_{\theta(i)}\circ\xi_i))$.
We define $\zeta_i=(\xi_{i,2}^\prime)^{-1}\circ\xi_{i,1}^\prime
({\bf a})\in T(V_{i,1},V_{i,2})$. Then by $(2.2.5)$, $(2.2.6)$, we have
$f_{i,2}=\zeta_i\circ f_{i,1}$ and $\rho_{ji,2}=\zeta_j\circ
\rho_{ji,1}\circ\zeta_i^{-1}$. The corresponding compositions with
$\tau$, $\tau\circ\sigma_l^\prime=(\{h_{i,l}\},\{\delta_{ji,l}\}):
\Gamma\{U_i\}\rightarrow\Gamma\{W_{i^\prime,l}\}$, $l=1,2$, which are
related by $h_{i,2}=\zeta_i^\prime\circ h_{i,1}$ and $\delta_{ji,2}
=\zeta_j^\prime\circ\delta_{ji,1}\circ (\zeta_i^\prime)^{-1}$ with
$\zeta_i^\prime=\eta_{\hat{\theta}_2(i)\hat{\theta}_1(i)}(\zeta_i)
\in T(W_{i,1},W_{i,2})$, are equivalent homomorphisms. Hence the choice on
$\{V_i\}$ is also irrelevant. Finally, suppose $\sigma$ is replaced by
a $\sigma_1$ which is induced by $\sigma$. Then the corresponding
induced homomorphism $\sigma_1^\prime$, which is from $\Gamma\{U_x\}$ to
$\Gamma\{V_{i^\prime}\}$ where $\{U_x\}$ is a refinement of $\{U_i\}$
through a mapping $x\mapsto i=\theta^\prime(x)$, and each $U_x$ is assigned
with $V_x=V_{\theta^\prime(x)}$, is induced by $\sigma^\prime:\Gamma\{U_i\}
\rightarrow\Gamma\{V_{i^\prime}\}$. It is easily seen that
$\tau\circ\sigma_1^\prime$ is induced by $\tau\circ\sigma^\prime$, so that
the choice of $\sigma$ is irrelevant. Hence the map $\tau\circ\Phi$
is well-defined.

We define $\Psi\circ\Phi=\tau\circ\Phi$. In order to see that
$\Psi\circ\Phi$ is independent of the choice on $\tau$,
we replace $\tau$ in the construction of $\tau\circ\Phi$ by a
homomorphism $\tau_1$ which is induced by $\tau$, say via the data
$(\iota,\{\xi_e\},\{\xi_e^\prime\})$. Here $\tau_1$ is from $\Gamma\{V_e\}$
to $\Gamma\{W_{e^\prime}\}$, where $\{V_e\}$ is a refinement of $\{V_a\}$
via a mapping $e\mapsto a=\iota(e)$. Now we let $\{U_x\}$ be the set of
connected components of $U_\alpha\cap f^{-1}(V_e)$ for all indices $\alpha,e$.
We define a refinement relation between $\{U_x\}$ and $\{U_i\}$ by the
following
rule: if $U_x$ is a connected component of $U_\alpha\cap f^{-1}(V_e)$,
then $U_i$ is the corresponding connected component of
$U_\alpha\cap f^{-1}(V_{\iota(e)})$ such that $U_x\subset U_i$.
Denote the corresponding mapping of indices by $x\mapsto i=\jmath(x)$.
We assign each $U_x$ with $V_x$, where $V_x=V_{\hat{\theta}(x)}$ for
some choice of mapping $x\mapsto e=\hat{\theta}(x)$. Then we assign each
$U_i$ with $V_i=V_{\iota(\hat{\theta}(x))}$. With these preparations, we
observe that a homomorphism $\sigma_1^\prime:\Gamma\{U_x\}\rightarrow
\Gamma\{V_{x^\prime}\}$ induced by $\sigma$ may be regarded as a
homomorphism induced by the induced homomorphism $\sigma^\prime:
\Gamma\{U_i\}\rightarrow\Gamma\{V_{i^\prime}\}$ of $\sigma$, say
via the data $(\jmath,\{\xi_x\},\{\xi_x^\prime\})$ where one may arrange to
have $\xi_x^\prime=\xi_{\hat{\theta}(x)}$ for some $\sigma_1^\prime$.
Now it is easy to see that $\tau_1\circ\sigma_1^\prime$ is induced by
$\tau\circ\sigma^\prime$ via the data $(\jmath,\{\xi_x\},
\{\xi_{\hat{\theta}(x)}^\prime\})$.
Thus we have verified that the composition $\Psi\circ\Phi$ is well-defined.

The compositions are associative:
$\tau\circ (\Psi\circ\Phi)
=(\tau\circ\Psi)\circ\Phi$ and
$\Xi\circ (\Psi\circ\Phi)
=(\Xi\circ\Psi)\circ\Phi$, which is clear from
the nature of construction.

\hfill $\Box$

\sectioni{Structure of mapping spaces}
\subsection{Some preliminary lemmas}

Let $\sigma=(\{f_\alpha\},\{\rho_{\beta\alpha}\}):\Gamma\{U_\alpha\}
\rightarrow\Gamma\{U_{\alpha^\prime}^\prime\}$ be any homomorphism
where $\{U_\alpha\}$ is a cover of $X$. We set
$$
G_\sigma=\left\{\begin{array}{lll}
g=\{g_\alpha\} & | & g_\alpha\in G_{U_\alpha^\prime} \mbox{ s.t.}\;
                   f_\alpha=g_\alpha\circ f_\alpha \mbox{ and }\\
               &   & \rho_{\beta\alpha}(\xi)=g_\beta\circ
                     \rho_{\beta\alpha}(\xi)\circ g_\alpha^{-1} \;
                     \forall \xi\in T(U_\alpha,U_\beta)
\end{array} \right\}. \leqno (3.1.1)
$$
The set $G_\sigma$ is naturally a group under $\{g_\alpha\}\{h_\alpha\}
=\{g_\alpha h_\alpha\}$. We call $G_\sigma$ the isotropy group
of $\sigma$. Denote by $f:X\rightarrow X^\prime$ the induced map
of $\sigma$ between the underlying spaces. For each connected component
$X_i$ of $X$, we pick a $q_i\in f(X_i)$, assuming $q_i\in U_{\alpha_i}^\prime$
for some index $\alpha_i$. Then there is an injective homomorphism $G_\sigma
\rightarrow\prod_i G_{q_i}$ defined by $\{g_\alpha\}\mapsto\prod_i
g_{\alpha_i}$, where $G_{q_i}$ stands for the isotropy group of
$q_i$.

\begin{lem}
To any pair $(\sigma,\tau)$ of equivalent homomorphisms, there is
associated a set $\Gamma_{\sigma\tau}$ with the following significance:
\begin{itemize}
\item [{(a)}] Each $\gamma\in\Gamma_{\sigma\tau}$ is assigned with
an isomorphism $\epsilon(\gamma):G_\tau\rightarrow G_\sigma$.
\item [{(b)}] $\Gamma_{\sigma\sigma}$ is canonically identified
with $G_\sigma$, under which $\epsilon(g)=Ad(g)$ for any $g\in G_\sigma$.
\item [{(c)}] There are mappings
$\Gamma_{\sigma\tau}\times\Gamma_{\tau\kappa}
\rightarrow\Gamma_{\sigma\kappa}$, denoted by
$(\gamma_2,\gamma_1)\mapsto\gamma_2\circ\gamma_1$, which are associative
and satisfy $\epsilon(\gamma_2\circ\gamma_1)=\epsilon(\gamma_2)
\circ\epsilon(\gamma_1)$, and when restricted to $\Gamma_{\sigma\sigma}$,
coincide with the multiplication in $G_\sigma$ under the canonical
identification $\Gamma_{\sigma\sigma}=G_\sigma$.
\item [{(d)}] The actions $G_\sigma\times\Gamma_{\sigma\tau}\rightarrow
\Gamma_{\sigma\tau}$ and $\Gamma_{\sigma\tau}\times G_{\tau}\rightarrow
\Gamma_{\sigma\tau}$ on $\Gamma_{\sigma\tau}$ are transitive.
\item [{(e)}] $\gamma\circ g=\epsilon(\gamma)(g)\circ\gamma$ holds
for any $g\in G_\tau$ and $\gamma\in\Gamma_{\sigma\tau}$.
\end{itemize}
\end{lem}

\pf
Let $\tau=(\{f_a\},\{\rho_{ba}\}):\Gamma\{U_a\}\rightarrow
\Gamma\{U_{a^\prime}^\prime\}$ and $\sigma=(\{f_\alpha\},
\{\rho_{\beta\alpha}\}):\Gamma\{U_\alpha\}\rightarrow\Gamma
\{U_{\alpha^\prime}^\prime\}$. We shall first prove the lemma for
the special case where $\tau$ is induced by $\sigma$, and for
any index $\alpha$, there is a subset $\{U_a|a\in I_\alpha\}$ of
$\{U_a\}$ such that
$$
U_\alpha=\bigcup_{a\in I_\alpha} U_a. \leqno (3.1.2)
$$
In this case,
we consider the set $\overline{\Gamma}_{\sigma\tau}$, which consists
of all $(\theta,\{\xi_a\},\{\xi_a^\prime\})$, where $\theta:a\mapsto
\alpha$ is a mapping of indices satisfying $U_a\subset U_{\theta(a)}$,
and $\xi_a\in T(U_a,U_{\theta(a)})$, $\xi_a^\prime\in T(U_a^\prime,
U_{\theta(a)}^\prime)$, such that the following equations
$$
f_a=(\phi_{\xi_a^\prime})^{-1}\circ f_{\theta(a)}\circ\phi_{\xi_a},
\hspace{2mm} \rho_{ba}(\eta)=(\xi^\prime_b)^{-1}\circ
(\rho_{\theta(b)\theta(a)}\circ\theta_{ba}(\eta))\circ\xi^\prime_a(x),
\leqno (3.1.3)
$$
are satisfied, where $\theta_{ba}:T(U_a,U_b)\rightarrow
T(U_{\theta(a)},U_{\theta(b)})$ is defined by $\eta\mapsto
\xi_b\circ\eta\circ\xi_a^{-1}$, and $x\in f_a(\mbox{Domain }(\phi_\eta))$,
cf. equations $(2.2.5)$ and $(2.2.6)$. The set $\Gamma_{\sigma\tau}$ is
$\overline{\Gamma}_{\sigma\tau}$ modulo the following equivalence relation:
$(\theta_1,\{\xi_{a,1}\},\{\xi_{a,1}^\prime\})$ and
$(\theta_2,\{\xi_{a,2}\},\{\xi_{a,2}^\prime\})$ are equivalent if and
only if
$$
\rho_{\theta_2(a)\theta_1(a)}(\xi_{a,2}\circ\xi_{a,1}^{-1})=
\xi_{a,2}^\prime\circ (\xi_{a,1}^\prime)^{-1}(x),\;
\forall x\in f_a(\widehat{U_a}). \leqno (3.1.4)
$$

One can easily verify that $(3.1.4)$ indeed defines an equivalence
relation, and that as in the proof of Lemma 2.2.4, $(3.1.3)$ is
preserved under this equivalence. When $\tau$ and $\sigma$ are identical,
$G_\tau$ is canonically identified with $\Gamma_{\tau\tau}$ by the
correspondence $\{g_a\}\mapsto (Id,\{1\},\{g_a\})$. Suppose
$\kappa:\Gamma\{U_i\}\rightarrow\Gamma\{U_{i^\prime}^\prime\}$ is
a homomorphism induced by $\tau$. We define the mappings $\Gamma_{\sigma\tau}
\times\Gamma_{\tau\kappa}\rightarrow\Gamma_{\sigma\kappa}$, denoted
by $(\gamma_2,\gamma_1)\mapsto \gamma_2\circ\gamma_1$, to be the ones
induced by
$$
((\theta,\{\xi_a\},\{\xi_a^\prime\}),(\iota,\{\eta_i\},\{\eta_i^\prime\}))
\mapsto (\theta\circ\iota,\{\xi_{\iota(i)}\circ\eta_i\},
\{\xi_{\iota(i)}^\prime\circ\eta_i^\prime\}). \leqno (3.1.5)
$$
It is a routine exercise to check that the equivalence defined by
$(3.1.4)$ is preserved under $(3.1.5)$. Hence the mappings
$\Gamma_{\sigma\tau}\times\Gamma_{\tau\kappa}\rightarrow
\Gamma_{\sigma\kappa}$ are well-defined, which are naturally
associative. Under the canonical identification $G_\tau=
\Gamma_{\tau\tau}$, we have $g_2\circ g_1=g_2g_1$
for any $g_1,g_2\in G_\tau$. The action $G_\sigma\times\Gamma_{\sigma\tau}
\rightarrow\Gamma_{\sigma\tau}$ is thus on the left and the action
$\Gamma_{\sigma\tau}\times G_\tau\rightarrow\Gamma_{\sigma\tau}$ is on the right.
The transitivity of the latter is part of Lemma 2.2.4, and the transitivity of
the former follows from the existence of isomorphisms $\epsilon(\gamma):G_\tau
\rightarrow G_\sigma$ which satisfy $\gamma\circ g
=\epsilon(\gamma)(g)\circ\gamma$
for any $g\in G_\tau$ and $\gamma\in\Gamma_{\sigma\tau}$. The construction
of $\epsilon(\gamma):G_\tau\rightarrow G_\sigma$ will occupy the next four
paragraphs.

For each $\gamma\in\Gamma_{\sigma\tau}$, we define the isomorphism
$\epsilon(\gamma):G_\tau\rightarrow G_\sigma$ as follows. First of
all, we fix a representative $\bar{\gamma}=(\theta,\{\xi_a\},
\{\xi_a^\prime\})$ of $\gamma$. Now suppose $\{g_a\}\in G_\tau$ is
any given element. Then for each $U_\alpha$, we assign each pair
$(U_a,\eta_a)$, where $U_a\in\{U_a|a\in I_\alpha\}$ and $\eta_a\in
T(U_a,U_\alpha)$, with an $\eta_a^\prime\in T(U_a^\prime,U_\alpha^\prime)$
defined by
$$
\eta_a^\prime=\rho_{\alpha\theta(a)}(\eta_a\circ\xi_a^{-1})
\circ\xi_a^\prime (x),\;\forall x\in f_a(\widehat{U_a}).
\leqno (3.1.6)
$$
With this understood, we define
$$
g(U_a,\eta_a)=\eta_a^\prime\circ g_a\circ (\eta_a^\prime)^{-1}(x),\;
\forall x\in f_a(\widehat{U_a}). \leqno (3.1.7)
$$
Note that $g(U_a,\eta_a)\in T(U_\alpha^\prime,U_\alpha^\prime)=
G_{U_\alpha^\prime}$. We shall prove that $g(U_a,\eta_a)$ depends
only on the index $\alpha$, whose common value is to be denoted by
$g_\alpha$, and that $\{g_\alpha\}$ is an element of $G_\sigma$.
We define $\epsilon(\gamma)$ by $\{g_a\}\mapsto \{g_\alpha\}$.

It is easy to see that $g(U_a,\eta_a)$ satisfies $g(U_a,\eta_a)
\circ f_\alpha|_{\mbox{Range }(\phi_{\eta_a})}
=f_\alpha|_{\mbox{Range }(\phi_{\eta_a})}$. Now let $(U_b,\eta_b)$
be another pair such that $\mbox{Range }(\phi_{\eta_a})\cap
\mbox{Range }(\phi_{\eta_b})\neq\emptyset$. Then
$\Lambda(\eta_a,\eta_b^{-1})\neq\emptyset$. Let ${\bf a}\in
\Lambda(\eta_a,\eta_b^{-1})$ be any element, and ${\bf a}^\prime\in
\Lambda(\eta_a^\prime,(\eta_b^\prime)^{-1})$ be the corresponding
element containing $f_a({\bf a})$. Then by the second equation in
$(3.1.3)$, we have
$$
\rho_{ba}(\eta_b^{-1}\circ\eta_a({\bf a}))=
(\eta_b^\prime)^{-1}\circ\rho_{\alpha\alpha}(\eta_b\circ (
\eta_b^{-1}\circ\eta_a({\bf a}))\circ\eta_a^{-1})\circ\eta_a^\prime
({\bf a}^\prime)=(\eta_b^\prime)^{-1}\circ\eta_a^\prime({\bf a}^\prime).
$$
With this in hand, the following holds on the appropriate components:
\begin{eqnarray*}
&   & (\eta_b^\prime)^{-1}\circ g(U_a,\eta_a)
=\rho_{ba}(\eta_b^{-1}\circ\eta_a)\circ g_a\circ (\eta_a^\prime)^{-1}\\
& = & g_b\circ\rho_{ba}(\eta_b^{-1}\circ\eta_a)\circ
(\eta_a^\prime)^{-1}=g_b\circ (\eta_b^\prime)^{-1}
=(\eta_b^\prime)^{-1}\circ g(U_b,\eta_b),
\end{eqnarray*}
which implies that $g(U_a,\eta_a)=g(U_b,\eta_b)$. Since
$U_\alpha=\bigcup_{a\in I_\alpha} U_a$ and $\widehat{U_\alpha}$ is
connected, we actually proved that $g(U_a,\eta_a)$ depends only on
$\alpha$, which is to be denoted by $g_\alpha$. It is clear that
$g_\alpha\circ f_\alpha=f_\alpha$. We have to check that
$\rho_{\beta\alpha}(\xi)=g_\beta\circ\rho_{\beta\alpha}(\xi)\circ
g_\alpha^{-1}$ for any $\xi\in T(U_\alpha,U_\beta)$.

Given any $\xi\in T_W(U_\alpha,U_\beta)$, since $U_\alpha=
\bigcup_{a\in I_\alpha} U_a$, there must be $a\in I_\alpha$
and $b\in I_\beta$ such that $U_a\cap U_b\neq\emptyset$ and has a
component contained in $W$, and there exist $\eta_a\in T(U_a,U_\alpha)$
and $\eta_b\in T(U_b,U_\beta)$ such that $\Lambda(\eta_a,\xi,\eta_b^{-1})
\neq\emptyset$. We set $\eta_{{\bf a}}=\eta_b^{-1}\circ\xi\circ
\eta_a({\bf a}),\; \forall {\bf a}\in\Lambda(\eta_a,\xi,\eta_b^{-1})$.
Let $\eta_a^\prime\in T(U_a^\prime,U_\alpha^\prime)$, $\eta_b^\prime\in
T(U_b^\prime,U_\beta^\prime)$ be the elements associated to
$\eta_a$, $\eta_b$ by $(3.1.6)$. Then we have $(\eta_b^\prime)^{-1}
\circ\rho_{\beta\alpha}(\xi)\circ\eta_a^\prime({\bf a}^\prime)
=\rho_{ba}(\eta_{{\bf a}})$, where ${\bf a}^\prime$ is the component
containing $f_a({\bf a})$. With $(3.1.7)$, $\rho_{\beta\alpha}(\xi)
=g_\beta\circ\rho_{\beta\alpha}(\xi)\circ g_\alpha^{-1}$ can be easily
deduced.

It is easy to see that $\epsilon(\gamma):G_\tau\rightarrow G_\sigma$
is independent of the choice made on the representative $\bar{\gamma}
=(\theta,\{\xi_a\},\{\xi_a^\prime\})$, and that $\epsilon(\gamma)$ is
an isomorphism. The other properties, i.e., (1) $\epsilon(g)=Ad(g)$ for
any $g\in G_\tau$ with $G_\tau$, $\Gamma_{\tau\tau}$ canonically
identified, (2) $\epsilon(\gamma_2\circ\gamma_1)=\epsilon(\gamma_2)
\circ\epsilon(\gamma_1)$, and (3) $\gamma\circ g=\epsilon(\gamma)(g)
\circ\gamma$ for any $g\in G_\tau$ and $\gamma\in\Gamma_{\sigma\tau}$,
can be easily verified from the construction. Thus we have completed the
proof of the lemma under the assumption that $\tau$ is induced by $\sigma$
and $(3.1.2)$ is satisfied.

To motivate the construction for the general case, we observe that
the mappings $\Gamma_{\sigma\tau}\times\Gamma_{\tau\kappa}\rightarrow
\Gamma_{\sigma\kappa}$, where $\tau$ is induced by $\sigma$, $\kappa$
is induced by $\tau$, and $(3.1.2)$ is satisfied for both pairs,
canonically identify $\Gamma_{\sigma\tau}$ with the orbit space
$(\Gamma_{\sigma\kappa}\times\Gamma_{\tau\kappa})/G_{\kappa}$ via
$\gamma\mapsto [\gamma\circ\gamma_0,\gamma_0],\; \forall\gamma_0\in
\Gamma_{\tau\kappa}$, where the action of $G_{\kappa}$ is given by
$(\gamma_1,\gamma_2)\cdot g=(\gamma_1\circ g,\gamma_2\circ g)$.
Now for any pair $(\sigma,\tau)$ of equivalent homomorphisms, let
$\kappa$ be a homomorphism induced by both $\tau$ and $\sigma$ such that
$(3.1.2)$ is satisfied. We define $\Gamma_{\sigma\tau}$ to be
$\Gamma_{\sigma\kappa}\times\Gamma_{\tau\kappa}$ modulo the
action of $G_\kappa$. Moreover, we define $\epsilon([\gamma_1,\gamma_2])
=\epsilon(\gamma_1)\circ\epsilon(\gamma_2)^{-1}$. The orbit space
$(\Gamma_{\sigma\kappa}\times\Gamma_{\tau\kappa})/G_\kappa$ for different
choices of $\kappa$ can be canonically identified, using the fact that
$\gamma\circ g=\epsilon(\gamma)(g)\circ\gamma$ for any $g\in G_\tau$ and
$\gamma\in\Gamma_{\sigma\tau}$ whenever $\tau$ is induced by $\sigma$ and
$(3.1.2)$ is satisfied. On the other hand, for any homomorphism $\kappa$
which is induced by $\tau$ with $(3.1.2)$ satisfied, we define a mapping
$\Gamma_{\tau\kappa}\times\Gamma_{\tau\kappa}\rightarrow G_{\kappa}$,
denoted by $(\gamma_1,\gamma_2)\mapsto \gamma_1\times\gamma_2$, as follows.
For any $\gamma_1,\gamma_2\in\Gamma_{\tau\kappa}$, we write
$\gamma_1=\gamma_0\circ g_1$ and $\gamma_2=\gamma_0\circ g_2$
for some $\gamma_0\in\Gamma_{\tau\kappa}$, $g_1,g_2\in G_\kappa$,
and then define $\gamma_1\times\gamma_2=g_1^{-1}g_2$. It follows
easily that $\gamma_1\times\gamma_2$ is well-defined and satisfies
$(\gamma_1\circ g_1)\times (\gamma_2\circ g_2)=g_1^{-1}
(\gamma_1\times\gamma_2)g_2$. Now to define the mappings
$\Gamma_{\sigma\tau}\times\Gamma_{\tau\kappa}\rightarrow
\Gamma_{\sigma\kappa}$ for the general case, we take a homomorphism
$\zeta$ induced by all $\sigma,\tau,\kappa$, and identify
$\Gamma_{\sigma\tau}$, $\Gamma_{\tau\kappa}$ and $\Gamma_{\sigma\kappa}$
as the orbit space $\Gamma_{\sigma\zeta}\times\Gamma_{\tau\zeta}/G_{\zeta}$,
$\Gamma_{\tau\zeta}\times\Gamma_{\kappa\zeta}/G_{\zeta}$ and
$\Gamma_{\sigma\zeta}\times\Gamma_{\kappa\zeta}/G_{\zeta}$ respectively.
Then for any $\gamma\in\Gamma_{\sigma\tau}$, $\gamma^\prime\in
\Gamma_{\tau\kappa}$, we write $\gamma=[\gamma_1,\gamma_2]$,
$\gamma^\prime=[\gamma_3,\gamma_4]$, and define
$\gamma\circ\gamma^\prime=[\gamma_1\circ (\gamma_2\times\gamma_3),
\gamma_4]$. We leave the verifications to the reader. (It is instructive
to think $[\gamma_1,\gamma_2]$ as $\gamma_1\circ\gamma_2^{-1}$
and $\gamma_1\times\gamma_2$ as $\gamma_1^{-1}\circ\gamma_2$.)

\hfill $\Box$

We remark that if we let $U(\Phi)=\{\sigma\mid\sigma \mbox{ is a
homomorphism representing } \Phi\}$ and $\Gamma(\Phi)=\bigcup_{(\sigma,\tau)}
\Gamma_{\sigma\tau}$, with $\alpha,\omega:\Gamma(\Phi)\rightarrow
U(\Phi)$ given by $\alpha(\gamma)=\tau$, $\omega(\gamma)=\sigma$,
$\forall \gamma\in\Gamma_{\sigma\tau}$, then the preceding lemma implies
that $\Gamma(\Phi)$ is a groupoid acting on $U(\Phi)$, with the set
of orbits $\Gamma(\Phi)\backslash U(\Phi)=\{\Phi\}$.

\begin{lem}
Let $(\{f_\alpha^{(1)}\},\{\rho_{\beta\alpha}^{(1)}\}),
(\{f_\alpha^{(2)}\},\{\rho_{\beta\alpha}^{(2)}\}):\Gamma\{U_\alpha\}
\rightarrow\Gamma\{U_{\alpha^\prime}^\prime\}$ be any pair of
equivalent homomorphisms, where the assignment $U_\alpha\mapsto
U_\alpha^\prime$ is common for both of them. Then $(\{f_\alpha^{(1)}\},
\{\rho_{\beta\alpha}^{(1)}\}),(\{f_\alpha^{(2)}\},
\{\rho_{\beta\alpha}^{(2)}\})$ are in fact conjugate to each
other.
\end{lem}

\pf
By the assumption that $(\{f_{\alpha}^{(1)}\},\{\rho_{\beta\alpha}^{(1)}\})$
and $(\{f_{\alpha}^{(2)}\},\{\rho_{\beta\alpha}^{(2)}\})$ are equivalent,
there exists a common induced homomorphism $(\{f_a\},\{\rho_{ba}\}):
\Gamma\{U_a\}\rightarrow\Gamma\{U_{a^\prime}^\prime\}$. Without loss
of generality, we assume that $(3.1.2)$ is satisfied: for any
index $\alpha$, there is a subset $\{U_a\mid a\in I_\alpha\}$ of $\{U_a\}$
such that $U_\alpha=\bigcup_{a\in I_\alpha} U_a$.

For $i=1,2$, we fix a set of data $(\theta_i,\{\xi_{a,i}\},
\{\xi_{a,i}^\prime\})$, where $\theta_i:a\mapsto\alpha$ is a
mapping of indices such that $U_a\subset U_{\theta_i(a)}$, and
$\xi_{a,i}\in T(U_a,U_{\theta_i(a)})$, $\xi_{a,i}^\prime
\in T(U_a^\prime,U_{\theta_i(a)}^\prime)$ so that
$$
f_a=(\xi_{a,i}^\prime)^{-1}\circ f_{\theta_i(a)}^{(i)}\circ\xi_{a,i}
$$
and
$$
\rho_{ba}(\eta)=(\xi_{b,i}^\prime)^{-1}\circ\rho_{\theta_i(b)\theta_i(a)}
^{(i)}(\xi_{b,i}\circ\eta\circ\xi_{a,i}^{-1})\circ\xi_{a,i}^\prime
(x) \;\forall x\in f_a(\mbox{Domain }(\phi_\eta)), \eta\in
T(U_a,U_b).
$$
Now given any $\alpha$ and $a\in I_\alpha$, we associate a
$\zeta_{a,i}^\prime\in T(U_a^\prime,U_\alpha^\prime)$ to each
$\zeta_a\in T(U_a,U_\alpha)$ by
$\zeta_{a,i}^\prime=\rho_{\alpha\theta_i(a)}^{(i)}
(\zeta_a\circ\xi_{a,i}^{-1})\circ\xi_{a,i}^\prime(x)$,
$\forall x\in f_a(\widehat{U_a})$. It is easy to
check that for any $a\in I_\alpha, b\in I_\beta$, one has
$$
f_a=(\zeta_{a,i}^\prime)^{-1}\circ f_{\alpha}^{(i)}\circ\zeta_{a}
$$
and
$$
\rho_{ba}(\eta)=(\zeta_{b,i}^\prime)^{-1}\circ\rho_{\beta\alpha}
^{(i)}(\zeta_{b}\circ\eta\circ\zeta_{a}^{-1})\circ\zeta_{a,i}^\prime
(x)\;\forall x\in f_a(\mbox{Domain }(\phi_\eta)), \eta\in
T(U_a,U_b).
$$
We define $g_\alpha(\zeta_a)=\zeta_{a,2}^\prime\circ (\zeta_{a,1}
^\prime)^{-1}(x)$, $\forall x\in f_a(\widehat{U_a})$. Then
$g_\alpha(\zeta_a)\in T(U_\alpha^\prime,U_\alpha^\prime)=
G_{U_\alpha^\prime}$. Observe that if $\mbox{Range }(\phi_{\zeta_a})
\cap\mbox{Range }(\phi_{\zeta_b})\neq\emptyset$, then $\Lambda(\zeta_a,
\zeta_b^{-1})\neq\emptyset$. We set $\eta_{\bf a}=\zeta_b^{-1}\circ
\zeta_a({\bf a})\in T(U_a,U_b)$ for any ${\bf a}\in\Lambda(\zeta_a,
\zeta_b^{-1})$. Then $\rho_{ba}(\eta_{\bf a})=(\zeta_{b,i}^\prime)^{-1}
\circ\zeta_{a,i}^\prime(x)$, $\forall x\in f_a({\bf a})$, for $i=1,2$.
This implies that $g_\alpha(\zeta_a)=g_\alpha(\zeta_b)$
whenever $\mbox{Range }(\phi_{\zeta_a})\cap\mbox{Range }(\phi_{\zeta_b})
\neq\emptyset$. But $\bigcup_{a\in I_\alpha} U_a=U_\alpha$, and
$\widehat{U_\alpha}$ is connected, so $g_\alpha(\zeta_a)$ is in fact
independent of $\zeta_a$. We define $g_\alpha=g_\alpha(\zeta_a)$.

It remains to verify that $(\{f_\alpha^{(2)}\},
\{\rho_{\beta\alpha}^{(2)}\})$ is conjugate to $(\{f_\alpha^{(1)}\},
\{\rho_{\beta\alpha}^{(1)}\})$ via $\{g_\alpha\}$. First,
$f_{\alpha}^{(2)}=g_\alpha\circ f_{\alpha}^{(1)}$ is clear from the
construction of $g_\alpha$. Second, given any $\xi\in T_W(U_\alpha,U_\beta)$,
since $U_\alpha=\bigcup_{a\in I_\alpha} U_a$, there must be $a\in I_\alpha$
and $b\in I_\beta$ such that $U_a\cap U_b\neq\emptyset$ and has a
connected component contained in $W$, and there exist $\zeta_a\in
T(U_a,U_\alpha)$ and $\zeta_b\in T(U_b,U_\beta)$ such that
$\Lambda(\zeta_a,\xi,\zeta_b^{-1})\neq\emptyset$. We set $\eta_{\bf a}
=\zeta_b^{-1}\circ\xi\circ\zeta_a({\bf a})$, $\forall {\bf a}\in
\Lambda(\zeta_a,\xi,\zeta^{-1}_b)$. Then $\rho_{ba}(\eta_{\bf a})=
(\zeta_{b,i}^\prime)^{-1}\circ\rho_{\beta\alpha}^{(i)}(\xi)\circ
\zeta_{a,i}^\prime(x)$, $\forall x\in f_a({\bf a})$, for $i=1,2$.
We deduce easily from these equations that
$\rho_{\beta\alpha}^{(2)}(\xi)=g_\beta\circ\rho_{\beta\alpha}^{(1)}(\xi)
\circ g_\alpha^{-1}$. Hence $(\{f_\alpha^{(2)}\},
\{\rho_{\beta\alpha}^{(2)}\})$ is conjugate to $(\{f_\alpha^{(1)}\},
\{\rho_{\beta\alpha}^{(1)}\})$ via $\{g_\alpha\}$.

\hfill $\Box$

\begin{lem}
Let $\Phi:X\rightarrow X^\prime$ be any map of orbispaces. Suppose
there are $\{U_\alpha\}$, $\{U_{\alpha^\prime}^\prime\}$ of local
charts with the following significance:
\begin{itemize}
\item $\{U_\alpha\}$ is a cover of $X$.
\item There is a correspondence $U_\alpha\mapsto U_\alpha^\prime\in
\{U_{\alpha^\prime}^\prime\}$ such that for each $U_\alpha$, the
restriction of $\Phi$ to the open subspace $U_\alpha$ of $X$ is
represented by a pair $(f_\alpha,\rho_\alpha):(\widehat{U_\alpha},
G_{U_\alpha})\rightarrow (\widehat{U_\alpha^\prime},G_{U_\alpha^\prime})$
where $f_\alpha$ is $\rho_\alpha$-equivariant.
\end{itemize}
Then there are mappings $\rho_{\beta\alpha}:T(U_\alpha,U_\beta)
\rightarrow T(U_\alpha^\prime,U_\beta^\prime)$ with
$\rho_{\alpha\alpha}=\rho_\alpha$ such that $(\{f_\alpha\},
\{\rho_{\beta\alpha}\}):\Gamma\{U_\alpha\}\rightarrow\Gamma\{
U_{\alpha^\prime}^\prime\}$ is a homomorphism whose equivalence
class is the given map $\Phi$.
\end{lem}

\pf
We pick a homomorphism $\tau=(\{f_a\},\{\rho_{ba}\}):\Gamma\{U_a\}
\rightarrow\Gamma\{U_{a^\prime}^\prime\}$ whose equivalence class
is $\Phi:X\rightarrow X^\prime$, where we may assume, by passing to
an induced homomorphism, that $(3.1.2)$ is satisfied: for each $\alpha$,
there is a subset $\{U_a|a\in I_\alpha\}\subset\{U_a\}$ such that
$U_\alpha=\bigcup_{a\in I_\alpha}U_a$. The restriction of $\tau$ to
$U_\alpha$, denoted by $\tau_\alpha=(\{f_a\},\{\rho_{ba}\})$
where $a,b\in I_\alpha$, represents the restriction of $\Phi$
to the open subspace $U_\alpha$. On the other hand, the restriction of
$\Phi$ to $U_\alpha$ is also represented by $(f_\alpha,\rho_\alpha)$.
We apply Lemma 2.2.4 first to obtain an induced homomorphism of
$(f_\alpha,\rho_\alpha)$ from $\Gamma\{U_a\}$ to
$\Gamma\{U_{a^\prime}^\prime\}$ where each $U_a\mapsto
U_a^\prime$, then apply Lemma 3.1.2 to conclude that $\tau_\alpha
=(\{f_a\},\{\rho_{ba}\})$ is conjugate to the induced one, so that
$\tau_\alpha$ is induced by $(f_\alpha,\rho_\alpha)$ also.
Hence there are $\xi_{a,\alpha}\in T(U_a,U_\alpha)$,
$\xi_{a,\alpha}^\prime\in T(U_a^\prime,U_\alpha^\prime)$, where
$a\in I_\alpha$, such that
$$
f_a=(\phi_{\xi_{a,\alpha}^\prime})^{-1}\circ f_\alpha\circ
\phi_{\xi_{a,\alpha}}\; \forall a\in I_\alpha, \leqno (3.1.8)
$$
and for any $\eta\in T(U_a,U_b)$, $a,b\in I_\alpha$,
$$
\rho_{ba}(\eta)=(\xi_{b,\alpha}^\prime)^{-1}\circ\rho_\alpha
(\xi_{b,\alpha}\circ\eta\circ\xi_{a,\alpha}^{-1})\circ
\xi_{a,\alpha}^\prime(x) \; \forall x\in f_a(\mbox{Domain }(\phi_\eta)).
\leqno (3.1.9)
$$
Note that if we replace $\xi_{a,\alpha}$ by $g\circ\xi_{a,\alpha}$ for
any $g\in G_{U_\alpha}$ and correspondingly replace
$\xi_{a,\alpha}^\prime$ by $\rho_\alpha(g)\circ\xi_{a,\alpha}^\prime$,
$(3.1.8)$ and $(3.1.9)$ continue to hold. In other words, $(3.1.8)$ and
$(3.1.9)$ hold with $\xi_{a,\alpha}$ allowed to be any element in
$T(U_a,U_\alpha)$ as long as $\xi_{a,\alpha}^\prime$ is chosen by the
correspondence $\xi_{a,\alpha}\mapsto\xi_{a,\alpha}^\prime$ described
above.

Now suppose $U_\alpha,U_\beta\in\{U_\alpha\}$ have nonempty intersection.
Given any $\zeta\in T(U_\alpha,U_\beta)$, since $U_\alpha
=\bigcup_{a\in I_\alpha}U_a$, $U_\beta=\bigcup_{b\in I_\beta}U_b$,
there are $a\in I_\alpha, b\in I_\beta$, with $\xi_{a,\alpha}\in
T(U_a,U_\alpha)$, $\xi_{b,\beta}\in T(U_b,U_\beta)$, and $\eta\in
T(U_a,U_b)$, such that $\zeta=\xi_{b,\beta}\circ\eta\circ
\xi_{a,\alpha}^{-1}$. We define
$$
\rho_{\beta\alpha}(\zeta)=\xi_{b,\beta}^\prime\circ\rho_{ba}(\eta)
\circ (\xi_{a,\alpha}^\prime)^{-1}(x) \; \forall
x\in\phi_{\xi_{a,\alpha}^\prime}(f_a(\mbox{Domain }(\phi_\eta))).
\leqno (3.1.10)
$$
In order to verify that $\rho_{\beta\alpha}(\zeta)$ is independent
of the choices made on $\xi_{a,\alpha},\eta,\xi_{b,\beta}$, suppose
we have two different choices $(\xi_{a,\alpha},\eta,\xi_{b,\beta})$
and $(\xi_{a^\prime,\alpha},\eta^\prime,\xi_{b^\prime,\beta})$. The
fact that $\mbox{Domain }(\phi_\zeta)$ and $\mbox{Range }(\phi_\zeta)$
are connected allows us to reduce the problem to the following special
cases: (1) $\xi_{a,\alpha}=\xi_{a^\prime,\alpha}$ and
$\mbox{Domain }(\phi_\eta)\cap\mbox{Domain }(\phi_{\eta^\prime})
\neq\emptyset$, (2) $\xi_{b,\beta}=\xi_{b^\prime,\beta}$ and
$\mbox{Range }(\phi_\eta)\cap\mbox{Range }(\phi_{\eta^\prime})
\neq\emptyset$. We shall only consider the first case, the second case
is completely parallel. Under the assumption that
$\mbox{Domain }(\phi_\eta)\cap\mbox{Domain }(\phi_{\eta^\prime})
\neq\emptyset$, there is an $\epsilon\in T(U_b,U_{b^\prime})$ such
that $\eta^\prime=\epsilon\circ\eta({\bf a})$ for some ${\bf a}\in
\Lambda(\eta,\epsilon)$. It follows from the fact $\xi_{b,\beta}
\circ\eta\circ\xi_{a,\alpha}^{-1}=\zeta=\xi_{b^\prime,\beta}\circ
\eta^\prime\circ\xi_{a^\prime,\alpha}^{-1}$ and the assumption
that $\xi_{a,\alpha}=\xi_{a^\prime,\alpha}$ that $\xi_{b,\beta}
=\xi_{b^\prime,\beta}\circ\epsilon(\phi_\eta({\bf a}))$, which
in turn implies that
$\rho_{b^\prime b}(\epsilon)=(\xi_{b^\prime,\beta}^\prime)^{-1}
\circ\xi_{b,\beta}^\prime(x)$, $\forall x\in f_b(\phi_\eta
({\bf a}))$, by $(3.1.9)$. With this relation,
$\xi_{b,\beta}^\prime\circ\rho_{ba}(\eta)\circ
(\xi_{a,\alpha}^\prime)^{-1}(\phi_{\xi_{a,\alpha}^\prime}(f_a(x)))
=\xi_{b^\prime,\beta}^\prime\circ\rho_{b^\prime a^\prime}(\eta^\prime)
\circ (\xi_{a^\prime,\alpha}^\prime)^{-1}(\phi_{\xi_{a^\prime,\alpha}^\prime}
(f_{a^\prime}(x)))$, $\forall x\in {\bf a}$, follows by recalling
$\rho_{b^\prime a^\prime}(\eta^\prime)=\rho_{b^\prime b}(\epsilon)\circ
\rho_{ba}(\eta)(\underline{\{f_a\}}({\bf a}))$. Hence $\rho_{\beta\alpha}
(\zeta)$ is well-defined.

It remains to verify that $(\{f_\alpha\},\{\rho_{\beta\alpha}\})
:\Gamma\{U_\alpha\}\rightarrow\Gamma\{U_{\alpha^\prime}^\prime\}$ is a
homomorphism whose equivalence class is $\Phi$, and $\rho_{\alpha\alpha}
=\rho_\alpha$ for each $\alpha$. The latter follows easily by
comparing $(3.1.9)$ with $(3.1.10)$. In order to show that $(\{f_\alpha\},
\{\rho_{\beta\alpha}\})$ is a homomorphism, we need to check (a)
$f_\beta\circ\phi_{\zeta}=\phi_{\rho_{\beta\alpha}(\zeta)}\circ
f_\alpha, \forall\zeta\in T(U_\alpha,U_\beta)$, (b) $\rho_{\gamma\alpha}
(\delta\circ\zeta({\bf a}))=\rho_{\gamma\beta}(\delta)\circ\rho_{\beta\alpha}
(\zeta)(\underline{\{f_\alpha\}}({\bf a})), \forall {\bf a}\in
\Lambda(\zeta,\delta), \zeta\in T(U_\alpha,U_\beta),\delta\in
T(U_\beta,U_\gamma)$. The former follows directly from $(3.1.10)$ and
$(3.1.8)$. For the latter, we first observe that we can choose
$a\in I_\alpha, b\in I_\beta, c\in I_\gamma$ such that $\zeta=
\xi_{b,\beta}\circ\eta\circ\xi_{a,\alpha}^{-1}, \delta=\xi_{c,\gamma}
\circ\epsilon\circ\xi_{b,\beta}^{-1}$ for some $\xi_{a,\alpha}\in
T(U_a,U_\alpha)$, $\xi_{b,\beta}\in T(U_b,U_\beta)$, $\xi_{c,\gamma}\in
T(U_c,U_\gamma)$ and some $\eta\in T(U_a,U_b)$, $\epsilon\in T(U_b,U_c)$
satisfying $\mbox{Range }(\phi_\eta)\cap\mbox{Domain }(\phi_\epsilon)
\neq\emptyset$. Moreover, there is an element ${\bf b}\in\Lambda(\eta,
\epsilon)$ contained in $\phi_{\xi_{a,\alpha}}^{-1}({\bf a})$ such that
$\delta\circ\zeta({\bf a})=\xi_{c,\gamma}\circ (\epsilon\circ\eta({\bf b}))
\circ \xi_{a,\alpha}^{-1}$. With this relation, one can easily deduce (b)
from $(3.1.10)$ and the equations $\rho_{ca}(\epsilon\circ\eta({\bf b}))
=\rho_{cb}(\epsilon)\circ\rho_{ba}(\eta)(\underline{\{f_a\}}({\bf b}))$.
Finally, it is straightforward from $(3.1.8)$, $(3.1.9)$ and $(3.1.10)$
that $\tau=(\{f_a\},\{\rho_{ba}\})$ is induced by $(\{f_\alpha\},
\{\rho_{\beta\alpha}\})$, hence the equivalence class of the latter is
also the given map $\Phi$.

\hfill $\Box$

\subsection{Proof of Theorem 1.3}

First of all, we give the precise definition of the assumption
on the domain orbispace $X$ in Theorem 1.3. An orbispace $X$
is said to be paracompact, locally compact and Hausdorff provided that
the following conditions are satisfied:
\begin{itemize}
\item The underlying space $X$ is paracompact, locally compact,
and Hausdorff.
\item Each local chart $\widehat{U_i}$ of $X$ is locally compact
and $\pi_{U_i}:\widehat{U_i}\rightarrow X$ is proper.
\item Each local chart $\widehat{U_i}$ of $X$ is Hausdorff.
\end{itemize}
These conditions will be imposed on $X$ throughout this subsection.
Finally, we shall fix the notation $[X;X^\prime]$ for the set of maps
of orbispaces from $X$ to $X^\prime$.

\begin{lem}
Any open cover of $X$ has a refinement
$\{U_\alpha\mid\alpha\in\Lambda\}$, which is locally finite and
the closure $\overline{U_\alpha}$ is compact. Furthermore, there
is an open cover $\{V_\alpha\mid\alpha\in\Lambda\}$ such that
$\overline{V_\alpha}\subset U_\alpha$ for all $\alpha\in\Lambda$.
\end{lem}

\pf
Let $\{U_i\}$ be any given open cover. Since the underlying space
$X$ is locally compact and Hausdorff, there is a refinement $\{U_a\}$
where each $U_a$ has a compact closure. On the other hand, the
underlying space $X$ is paracompact, so $\{U_a\}$ has a locally finite
refinement $\{U_\alpha|\alpha\in\Lambda\}$. It is easy to see that
$\{U_\alpha|\alpha\in\Lambda\}$ is a desired refinement of $\{U_i\}$.

We construct $\{V_\alpha\}$ as follows. Let $S$ be the set of the subsets
$I$ of $\Lambda$ satisfying the following condition: for each $\alpha\in I$,
there exists a $V_\alpha$ such that $\overline{V_\alpha}\subset U_\alpha$
and $X=\cup_{\alpha\in I} V_\alpha\cup_{\alpha\in \Lambda\setminus I}
U_\alpha$. We shall prove that $\Lambda\in S$.

First of all, $S$ is nonempty. To see this, we pick an $\alpha_0\in\Lambda$
and set $W_{\alpha_0}=X\setminus\cup_{\alpha\neq\alpha_0} U_\alpha$,
which is a compact subset of $U_{\alpha_0}$. Since the underlying space
$X$ is Hausdorff, and $\overline{U_{\alpha_0}}\setminus U_{\alpha_0}$ is
compact, we see that for any $p\in W_{\alpha_0}$, there is a neighborhood
$O_p$ of $p$ with $\overline{O_p}\subset U_{\alpha_0}$. The compactness of
$W_{\alpha_0}$ implies that there are finitely many $p_1,\cdots, p_n
\in W_{\alpha_0}$ such that $W_{\alpha_0}\subset\cup_{i=1}^n O_{p_i}$.
We take $V_{\alpha_0}=\cup_{i=1}^n O_{p_i}$. Then clearly
$\overline{V_{\alpha_0}}\subset U_{\alpha_0}$ and
$X=\cup_{\alpha\neq\alpha_0} U_\alpha\cup V_{\alpha_0}$.
Hence $\{\alpha_0\}\in S$.

Secondly, $S$ is partially ordered: $I\leq I^\prime$ iff $I\subset
I^\prime$. Let $T\subset S$ be a linearly ordered subset. We shall
prove that $I(T)=\cup_{I\in T}I\in S$, or equivalently, $X=\cup_{\alpha
\in I(T)} V_\alpha\cup_{\alpha\in \Lambda\setminus I(T)} U_\alpha$.
Given any $p\in X$, since $\{U_\alpha|\alpha\in\Lambda\}$ is locally
finite, there are only finitely many elements, say $U_{\alpha_1},\cdots,
U_{\alpha_k}$, which contains $p$. We need to show that if
$\alpha_i\in I(T)$ for all $1\leq i\leq k$, then $p\in\cup_{\alpha\in
I(T)} V_{\alpha}$. To this end, we observe that since $T$ is
linearly ordered, there exists an $I\in T$ containing all $\alpha_i,
1\leq i\leq k$. Then $p\in X=\cup_{\alpha\in I} V_\alpha\cup_{\alpha\in
\Lambda\setminus I} U_\alpha$ implies that $p\in\cup_{\alpha\in I}
V_\alpha\subset\cup_{\alpha\in I(T)} V_\alpha$. Hence $I(T)\in S$.

Finally, we apply Zorn's Lemma to $S$ to conclude that there is a
maximal element $I_0$ in $S$. Then a verbatim proof of $S\neq\emptyset$
shows that $I_0=\Lambda$. Hence $\Lambda\in S$.

\hfill $\Box$

We need to consider a class of special local charts $\{U_i\}$ on $X$,
which are dictated by the following conditions: (1) the closure
$\overline{U_i}$ is a compact subset, (2) there is a space
$\widehat{\overline{U_i}}$ containing $\widehat{U_i}$ as its interior,
with an action of $G_{U_i}$ extending the action of $G_{U_i}$ on
$\widehat{U_i}$, and a map $\pi_{\overline{U_i}}:\widehat{\overline{U_i}}
\rightarrow X$ extending $\pi_{U_i}$, such that
$\widehat{\overline{U_i}}/G_{U_i}$ is homeomorphic to $\overline{U_i}$
under $\pi_{\overline{U_i}}$, and (3) for any pair of such $U_i,U_j$
with $U_i\cap U_j\neq\emptyset$, the homeomorphism $\phi_\xi$ associated
to any $\xi\in T(U_i,U_j)$ can be extended to a homeomorphism between
the closures of its domain and range. We call this class of local charts
`admissible'.

In light of Lemma 3.2.1, we may assume without loss of generality for
the remaining subsections that in any homomorphisms $(\{f_\alpha\},
\{\rho_{\beta\alpha}\}):\Gamma\{U_\alpha\}\rightarrow
\Gamma\{U_{\alpha^\prime}^\prime\}$ under consideration, the following
conditions are satisfied: (1) the cover $\{U_\alpha\}$ is locally
finite and each local chart $U_\alpha$ is admissible, (2) each $f_\alpha$
can be extended over (uniquely) to the compact closure
$\widehat{\overline{U_\alpha}}$. We call this class of homomorphisms
`admissible'.

Let $(\{f_\alpha^{(0)}\},\{\rho_{\beta\alpha}\}):\Gamma\{U_\alpha\}
\rightarrow\Gamma\{U_{\alpha^\prime}^\prime\}$ be any admissible
homomorphism. We define
$$
\O_{\{\rho_{\beta\alpha}\}}=
\{\sigma\mid\sigma=(\{f_\alpha\},\{\rho_{\beta\alpha}\}):
\Gamma\{U_\alpha\}\rightarrow\Gamma\{U_{\alpha^\prime}^\prime\}
\mbox{ is admissible}\} \leqno (3.2.1)
$$
and
$$
G_{\{\rho_{\beta\alpha}\}}=
\left\{\begin{array}{lll}
g=\{g_\alpha\} & | & g_\alpha\in G_{U_\alpha^\prime},\;
               \rho_{\beta\alpha}(\xi)=g_\beta\circ
               \rho_{\beta\alpha}(\xi)\circ g_\alpha^{-1},\\
               &   & \forall\xi\in T(U_\alpha,U_\beta)\\
\end{array}  \right\}.\leqno (3.2.2)
$$
Note that $G_{\{\rho_{\beta\alpha}\}}$ is a group under
$\{g_\alpha\}\{h_\alpha\}=\{g_\alpha h_\alpha\}$. We give a
topology on $\O_{\{\rho_{\beta\alpha}\}}$ as follows. For each
index $\alpha$, let $\underline{K_\alpha}$ be a finite set
of compact subsets of $\widehat{U_\alpha}$, denoted by
$\{K_{\alpha,i}|i\in I_\alpha\}$, and let $\underline{O_\alpha}$ be a
finite set of open subsets of $\widehat{U_\alpha^\prime}$, denoted by
$\{O_{\alpha,i}|i\in I_\alpha\}$.
We define
$$
\O_{\{\rho_{\beta\alpha}\}}(\{\underline{K_\alpha}\},
\{\underline{O_\alpha}\})=\{\{f_\alpha\}\in\O_{\{\rho_{\beta\alpha}\}}\mid
f_\alpha(K_{\alpha,i})\subset O_{\alpha,i}, \forall i\in I_\alpha\}.
\leqno (3.2.3)
$$
Note that $\bigcap_j\O_{\{\rho_{\beta\alpha}\}}
(\{\underline{K_\alpha^{(j)}}\},\{\underline{O_\alpha^{(j)}}\})=
\O_{\{\rho_{\beta\alpha}\}}(\{\bigcup_j\underline{K_\alpha^{(j)}}\},
\{\bigcup_j\underline{O_\alpha^{(j)}}\})$. Hence the set of $(3.2.3)$ forms
a base of a topology on $\O_{\{\rho_{\beta\alpha}\}}$. The group
$G_{\{\rho_{\beta\alpha}\}}$ acts on $\O_{\{\rho_{\beta\alpha}\}}$
continuously by the formula $\{g_\alpha\}\cdot \{f_\alpha\}=
\{g_\alpha\circ f_\alpha\}$, whose orbit consists of mutually
conjugate homomorphisms. We denote the space of orbits by
$[\O_{\{\rho_{\beta\alpha}\}}]$, which is given with the quotient
topology. By Lemma 3.1.2, each $[\O_{\{\rho_{\beta\alpha}\}}]$ can
be naturally regarded as a subset of $[X;X^\prime]$. The image of
$\O_{\{\rho_{\beta\alpha}\}}(\{\underline{K_\alpha}\},
\{\underline{O_\alpha}\})$ in $[\O_{\{\rho_{\beta\alpha}\}}]$ is
denoted by $[\O_{\{\rho_{\beta\alpha}\}}(\{\underline{K_\alpha}\},
\{\underline{O_\alpha}\})]$. Note that the set of all
$[\O_{\{\rho_{\beta\alpha}\}}(\{\underline{K_\alpha}\},
\{\underline{O_\alpha}\})]$ forms a base of the quotient topology
on $[\O_{\{\rho_{\beta\alpha}\}}]$.

\begin{lem}
Let $\tau\in\O_{\{\rho_{ba}\}}$ and $\sigma\in
\O_{\{\rho_{\beta\alpha}\}}$ be any pair of equivalent homomorphisms.

\begin{itemize}
\item To any representative $\bar{\gamma}$ of an element $\gamma
\in\Gamma_{\sigma\tau}$, there is associated a local homeomorphism
$\phi_{\bar{\gamma}}$ from an open neighborhood of $\tau$ in
$\O_{\{\rho_{ba}\}}$ onto an open neighborhood of $\sigma$ in
$\O_{\{\rho_{\beta\alpha}\}}$ such that $\phi_{\bar{\gamma}}(\tau)=
\sigma$.
\item The germ of $\phi_{\bar{\gamma}}$ at $\tau$,
denoted by $\phi_\gamma$, depends only on the equivalence class
$\gamma\in\Gamma_{\sigma\tau}$, and for any $\gamma\in\Gamma_{\sigma\tau}$
and $\gamma^\prime\in\Gamma_{\tau\kappa}$,
$$
\phi_{\gamma}\circ\phi_{\gamma^\prime}=\phi_{\gamma\circ\gamma^\prime}.
\leqno (3.2.4)
$$
Moreover, if $\sigma=\phi_{\bar{\gamma}}(\tau)$ is also in
$\O_{\{\rho_{ba}\}}$, then $\Gamma_{\sigma\tau}\subset G_{\{\rho_{ba}\}}$
and $\phi_\gamma$ is induced by the action of $G_{\{\rho_{ba}\}}$
on $\O_{\{\rho_{ba}\}}$.
\item For any $\tau^\prime\in\mbox{Domain }(\phi_{\bar{\gamma}})$,
$\sigma^\prime=\phi_{\bar{\gamma}}(\tau^\prime)$ is equivalent to
$\tau^\prime$, and $\gamma\in\Gamma_{\sigma\tau}$ is canonically
associated with a $\gamma^\prime\in\Gamma_{\sigma^\prime\tau^\prime}$
such that the germ of $\phi_{\bar{\gamma}}$ at $\tau^\prime$ equals
$\phi_{\gamma^\prime}$.
\end{itemize}
\end{lem}

The proof of Lemma 3.2.2 is postponed to the end of this subsection.

\vspace{1.5mm}

\noindent{\bf Proof of Theorem 1.3:}

\vspace{1.5mm}

We give $[X;X^\prime]$ a topology which is generated by the set of
$[\O_{\{\rho_{\beta\alpha}\}}(\{\underline{K_\alpha}\},
\{\underline{O_\alpha}\})]$ for all possible data $\{\rho_{\beta\alpha}\}$,
$\{\underline{K_\alpha}\}$, and $\{\underline{O_\alpha}\}$. By
the existence of local homeomorphisms $\phi_{\bar{\gamma}}$, each
$[\O_{\{\rho_{\beta\alpha}\}}]$ with the quotient topology is an
open subset of $[X;X^\prime]$. Thus we obtained an open cover
$\{[\O_{\{\rho_{\beta\alpha}\}}]\}$ of $[X;X^\prime]$, where for each
$[\O_{\{\rho_{\beta\alpha}\}}]$, there is a space $\O_{\{\rho_{\beta\alpha}\}}$
and a discrete group $G_{\{\rho_{\beta\alpha}\}}$ acting
continuously on $\O_{\{\rho_{\beta\alpha}\}}$ such that
$[\O_{\{\rho_{\beta\alpha}\}}]=\O_{\{\rho_{\beta\alpha}\}}/
G_{\{\rho_{\beta\alpha}\}}$.

Consider the groupoid $\Gamma=\bigcup_{(\sigma,\tau)}\Gamma_{\sigma\tau}$,
where $\sigma,\tau$ are running over the set of all admissible
homomorphisms. The disjoint union $U=\bigsqcup\O_{\{\rho_{\beta\alpha}\}}$
is embedded into $\Gamma$ as the space of units. The mappings $\alpha,\omega:
\Gamma\rightarrow U$ are given by $\alpha(\gamma)=\tau$, $\omega(\gamma)=\sigma$
for any $\gamma\in\Gamma_{\sigma\tau}$. Clearly the restriction of $\Gamma$
to each $\O_{\{\rho_{\beta\alpha}\}}$ is the product groupoid
$G_{\{\rho_{\beta\alpha}\}}\times\O_{\{\rho_{\beta\alpha}\}}$
(cf. Lemma 3.1.2), and the space of $\Gamma$-orbits is isomorphic to
$[X;X^\prime]$. It remains to put a topology on $\Gamma$ so that it becomes
an \'{e}tale topological groupoid, such that the induced topology on
$U=\bigsqcup\O_{\{\rho_{\beta\alpha}\}}$ is defined by $(3.2.3)$. By
Lemma 3.2.2, $\{\phi_{\bar{\gamma}}\}$ generates a pseudogroup acting on
$U=\bigsqcup\O_{\{\rho_{\beta\alpha}\}}$ by local homeomorphisms,
where $U$ is given the topology defined by $(3.2.3)$, and the corresponding
space of germs of the elements in the pseudogroup is naturally
identified with $\Gamma$. Thus $\Gamma$, when equipped with the \'{e}tale
topology, is an \'{e}tale topological groupoid under which $[X;X^\prime]$
canonically becomes an orbispace.

\hfill $\Box$

\begin{prop}
Let $X,Y$ be paracompact, locally compact and Hausdorff, and $Z$ be any
orbispace. Then the mapping $(\Phi,\Psi)\mapsto \Psi\circ\Phi$,
$\forall \Phi\in [X;Y],\Psi\in [Y;Z]$, defines a continuous
map between the underlying spaces of $[X;Y]\times [Y;Z]$ and $[X;Z]$.
\end{prop}

\pf
We may represent $\Phi,\Psi$ by admissible homomorphisms $\sigma
=(\{f_{\alpha}\},\{\rho_{\beta\alpha}\})$, $\tau=(\{g_a\},
\{\eta_{ba}\})$, such that $\Psi\circ\Phi$ is represented by
$\sigma\circ\tau=(\{h_a\},\{\delta_{ba}\})$, which is also admissible.
Here $h_a=f_{\theta(a)}\circ g_a$, $\delta_{ba}=
\rho_{\theta(b)\theta(a)}\circ\eta_{ba}$ for some mapping of indices
$\theta:a\mapsto\alpha$ (cf. proof of Theorem 2.2.7).
Let $\O_{\{\delta_{ba}\}}(\{\underline{K_a}\},
\{\underline{O_a}\})$ be any given neighborhood of
$\sigma\circ\tau$, where $\underline{K_a}=\{K_{a,i}\mid i\in I_a\}$
and $\underline{O_a}=\{O_{a,i}\mid i\in I_a\}$. We set
$\underline{O_\alpha}=\emptyset$ if $\alpha\neq\theta(a)$
for any index $a$, and set $\underline{O_\alpha}=\{O_{a,i}\mid i\in I_a\}$
if $\alpha=\theta(a)$ for some index $a$. Since $Y$ is locally
compact and Hausdorff, we may choose $Q_{a,i}$, $i\in I_a$, such
that $g_a(K_{a,i})\subset Q_{a,i}$, and the closure $\overline{Q_{a,i}}
=L_{a,i}$ is compact, satisfying $f_{\theta(a)}(L_{a,i})\subset
O_{a,i}$. Now set $\underline{L_\alpha}=\emptyset$ if $\alpha\neq\theta(a)$
for any index $a$, $\underline{L_\alpha}=\{L_{a,i}\mid i\in I_a\}$
if $\alpha=\theta(a)$ for some index $a$, and $\underline{Q_a}
=\{Q_{a,i}\mid i\in I_a\}$. Then $\O_{\{\rho_{\beta\alpha}\}}
(\{\underline{L_\alpha}\},\{\underline{O_\alpha}\})$ and $\O_{\{\eta_{ba}\}}
(\{\underline{K_a}\},\{\underline{Q_a}\})$ are neighborhoods of
$\sigma$ and $\tau$ respectively, such that for any $\sigma^\prime\in
\O_{\{\rho_{\beta\alpha}\}}(\{\underline{L_\alpha}\},
\{\underline{O_\alpha}\})$, $\tau^\prime\in\O_{\{\eta_{ba}\}}(\{
\underline{K_a}\},\{\underline{Q_a}\})$, $\sigma^\prime\circ
\tau^\prime$ is defined and lies in
$\O_{\{\delta_{ba}\}}(\{\underline{K_a}\},\{\underline{O_a}\})$.
Thus the mapping $(\Phi,\Psi)\mapsto\Psi\circ\Phi$
is continuous.

\hfill $\Box$

\noindent{\bf Proof of Lemma 3.2.2}

\vspace{1.5mm}

Let $\tau=\{f_a\}$, $\sigma=\{f_\alpha\}$.
We shall first construct the local homeomorphism $\phi_{\bar{\gamma}}$
for the special circumstance where $\tau$ is induced by $\sigma$,
and the following condition
$$
U^\prime_a\subset U_\alpha^\prime \mbox{ whenever } U_a\subset
U_\alpha \leqno (3.2.5)
$$
is satisfied. Note that $(3.2.5)$ can be always arranged. This is
because $\{U_\alpha\}$ is locally finite, therefore for each $U_a$,
there are only finitely many $U_\alpha$'s containing $U_a$, so that
one may replace $U_a^\prime$ by a smaller one $U_a^\prime\cap
U_\alpha^\prime$ whenever $U_a\subset U_\alpha$ occurs, and then
work with the homomorphism induced by $\tau$.

On the other hand, in light of Lemma 3.2.1, we may assume, by passing to
an induced homomorphism of $\tau$, that for each index $a$, there is an
$\alpha$ such that $\overline{U_a}\subset U_{\alpha}$. Moreover, there exists
an open cover $\{V_a\}$ of $X$, which will be fixed throughout, such that
$\overline{V_a}\subset U_a$. We set $K_a=\pi_{U_a}^{-1}(\overline{V_a})$.

Now given any representative $\bar{\gamma}=(\theta,\{\xi_a\},
\{\xi_a^\prime\})$ of $\gamma\in\Gamma_{\sigma\tau}$, where
$\theta:a\mapsto\alpha$ is a mapping of indexes such that
$$
\overline{U_a}\subset U_{\theta(a)}, \; \forall a,\leqno (3.2.6)
$$
and $\xi_a\in T(U_a,U_{\theta(a)})$, $\xi_a^\prime\in
T(U_a^\prime,U_{\theta(a)}^\prime)$, with
$$
f_a=(\phi_{\xi_{a}^\prime})^{-1}\circ f_{\theta(a)}\circ\phi_{\xi_{a}},
\;\forall a, \leqno (3.2.7\;a)
$$
and for any $\eta\in T(U_a,U_b)$,
$$
\rho_{ba}(\eta)=(\xi_{b}^\prime)^{-1}\circ\rho_{\theta(b)\theta(a)}
(\xi_{b}\circ\eta\circ\xi_{a}^{-1})\circ\xi_{a}^\prime
(x)\;\forall x\in f_a(\mbox{Domain }(\phi_\eta)),\leqno (3.2.7\;b)
$$
we shall construct an open embedding $\phi_{\bar{\gamma}}$ from an
open neighborhood $\O_{\{\rho_{ba}\}}(\{\underline{L_a^0}\},
\{\underline{Q_a^0}\})$ of $\{f_a\}$ in $\O_{\{\rho_{ba}\}}$
onto an open neighborhood $\O_{\{\rho_{\beta\alpha}\}}
(\{\underline{K_\alpha^0}\},\{\underline{O_\alpha^0}\})$ of
$\{f_\alpha\}$ in $\O_{\{\rho_{\beta\alpha}\}}$, such that
$\phi_{\bar{\gamma}}(\{f_a\})=\{f_\alpha\}$, and for
any $\{f_a^\prime\}\in\O_{\{\rho_{ba}\}}(\{\underline{L_a^0}\},
\{\underline{Q_a^0}\})$, $\{f_a^\prime\}$ is induced by
$\{f_\alpha^\prime\}=\phi_{\bar{\gamma}}(\{f_a^\prime\})$ via
the data $\bar{\gamma}=(\theta,\{\xi_a\},\{\xi_a^\prime\})$.

We begin by defining a mapping from $T(U_a,U_\alpha)$ to
$T(U_a^\prime,U_\alpha^\prime)$ for each pair $(a,\alpha)$
with $U_a\cap U_\alpha\neq\emptyset$, which is denoted by
$\xi\mapsto\xi^\prime$, by the formula
$$
\xi^\prime=\rho_{\alpha\theta(a)}(\xi\circ\xi_a^{-1})\circ
\xi_a^\prime(x),\; \forall x\in f_a(\mbox{Domain }(\phi_\xi)).
\leqno (3.2.8)
$$
Note that $\xi_a\mapsto \xi_a^\prime$, $\forall a$, under
$\xi\mapsto \xi^\prime$ defined by $(3.2.8)$. Moreover, since
$U_a^\prime\subset U^\prime_{\theta(a)}$ by $(3.2.5)$, one has
$$
\mbox{Domain }(\phi_{\xi^\prime})\subset
\mbox{Domain }(\phi_{\rho_{\alpha\theta(a)}(\xi\circ\xi_a^{-1})}
\circ\phi_{\xi_a^\prime}). \leqno (3.2.9)
$$

We determine the data $\{\underline{L_a^0}\},\{\underline{Q_a^0}\}$
as follows. Let index $a$ be given. (1) Suppose $U_a\cap U_\alpha\neq
\emptyset$ for some index $\alpha$. Since $f_a(\overline{\pi_{U_a}^{-1}
(U_a\cap U_\alpha)})\subset \widehat{U_a^\prime}$ is compact, there are
only finitely many components $\{Q_{a,\alpha}^i\mid i\in\Lambda(\alpha)\}$
of $\pi_{U_a^\prime}^{-1}(U_a^\prime\cap U_\alpha^\prime)$ which
contain $f_a(\overline{\pi_{U_a}^{-1}(U_a\cap U_\alpha)})$. We let
$L_{a,\alpha}^i\subset \widehat{\overline{U_a}}$ be the union of
$\overline{\mbox{Domain }(\phi_\xi)}$ for all $\xi\in T(U_a,U_\alpha)$
such that $\mbox{Domain }(\phi_{\xi^\prime})=Q_{a,\alpha}^i$. We set
$\underline{Q_a^{(1)}}=\{Q_{a,\alpha}^i\mid i\in\Lambda(\alpha),
U_a\cap U_\alpha\neq\emptyset\}$, and set $\underline{L_a^{(1)}}
=\{L_{a,\alpha}^i\cap K_a\mid i\in\Lambda(\alpha), U_a\cap U_\alpha\neq
\emptyset\}$. Since $\{U_\alpha\}$ is locally finite, both
$\underline{L_a^{(1)}},\underline{Q_a^{(1)}}$ are finite sets.
(2) Suppose $U_a\cap U_\alpha\cap U_b\neq\emptyset$ for some indexes
$\alpha,b$. Since $f_{\theta(a)}(\overline{\pi_{U_{\theta(a)}}^{-1}
(U_{\theta(a)}\cap U_\alpha\cap U_{\theta(b)})})$ is a compact subset
in $\pi_{U_{\theta(a)}^\prime}^{-1}(U_{\theta(a)}^\prime\cap U_\alpha^\prime
\cap U_{\theta(b)}^\prime)$, it is contained in finitely many components
$\{W_{a,(\alpha,b)}^i\mid i\in\Lambda(\alpha,b)\}$. We set
$Q_{a,(\alpha,b)}^i=(\phi_{\xi_a^\prime})^{-1}(W_{a,(\alpha,b)}^i)$,
$i\in\Lambda(\alpha,b)$, and let $L_{a,(\alpha,b)}^i$ be the union of
the closure of ${\bf a}\in\Lambda(\zeta_a,\zeta_b^{-1})$ such that
$f_{\theta(a)}\circ\phi_{\xi_a}({\bf a})\subset W_{a,(\alpha,b)}^i$,
where $\zeta_a\in T(U_a,U_\alpha)$, $\zeta_b\in T(U_b,U_\alpha)$. We
define $\underline{Q_a^{(2)}}=\{Q_{a,(\alpha,b)}^i\mid i\in\Lambda(\alpha,b),
U_a\cap U_\alpha\cap U_b\neq\emptyset\}$, $\underline{L_a^{(2)}}=
\{L_{a,(\alpha,b)}^i\cap K_a\mid i\in\Lambda(\alpha,b),U_a\cap U_\alpha\cap
U_b\neq\emptyset\}$. Both $\underline{L_a^{(2)}},\underline{Q_a^{(2)}}$ are
finite sets because $\{U_\alpha\},\{U_b\}$ are locally finite. (3) Suppose
$U_a\cap U_\alpha\cap U_\beta\neq\emptyset$ for some indexes $\alpha,\beta$.
Since $f_{\theta(a)}(\overline{\pi_{U_{\theta(a)}}^{-1}(U_{\theta(a)}\cap
U_\alpha\cap U_\beta)})$ is a compact subset in $\pi_{U_{\theta(a)}^\prime}
^{-1}(U_{\theta(a)}^\prime\cap U_\alpha^\prime\cap U_\beta^\prime)$, it is
contained in finitely many components $\{W_{a,(\alpha,\beta)}^i\mid i\in
\Lambda(\alpha,\beta)\}$. We set $Q_{a,(\alpha,\beta)}^i=
(\phi_{\xi_a^\prime})^{-1}(W_{a,(\alpha,\beta)}^i)$, $i\in\Lambda(\alpha,
\beta)$, and let $L_{a,(\alpha,\beta)}^i$ be the union of the closure of
${\bf a}\in\Lambda(\xi,\eta)$ such that $f_{\theta(a)}\circ\phi_{\xi_a}
({\bf a})\subset W_{a,(\alpha,\beta)}^i$, where $\xi\in T(U_a,U_\alpha)$,
$\eta\in T(U_\alpha, U_\beta)$. We define $\underline{Q_a^{(3)}}=
\{Q_{a,(\alpha,\beta)}^i\mid i\in\Lambda(\alpha,\beta), U_a\cap U_\alpha
\cap U_\beta\neq\emptyset\}$, $\underline{L_a^{(3)}}=\{L_{a,(\alpha,\beta)}^i
\cap K_a\mid i\in\Lambda(\alpha,\beta),U_a\cap U_\alpha\cap U_\beta\neq
\emptyset\}$. Both $\underline{L_a^{(3)}},\underline{Q_a^{(3)}}$ are finite
sets because $\{U_\alpha\}$ is locally finite. (4) Suppose $U_a\cap U_\alpha
\cap U_\beta\cap U_\gamma\neq\emptyset$ for some indexes $\alpha,\beta,\gamma$.
Since $f_a(\overline{\pi_{U_a}^{-1}(U_a\cap U_\alpha\cap U_\beta\cap
U_\gamma)})$ is a compact subset in $\pi_{U_a^\prime}^{-1}(U_a^\prime\cap
U_\alpha^\prime\cap U_\beta^\prime\cap U_\gamma^\prime)$, it is contained
in finitely many components $\{Q_{a,(\alpha,\beta,\gamma)}^i\mid i\in
I(\alpha,\beta,\gamma)\}$. We define $L_{a,(\alpha,\beta,\gamma)}^i=\{x\in
\overline{\pi_{U_a}^{-1}(U_a\cap U_\alpha\cap U_\beta\cap U_\gamma)}
\cap K_a\mid f_a(x)\in Q_{a,(\alpha,\beta,\gamma)}^i\}$, and set
$\underline{L_a^{(4)}}=\{L_{a,(\alpha,\beta,\gamma)}^i\mid i\in I(\alpha,
\beta,\gamma),U_a\cap U_\alpha\cap U_\beta\cap U_\gamma\neq\emptyset\}$,
$\underline{Q_a^{(4)}}=\{Q_{a,(\alpha,\beta,\gamma)}^i\mid i\in I(\alpha,
\beta,\gamma),U_a\cap U_\alpha\cap U_\beta\cap U_\gamma\neq\emptyset\}$.
Then both $\underline{L_a^{(4)}}$, $\underline{Q_a^{(4)}}$ are finite
sets because $\{U_\alpha\}$ is locally finite. We define
$\underline{L_a^0}=\cup_{k=1}^4\underline{L_a^{(k)}}$,
$\underline{Q_a^0}=\cup_{k=1}^4\underline{Q_a^{(k)}}$. Clearly $\{f_a\}$
is contained in $\O_{\{\rho_{ba}\}}(\{\underline{L_a^0}\},
\{\underline{Q_a^0}\})$.

We define the map $\phi_{\bar{\gamma}}:\O_{\{\rho_{ba}\}}
(\{\underline{L_a^0}\},\{\underline{Q_a^0}\})\rightarrow
\O_{\{\rho_{\beta\alpha}\}}$. Given any $\{f_a^\prime\}\in\O_{\{\rho_{ba}\}}
(\{\underline{L_a^0}\},\{\underline{Q_a^0}\})$, $\{f_\alpha^\prime\}=
\phi_{\bar{\gamma}}(\{f_a^\prime\})$ is constructed as follows. For any
$x\in\widehat{\overline{U_\alpha}}$, we pick a $U_a$ and a $\xi\in T(U_a,
U_\alpha)$ such that $x\in\overline{\mbox{Range }(\phi_\xi)}$ and
$(\phi_\xi)^{-1}(x)\in K_a$. Then by the assumption that $\{f_a^\prime\}
\in\O_{\{\rho_{ba}\}}(\{\underline{L_a^0}\},\{\underline{Q_a^0}\})$, we
have $f_a^\prime\circ (\phi_\xi)^{-1}(x)\subset\mbox{Domain }
(\phi_{\xi^\prime})$. We define
$$
f_\alpha^\prime(x)=\phi_{\xi^\prime}\circ f_a^\prime\circ
(\phi_{\xi})^{-1}(x). \leqno (3.2.10)
$$
It remains to show that (1) each $f_\alpha^\prime$
is well-defined, (2) $\phi_{\rho_{\beta\alpha}(\eta)}\circ
f_\alpha^\prime=f_\beta^\prime\circ\phi_\eta$ holds for
any $\alpha,\beta$ and $\eta\in T(U_\alpha,U_\beta)$, (3)
$\rho_{\gamma\alpha}(\eta\circ\zeta({\bf a}))=\rho_{\gamma\beta}(\eta)
\circ\rho_{\beta\alpha}(\zeta)(\underline{\{f_\alpha^\prime\}}({\bf a}))$
holds for any indexes $\alpha,\beta,\gamma$, and for any $\zeta\in
T(U_\alpha,U_\beta)$, $\eta\in T(U_\beta,U_\gamma)$, and ${\bf a}\in
\Lambda(\zeta,\eta)$, (4) $\{f_a^\prime\}$ is induced by
$\{f_\alpha^\prime\}$ via the data $\bar{\gamma}=(\theta,\{\xi_a\},
\{\xi_a^\prime\})$, and (5) $\phi_{\bar{\gamma}}(\{f_a\})=\{f_\alpha\}$.

The problem can be easily reduced to the following two types of
identities. The first one is that for any $\xi\in T(U_a,U_\alpha)$,
$$
\phi_{\rho_{\alpha\theta(a)}(\xi\circ\xi_a^{-1})}\circ
(\phi_{\xi_a^\prime}\circ f_a^\prime\circ (\phi_{\xi_a})^{-1})
=(\phi_{\xi^\prime}\circ f_a^\prime\circ (\phi_{\xi})^{-1})\circ
\phi_{\xi\circ\xi_a^{-1}}
\leqno (3.2.11)
$$
on $\phi_{\xi_a}(K_a\cap\mbox{Domain }(\phi_\xi))$, which follows easily
from $(3.2.9)$. The second one is that for any $\xi\in T(U_a,U_b)$,
$$
\phi_{\rho_{\theta(b)\theta(a)}(\xi_b\circ\xi\circ\xi_a^{-1})}\circ
(\phi_{\xi_a^\prime}\circ f_a^\prime\circ (\phi_{\xi_a})^{-1})=
(\phi_{\xi_b^\prime}\circ f_b^\prime\circ (\phi_{\xi_b})^{-1})
\circ\phi_{\xi_b\circ\xi\circ\xi_a^{-1}} \leqno (3.2.12)
$$
on $\phi_{\xi_a}(\mbox{Domain }(\phi_\xi))$, which follows from the
fact that $\phi_{\rho_{ba}(\xi)}\circ f_a^\prime=f_b^\prime\circ\phi_{\xi}$
and that $\mbox{Domain }(\phi_{\rho_{ba}(\xi)})$ is contained in
$\mbox{Domain }((\phi_{\xi_b^\prime})^{-1}\circ
\phi_{\rho_{\theta(b)\theta(a)}(\xi_b\circ\xi\circ\xi_a^{-1})}
\circ\phi_{\xi_a^\prime})$ because of $(3.2.5)$.

To see that each $f_\alpha^\prime$ is well-defined, suppose a point
$x\in\widehat{\overline{U_\alpha}}$ lies in $\overline{\mbox{Range }
(\phi_{\zeta_a})}\cap\overline{\mbox{Range }(\phi_{\zeta_b})}$ for some
$\zeta_a\in T(U_a,U_\alpha)$, $\zeta_b\in T(U_b,U_\alpha)$, and
$(\phi_{\zeta_a})^{-1}(x)\in K_a$, $(\phi_{\zeta_b})^{-1}(x)\in K_b$.
We need to show that $\phi_{\zeta_a^\prime}\circ f_a^\prime\circ
(\phi_{\zeta_a})^{-1}(x)=\phi_{\zeta_b^\prime}\circ f_b^\prime
\circ (\phi_{\zeta_b})^{-1}(x)$, which can be derived as follows.
For the convenience of expression, we set $f_{\alpha,a}^\prime
=\phi_{\zeta_a^\prime}\circ f_a^\prime\circ (\phi_{\zeta_a})^{-1}$,
$f_{\alpha,b}^\prime=\phi_{\zeta_b^\prime}\circ f_b^\prime\circ
(\phi_{\zeta_b})^{-1}$, $f_{\theta(a)}^\prime=\phi_{\xi_a^\prime}
\circ f_a^\prime\circ (\phi_{\xi_a})^{-1}$ and $f_{\theta(b)}^\prime=
\phi_{\xi_b^\prime}\circ f_b^\prime\circ (\phi_{\xi_b})^{-1}$.
Let ${\bf a}_x\in\Lambda(\zeta_a,\zeta_b^{-1})$ be an element whose
closure contains $(\phi_{\zeta_a})^{-1}(x)$. Suppose ${\bf a}_x\subset
L_{a,(\alpha,b)}^i$ for some $i\in\lambda(\alpha,b)$. Then
$$
\phi_{\rho_{\alpha\theta(a)}(\zeta_a\circ\xi_a^{-1})}
=\phi_{\rho_{\alpha\theta(b)}(\zeta_b\circ\xi_b^{-1})}\circ
\phi_{\rho_{\theta(b)\theta(a)}(\xi_b\circ(\zeta_b^{-1}\circ
\zeta_a({\bf a}_x))\circ\xi_a^{-1})}\mbox { on }
W_{a,(\alpha,b)}^i. \leqno (3.2.13)
$$
By the assumption that $\{f_a^\prime\}\in\O_{\{\rho_{ba}\}}
(\{\underline{L_a^0}\},\{\underline{Q_a^0}\})$, we have
$f_{\theta(a)}^\prime(\phi_{\xi_a\circ\zeta_a^{-1}}(x))\in
W_{a,(\alpha,b)}^i$. Thus by $(3.2.11)$, $(3.2.13)$ and then $(3.2.12)$,
\begin{eqnarray*}
 &   & f_{\alpha,a}^\prime(x)=\phi_{\rho_{\alpha\theta(a)}
 (\zeta_a\circ\xi_a^{-1})}\circ f_{\theta(a)}^\prime\circ
 \phi_{\xi_a\circ\zeta_a^{-1}}(x)\\
 & = & \phi_{\rho_{\alpha\theta(b)}(\zeta_b\circ\xi_b^{-1})}
 \circ\phi_{\rho_{\theta(b)\theta(a)}(\xi_b\circ(\zeta_b^{-1}
 \circ\zeta_a({\bf a}_x))\circ\xi_a^{-1})}\circ f_{\theta(a)}^\prime
 \circ\phi_{\xi_a\circ\zeta_a^{-1}}(x)\\
 & = & \phi_{\rho_{\alpha\theta(b)}(\zeta_b\circ\xi_b^{-1})}\circ
 f_{\theta(b)}^\prime\circ\phi_{\xi_b\circ (\zeta_b^{-1}\circ
 \zeta_a({\bf a}_x))\circ\xi_a^{-1}}\circ
\phi_{\xi_a\circ\zeta_a^{-1}}(x)=f_{\alpha,b}^\prime(x).
\end{eqnarray*}

To verify $\phi_{\rho_{\beta\alpha}(\eta)}\circ f_\alpha^\prime
=f_\beta^\prime\circ\phi_\eta$, let $x$ be any point in
$\mbox{Domain }(\phi_\eta)$. We take a $\xi\in T(U_a,U_\alpha)$
for some index $a$ such that $x\in\mbox{Range }(\phi_\xi)$ and
$(\phi_\xi)^{-1}(x)\in K_a$. Set ${\bf a}_x\in\Lambda(\xi,\eta)$
to be the element containing $(\phi_{\xi})^{-1}(x)$, and suppose
${\bf a}_x\subset L_{a,(\alpha,\beta)}^i$ for some
$i\in\Lambda(\alpha,\beta)$. Then on $W_{a,(\alpha,\beta)}^i$
we have $\phi_{\rho_{\beta\theta(a)}(\eta\circ\xi({\bf a}_x)
\circ\xi_a^{-1})}=\phi_{\rho_{\beta\alpha}(\eta)}\circ
\phi_{\rho_{\alpha\theta(a)}(\xi\circ\xi_a^{-1})}$. On the other
hand, the assumption $\{f_a^\prime\}\in\O_{\{\rho_{ba}\}}
(\{\underline{L_a^0}\},\{\underline{Q_a^0}\})$ implies that
$f_{\theta(a)}^\prime\circ\phi_{\xi_a\circ\xi^{-1}}(x)\in
W_{a,(\alpha,\beta)}^i$. Hence by $(3.2.11)$,
we have
\begin{eqnarray*}
&   & \phi_{\rho_{\beta\alpha}(\eta)}\circ f_\alpha^\prime(x)
= \phi_{\rho_{\beta\alpha}(\eta)}\circ\phi_{\rho_{\alpha\theta(a)}
(\xi\circ\xi_a^{-1})}\circ f_{\theta(a)}^\prime\circ\phi_{\xi_a
\circ\xi^{-1}}(x)\\
& = & \phi_{\rho_{\beta\theta(a)}(\eta\circ\xi({\bf a}_x)\circ\xi_a^{-1})}
\circ f_{\theta(a)}^\prime\circ\phi_{\xi_a\circ\xi^{-1}}(x)=
 f_\beta^\prime\circ\phi_{\eta\circ\xi({\bf a}_x)\circ\xi_a^{-1}}
\circ\phi_{\xi_a\circ\xi^{-1}}(x)\\
& = & f_\beta^\prime\circ\phi_{\eta}(x), \forall x\in\mbox{Domain }
(\phi_\eta).
\end{eqnarray*}

To verify $\rho_{\gamma\alpha}(\eta\circ\zeta({\bf a}))=
\rho_{\gamma\beta}(\eta)\circ\rho_{\beta\alpha}(\zeta)
(\underline{\{f_\alpha^\prime\}}({\bf a}))$, it suffices to
show that $f_\alpha^\prime(x)\in\underline{\{f_\alpha\}}({\bf a})$
if $x\in {\bf a}$. To this end, we pick a $U_a$ and a $\xi\in
T(U_a,U_\alpha)$ such that $x\in\mbox{Range }(\phi_\xi)$ and
$(\phi_{\xi})^{-1}(x)\in K_a$. Now observe that $\phi_{\xi^\prime}
\circ f_a\circ (\phi_{\xi})^{-1}(x)=f_\alpha(x)\in\underline{\{f_\alpha\}}
({\bf a})$, so that $f_a\circ(\phi_{\xi})^{-1}(x)\in
(\phi_{\xi^\prime})^{-1}(\underline{\{f_\alpha\}}({\bf a}))$.
Let $(\phi_{\xi^\prime})^{-1}(\underline{\{f_\alpha\}}({\bf a}))
=Q_{a,(\alpha,\beta,\gamma)}^i$ for some $i\in I(\alpha,\beta,\gamma)$.
Then $(\phi_{\xi})^{-1}(x)$ lies in $L_{a,(\alpha,\beta,\gamma)}^i$
by the virtue of definition. Now we appeal to the assumption that
$\{f_a^\prime\}\in\O_{\{\rho_{ba}\}}(\{\underline{L_a^0}\},
\{\underline{Q_a^0}\})$ to conclude that $f_a^\prime\circ
(\phi_{\xi})^{-1}(x)$ lies in $Q_{a,(\alpha,\beta,\gamma)}^i$,
which implies that $f_\alpha^\prime(x)=\phi_{\xi^\prime}\circ
f_a^\prime\circ (\phi_{\xi})^{-1}(x)\in\phi_{\xi^\prime}(Q_{a,(\alpha,
\beta,\gamma)}^i)\subset\underline{\{f_\alpha\}}({\bf a})$.

To verify that $\{f_a^\prime\}$ is induced by $\{f_\alpha^\prime\}=
\phi_{\bar{\gamma}}(\{f_a^\prime\})$ via the data $\bar{\gamma}=
(\theta,\{\xi_a\},\{\xi_a^\prime\})$, we observe that: (1) For $(3.2.7\;a)$,
$f_a^\prime(x)=(\phi_{\xi_a^\prime})^{-1}\circ f_{\theta(a)}^\prime
\circ\phi_{\xi_a}(x)$ holds for any $x\in K_a$ by the nature of construction;
for the case when $x\in\widehat{U_a}\setminus K_a$, we pick a $U_b$
such that $x=\phi_{\xi}(y)$ for some $y\in K_b$ and $\xi\in T(U_b,U_a)$,
and then apply $(3.2.11)$ and $(3.2.12)$. Hence $(3.2.7\;a)$ holds for
$\{f_a^\prime\}$, $\{f_\alpha^\prime\}$. (2) For $(3.2.7\;b)$, we need to
verify that $f_a^\prime(\mbox{Domain }(\phi_\eta))$ is contained in
${\bf a}_\eta^{\bar{\gamma}}$, which is the unique element of
$\Lambda(\xi_a^\prime,\rho_{\theta(b)\theta(a)}(\xi_b\circ\eta\circ
\xi_a^{-1}),(\xi_b^\prime)^{-1})$ that contains
$f_a(\mbox{Domain }(\phi_\eta))$. But this follows from the fact that
$\mbox{Domain }(\phi_{\rho_{ba}(\eta)})$ is contained in
$\mbox{Domain }((\phi_{\xi_b^\prime})^{-1}\circ
\phi_{\rho_{\theta(b)\theta(a)}(\xi_b\circ\eta\circ\xi_a^{-1})}
\circ\phi_{\xi_a^\prime})$ because of $(3.2.5)$, so that $\mbox{Domain }
(\phi_{\rho_{ba}(\eta)})={\bf a}_\eta^{\bar{\gamma}}$, and the fact that
$\phi_{\rho_{ba}(\eta)}\circ f_a^\prime=f_b^\prime\circ\phi_\eta$
so that $f_a^\prime(\mbox{Domain }(\phi_\eta))\subset\mbox{Domain }
(\phi_{\rho_{ba}(\eta)})$. Hence $\{f_a^\prime\}$ is induced by
$\{f_\alpha^\prime\}$ via the data $\bar{\gamma}=(\theta,\{\xi_a\},
\{\xi_a^\prime\})$. Note that as a corollary, $\phi_{\bar{\gamma}}$
is injective.

To verify $\phi_{\bar{\gamma}}(\{f_a\})=\{f_\alpha\}$, we observe
that if we let $\{\bar{f}_\alpha\}=\phi_{\bar{\gamma}}(\{f_a\})$,
then since both $\{f_\alpha\}$, $\{\bar{f}_\alpha\}$ induce $\{f_a\}$
via the data $\bar{\gamma}=(\theta,\{\xi_a\},\{\xi_a^\prime\})$,
by Lemma 3.1.2, $\{f_\alpha\}$, $\{\bar{f}_\alpha\}$ are conjugate
to each other via some $\{g_\alpha\}$, where from the proof of Lemma
3.1.2, we notice that $g_\alpha=g_\alpha(\zeta_a)=\zeta_{a,2}^\prime
\circ (\zeta_{a,1}^\prime)^{-1}(x), \forall x\in f_a(\widehat{U_a})$.
In the present case, one may take $\zeta_a=\xi_a$, and hence
$\zeta_{a,1}^\prime=\zeta_{a,2}^\prime=\xi_a^\prime$. So $g_\alpha=1$,
$\forall\alpha$. Thus $\{f_\alpha\}=\{\bar{f}_\alpha\}
=\phi_{\bar{\gamma}}(\{f_a\})$. (In this argument, we may accommodate
the hypothesis in the proof of Lemma 3.1.2 that each $U_\alpha=
\bigcup_{a\in I_\alpha}U_a$ for a subset $\{U_a\mid a\in I_\alpha\}$
of $\{U_a\}$ by passing to an induced homomorphism of $\{f_a\}$,
which may not be admissible.)

It remains to show that $\phi_{\bar{\gamma}}$ is a continuous, open map.

Given any open subset $\O_{\{\rho_{\beta\alpha}\}}
(\{\underline{K_\alpha}\},\{\underline{O_\alpha}\})$
which contains $\{f_{\alpha,0}^\prime\}=\phi_{\bar{\gamma}}
(\{f_{a,0}^\prime\})$, we shall construct an open neighborhood
$\O_{\{\rho_{ba}\}}(\{\underline{L_a}\},\{\underline{Q_a}\})$ of
$\{f_{a,0}^\prime\}$, such that the open neighborhood $\O_{\{\rho_{ba}\}}
(\{\underline{L_a}\cup\underline{L_a^0}\},\{\underline{Q_a}\cup
\underline{Q_a^0}\})$ of $\{f_{a,0}^\prime\}$ is mapped into
$\O_{(\{\rho_{\beta\alpha}\}}(\{\underline{K_\alpha}\},
\{\underline{O_\alpha}\})$ under $\phi_{\bar{\gamma}}$.
Write $\underline{K_\alpha}=\{K_{\alpha,i}\mid i\in I_\alpha\}$,
$\underline{O_\alpha}=\{O_{\alpha,i}\mid i\in I_\alpha\}$. Let
$a$ be any given index. Suppose $U_a\cap U_\alpha\neq\emptyset$
for some index $\alpha$. Then for each $i\in I_\alpha$,
$\overline{V_a}\cap\pi_{U_\alpha}(K_{\alpha,i})$,
where $\overline{V_a}$ is given in $K_a=\pi_{U_a}^{-1}
(\overline{V_a})$, is contained in finitely many components
$\{U_{a,(\alpha,i)}^j|j\in\Lambda(\alpha,i)\}$ of $U_a\cap U_\alpha$.
This is because $\overline{V_a}\cap\pi_{U_\alpha}(K_{\alpha,i})$
is a compact subset of $U_a\cap U_\alpha$. Moreover, both $\cup_j\pi_{U_a}^{-1}
(U_{a,(\alpha,i)}^j)$ and $\cup_j\pi_{U_\alpha}^{-1}(U_{a,(\alpha,i)}^j)$
have finitely many components because $U_a, U_\alpha$ are locally compact,
Hausdorff, and the maps $\pi_{U_a},\pi_{U_\alpha}$ are proper.
Hence there exist finitely many $\xi_k, k\in I_{(\alpha,i)}$,
such that if $x\in K_{\alpha,i}$ with $\pi_{U_\alpha}(x)\in\overline{V_a}
=\pi_{U_a}(K_a)$, then there is a $\xi_k$ satisfying $x\in
\mbox{Range }(\phi_{\xi_k})$ for some $k\in I_{(\alpha,i)}$.
Set $L_{a,(\alpha,i)}^k=(\phi_{\xi_k})^{-1}(K_{\alpha,i})\cap
K_a$, $Q_{a,(\alpha,i)}^k=(\phi_{\xi_k^\prime})^{-1}(O_{\alpha,i})$
for each $k\in I_{(\alpha,i)}$. We define $\underline{L_a}=
\{L_{a,(\alpha,i)}^k\mid k\in I_{(\alpha,i)}, i\in I_\alpha,U_a
\cap U_\alpha\neq\emptyset\}$ and $\underline{Q_a}=\{Q_{a,(\alpha,i)}^k
\mid k\in I_{(\alpha,i)},i\in I_\alpha,U_a\cap U_\alpha\neq\emptyset\}$.
Both $\underline{L_a},\underline{Q_a}$ are finite sets because
$\{U_\alpha\}$ is locally finite. Now the fact that
$\{f_{\alpha,0}^\prime\}=\phi_{\bar{\gamma}}(\{f_{a,0}^\prime\})$
is contained in $\O_{\{\rho_{\beta\alpha}\}}
(\{\underline{K_\alpha}\},\{\underline{O_\alpha}\})$ implies,
by the nature of construction, that $\{f_{a,0}^\prime\}$ is contained
in $\O_{\{\rho_{ba}\}}(\{\underline{L_a}\},\{\underline{Q_a}\})$.

Now suppose $\{f_a^\prime\}$ is in $\O_{\{\rho_{ba}\}}
(\{\underline{L_a}\cup\underline{L_a^0}\},\{\underline{Q_a}
\cup\underline{Q_a^0}\})$. Given any index $\alpha$, let $x\in
K_{\alpha,i}$ be any point. Since $\{V_a\}$ is a cover of $X$,
there is a $U_a$ with a $\xi_k\in T(U_a,U_\alpha)$
for some $k\in I_{(\alpha,i)}$, such that $x\in \mbox{Range }
(\phi_{\xi_k})$ and $(\phi_{\xi_k})^{-1}(x)\in K_a$. By definition
$(\phi_{\xi_k})^{-1}(x)\in L_{a,(\alpha,i)}^k$, hence $f_a^\prime
\circ (\phi_{\xi_k})^{-1}(x)\in Q_{a,(\alpha,i)}^k
=(\phi_{\xi_k^\prime})^{-1}(O_{\alpha,i})$, which implies that
$f_\alpha^\prime(x)=\phi_{\xi_k^\prime}\circ f_a^\prime\circ
(\phi_{\xi_k})^{-1}(x)\in O_{\alpha,i}$. Thus $\{f_\alpha^\prime\}
=\phi_{\bar{\gamma}}(\{f_a^\prime\})$ lies in
$\O_{\{\rho_{\beta\alpha}\}}(\{\underline{K_\alpha}\},
\{\underline{O_\alpha}\})$. In other words, the open neighborhood
$\O_{\{\rho_{ba}\}}(\{\underline{L_a}\cup\underline{L_a^0}\},
\{\underline{Q_a}\cup\underline{Q_a^0}\})$ of $\{f_{a,0}^\prime\}$
is mapped into the given open neighborhood $\O_{\{\rho_{\beta\alpha}\}}
(\{\underline{K_\alpha}\},\{\underline{O_\alpha}\})$ of
$\{f_{\alpha,0}^\prime\}$ under the map $\phi_{\bar{\gamma}}$,
hence $\phi_{\bar{\gamma}}$ is continuous.

To see that $\phi_{\bar{\gamma}}$ is an open map, we first show that
the image of $\phi_{\bar{\gamma}}$ is an open subset of
$\O_{\{\rho_{\beta\alpha}\}}$. To this end, suppose $\{f_\alpha^\prime\}
\in\O_{\{\rho_{\beta\alpha}\}}$ satisfies $f_{\theta(a)}^\prime
(\overline{\mbox{Range }(\phi_{\xi_a})})\subset\mbox{Range }
(\phi_{\xi_a^\prime})$ for each index $a$. Then an induced $\{\bar{f}_a\}$
of $\{f_\alpha^\prime\}$ is defined via $\bar{\gamma}=(\theta,\{\xi_a\},
\{\xi_a^\prime\})$. If furthermore, $\{\bar{f}_a\}$ lies in
$\O_{\{\rho_{ba}\}}(\{\underline{L_a^0}\},\{\underline{Q_a^0}\})$,
then $\{\bar{f}_\alpha\}=\phi_{\bar{\gamma}}(\{\bar{f}_a\})$ must
equal $\{f_\alpha^\prime\}$, since both $\{\bar{f}_\alpha\}$,
$\{f_\alpha^\prime\}$ induce $\{\bar{f}_a\}$ via the same data
$\bar{\gamma}=(\theta,\{\xi_a\},\{\xi_a^\prime\})$, cf. Lemma 3.1.2.

Note that by $(3.2.6)$, $\overline{\mbox{Range }(\phi_{\xi_a})}
\subset\widehat{U_{\theta(a)}}$ is compact. We set
$\underline{K_\alpha^0}=\{\overline{\mbox{Range }(
\phi_{\xi_a})},\phi_{\xi_a}(L)\mid L\in\underline{L_a^0}\}$,
$\underline{O_\alpha^0}=\{\mbox{Range }(\phi_{\xi_a^\prime}),
\phi_{\xi_a^\prime}(Q)\mid Q\in\underline{Q_a^0}\}$ if $\alpha
=\theta(a)$, and define both to be empty sets otherwise. Then
the argument in the preceding paragraph implies that
$\O_{\{\rho_{\beta\alpha}\}}(\{\underline{K_\alpha^0}\},
\{\underline{O_\alpha^0}\})$ is the image of the map
$\phi_{\bar{\gamma}}$. Now for any $\O_{\{\rho_{ba}\}}
(\{\underline{L_a}\},\{\underline{Q_a}\})$, we define
$\underline{K_\alpha}=\{\phi_{\xi_a}(L)\mid L\in\underline{L_a}\}$,
$\underline{O_\alpha}=\{\phi_{\xi_a^\prime}(Q)\mid Q\in\underline{Q_a}\}$
if $\alpha=\theta(a)$, and define both to be empty sets otherwise.
Then $\O_{\{\rho_{ba}\}}(\{\underline{L_a}\cup\underline{L_a^0}\},
\{\underline{Q_a}\cup\underline{Q_a^0}\})$ is mapped onto
$\O_{\{\rho_{\beta\alpha}\}}(\{\underline{K_\alpha}\cup
\underline{K_\alpha^0}\},\{\underline{O_\alpha}\cup
\underline{O_\alpha^0}\})$ under $\phi_{\bar{\gamma}}$.
Hence $\phi_{\bar{\gamma}}$ is an open map.

We summarize for the special case where $\tau$ is induced by $\sigma$
and $(3.2.5)$ is satisfied, and $\overline{\Gamma}_{\sigma\tau}=
\{\bar{\gamma}=(\theta,\{\xi_a\},\{\xi_a^\prime\})\mid (3.2.6),
(3.2.7\;a), (3.2.7\;b) \mbox{ are satisfied }\}\neq\emptyset$. Then each
element $\bar{\gamma}$ of $\overline{\Gamma}_{\sigma\tau}$ is assigned with
a local homeomorphism $\phi_{\bar{\gamma}}$, which clearly depends only
on the equivalence class of $\bar{\gamma}$ defined by $(3.1.4)$.
According to Lemma 3.1.1, the set $\Gamma_{\sigma\tau}$ is defined
to be the orbit space $(\Gamma_{\sigma\kappa}\times\Gamma_{\tau\kappa})
/G_{\kappa}$ where $\kappa$ is induced by $\tau$ with $(3.1.2)$
satisfied. There is a natural surjective mapping $\overline{\Gamma}
_{\sigma\tau}\rightarrow\Gamma_{\sigma\tau}$ given by $\bar{\gamma}
\mapsto [\bar{\gamma}\circ\gamma_0,\gamma_0]$. The assignment
$\bar{\gamma}\mapsto\phi_{\bar{\gamma}}$ factors through it, so
that for any $\bar{\gamma}\in\overline{\Gamma}_{\sigma\tau}$,
$\phi_{\bar{\gamma}}$ depends only on its image $\gamma\in
\Gamma_{\sigma\tau}$.

More generally for any pair $(\sigma,\tau)$ of equivalent
admissible homomorphisms, there are natural surjective mappings
$(\overline{\Gamma}_{\sigma\hat{\kappa}}\times
\overline{\Gamma}_{\tau\hat{\kappa}})/G_{\hat{\kappa}}\rightarrow
\Gamma_{\sigma\tau}$, where $\hat{\kappa}$ is any admissible homomorphism
induced by $\sigma$, $\tau$ with $(3.2.5), (3.2.6)$ satisfied. These
mappings are given by $[\bar{\gamma}_1,\bar{\gamma}_2]\mapsto
[\bar{\gamma}_1\circ\gamma_0,\bar{\gamma}_2\circ\gamma_0]\in
(\Gamma_{\sigma\kappa}\times\Gamma_{\tau\kappa})/G_{\kappa}=
\Gamma_{\sigma\tau}$, where $\kappa$ is any homomorphism induced by
$\hat{\kappa}$ with $(3.1.2)$ satisfied with respect to $\sigma$ and $\tau$,
and $\forall \gamma_0\in\Gamma_{\hat{\kappa}\kappa}$. Now given any
$\gamma\in\Gamma_{\sigma\tau}$ which is the image of $[\bar{\gamma}_1,
\bar{\gamma}_2]\in (\Gamma_{\sigma\hat{\kappa}}\times
\Gamma_{\tau\hat{\kappa}})/G_{\hat{\kappa}}$ under the above mappings,
we define
$$
\phi_{\gamma}=\mbox{ germ of }\phi_{{\bar{\gamma}}_1}\circ
\phi_{{\bar{\gamma}}_2}^{-1}. \leqno (3.2.14)
$$

In order to verify $(3.2.4)$, we need one more piece of identities.
Suppose $\bar{\gamma}_1, \bar{\gamma}_2\in\overline{\Gamma}_{\sigma\tau}$
where $\tau$ is induced by $\sigma$ and $(3.2.5), (3.2.6)$ is satisfied.
Then $\phi_{\bar{\gamma}_1}^{-1}\circ\phi_{\bar{\gamma}_2}$ is a local
homeomorphism from an open neighborhood of $\tau$ in $\O_{\{\rho_{ba}\}}$
to an open neighborhood of $\tau$ in $\O_{\{\rho_{ba}\}}$, leaving
$\tau$ fixed. By Lemma 3.1.2, this must be induced by the action
of an element in $G_\tau$. It can be computed as follows. Let
$\bar{\gamma}_1=\bar{\gamma}\circ g_1$, $\bar{\gamma}_2
=\bar{\gamma}\circ g_2$, where $g_1,g_2\in G_{\tau}$. Then
$\phi_{\bar{\gamma}_1}^{-1}\circ\phi_{\bar{\gamma}_2}=\phi_{g_1^{-1}g_2}$,
where $\phi_{g_1^{-1}g_2}$ is the action on $\O_{\{\rho_{ba}\}}$
induced by $g_1^{-1}g_2\in G_{\tau}$. On the other hand, for any homomorphism
$\kappa$ induced by $\sigma$, $\tau$ with $(3.1.2)$ satisfied,
$(\bar{\gamma}_1\circ\gamma_0)\times (\bar{\gamma}_2\circ\gamma_0)
=(\bar{\gamma}\circ g_1\circ\gamma_0)\times (\bar{\gamma}\circ
g_2\circ \gamma_0)=(\bar{\gamma}\circ\gamma_0\circ\epsilon(\gamma_0)^{-1}(g_1))
\times (\bar{\gamma}\circ\gamma_0\circ\epsilon(\gamma_0)^{-1}(g_2))
=(\epsilon(\gamma_0)^{-1}(g_1))^{-1}\epsilon(\gamma_0)^{-1}(g_2)
=\epsilon(\gamma_0)^{-1}(g_1^{-1}g_2)\in G_\kappa$ for any $\gamma_0\in
\Gamma_{\tau\kappa}$. Thus
$$
\phi_{\bar{\gamma}_1}^{-1}\circ\phi_{\bar{\gamma}_2}=
\phi_{\epsilon(\gamma_0)((\bar{\gamma}_1\circ\gamma_0)\times
(\bar{\gamma}_2\circ\gamma_0))},\; \forall \bar{\gamma}_1,\bar{\gamma}_2\in
\overline{\Gamma}_{\sigma\tau}, \gamma_0\in\Gamma_{\tau\kappa}.
 \leqno (3.2.15)
$$
To verify $(3.2.4)$, recall that if we write $\gamma=[\bar{\gamma}_1\circ
\gamma_0,\bar{\gamma}_2\circ\gamma_0]$ and $\gamma^\prime=[\bar{\gamma}_3
\circ\gamma_0,\bar{\gamma}_4\circ\gamma_0]$, we have $\gamma\circ\gamma^\prime
=[(\bar{\gamma}_1\circ\gamma_0)\circ ((\bar{\gamma}_2\circ\gamma_0)\times
(\bar{\gamma}_3\circ\gamma_0)),(\bar{\gamma}_4\circ\gamma_0)]
=[(\bar{\gamma}_1\circ\epsilon(\gamma_0)((\bar{\gamma}_2\circ\gamma_0)\times
(\bar{\gamma}_3\circ\gamma_0)))\circ\gamma_0,\bar{\gamma}_4\circ\gamma_0]$.
Hence we deduce with the aid of $(3.2.15)$ that
\begin{eqnarray*}
\phi_{\gamma\circ\gamma^\prime}
& = &\mbox{ germ of }
      \phi_{\bar{\gamma}_1\circ\epsilon(\gamma_0)((\bar{\gamma}_2\circ\gamma_0)
      \times (\bar{\gamma}_3\circ\gamma_0))}\circ \phi_{\bar{\gamma}_4}^{-1}\\
& = &\mbox{ germ of }
      \phi_{\bar{\gamma}_1}\circ\phi_{\epsilon(\gamma_0)((\bar{\gamma}_2\circ
       \gamma_0)\times (\bar{\gamma}_3\circ\gamma_0))}\circ
       \phi_{\bar{\gamma}_4}^{-1}\\
& = &\mbox{ germ of }
      \phi_{\bar{\gamma}_1}\circ (\phi_{\bar{\gamma}_2}^{-1}\circ
       \phi_{\bar{\gamma}_3})\circ\phi_{\bar{\gamma}_4}^{-1}\\
& = &\mbox{ germ of }
      (\phi_{\bar{\gamma}_1}\circ\phi_{\bar{\gamma}_2}^{-1})
       \circ (\phi_{\bar{\gamma}_3}\circ\phi_{\bar{\gamma}_4}^{-1})
       =\phi_{\gamma}\circ\phi_{\gamma^\prime}
\end{eqnarray*}
for any $\gamma\in\Gamma_{\sigma\tau}$ and
$\gamma^\prime\in\Gamma_{\tau\kappa}$.

It is easy to check that if $\sigma=\phi_{\bar{\gamma}}(\tau)$ is also in
$\O_{\{\rho_{ba}\}}$, then $\Gamma_{\sigma\tau}\subset G_{\{\rho_{ba}\}}$
and $\phi_\gamma$ is induced by the action of $G_{\{\rho_{ba}\}}$
on $\O_{\{\rho_{ba}\}}$.

Finally, for any $\tau^\prime\in\mbox{Domain }(\phi_{\bar{\gamma}})$,
$\sigma^\prime=\phi_{\bar{\gamma}}(\tau^\prime)$, the equivalence
class $\gamma$ of $\bar{\gamma}$ is assigned with a $\gamma^\prime\in
\Gamma_{\sigma^\prime\tau^\prime}$ such that the germ of
$\phi_{\bar{\gamma}}$ at $\tau^\prime$ equals $\phi_{\gamma^\prime}$.
Clearly, $\gamma^\prime$ is the equivalence class of $\bar{\gamma}$ in
$\Gamma_{\sigma^\prime\tau^\prime}$.

\hfill $\Box$

\subsection{Proof of Theorem 1.4}

Let $X, X^\prime$ be smooth orbifolds in the more general sense,
ie., the group actions on a local uniformizing system need not to be
effective. Assume without loss of generality that $X, X^\prime$ are
connected. Note that $X,X^\prime$ are paracompact, locally compact
and Hausdorff as defined in the beginning of \S 3.2. In particular,
the space of maps $[X;X^\prime]$ is canonically an orbispace.

We shall consider in this subsection, for any $r\geq 1$, the subset
$[X;X^\prime]^r\subset [X;X^\prime]$ of maps of $C^r$ class, i.e.,
the set of equivalence classes of (admissible) homomorphisms
$(\{f_\alpha\},\{\rho_{\beta\alpha}\})$ where each $f_\alpha:
\widehat{U_\alpha}\rightarrow\widehat{U_\alpha^\prime}$ is of $C^r$
class between domains in the Euclidean spaces.

First of all, we shall give a compatible topology to $[X;X^\prime]^r$,
stronger than the induced topology from $[X;X^\prime]$. For this
purpose, we introduce
$$
\O_{\{\rho_{\beta\alpha}\}}^r=\{\sigma\mid\sigma=(\{f_\alpha\},
\{\rho_{\beta\alpha}\})\mbox{ is admissible and } f_\alpha \mbox{
is } C^r, \forall\alpha\}. \leqno (3.3.1)
$$
We give each $\O_{\{\rho_{\beta\alpha}\}}^r$ a topology that is
generated by the subsets of the following type:
$$
\O_{\{\rho_{\beta\alpha}\}}^r(\{f_\alpha\},\{K_\alpha\},
\{\epsilon_\alpha\})=\{\{f_\alpha^\prime\}\mid ||f_\alpha^\prime
-f_\alpha||_{K_\alpha}^{C^r}<\epsilon_\alpha,\forall\alpha\},
\leqno (3.3.2)
$$
where $\{f_\alpha\}\in\O_{\{\rho_{\beta\alpha}\}}^r$, each
$K_\alpha$ is a compact subset of $\widehat{U_\alpha}$, each
$\epsilon_\alpha$ is a positive real number, and
$||\;||_{K_\alpha}^{C^r}$ is the $C^r$-norm of a $C^r$-map between
domains in the Euclidean spaces, defined over $K_\alpha$. Note that
the set of subsets in $(3.3.2)$ is actually a base of the topology on
$\O_{\{\rho_{\beta\alpha}\}}^r$. The group $G_{\{\rho_{\beta\alpha}\}}$
in $(3.2.2)$ acts continuously on $\O_{\{\rho_{\beta\alpha}\}}^r$
with the given topology, and the orbit space
$[\O_{\{\rho_{\beta\alpha}\}}^r]$ is a subset of $[X;X^\prime]^r$ by
Lemma 3.1.2. We denote the set of orbits of the subsets in $(3.3.2)$ by
$[\O_{\{\rho_{\beta\alpha}\}}^r(\{f_\alpha\},\{K_\alpha\},
\{\epsilon_\alpha\})]$, which all together form a base of the quotient
topology on $[\O_{\{\rho_{\beta\alpha}\}}^r]$. Finally, $[X;X^\prime]^r$
is given a topology which is generated by the set of subsets
$$
[\O_{\{\rho_{\beta\alpha}\}}^r(\{f_\alpha\},\{K_\alpha\},
\{\epsilon_\alpha\})]
$$
for all possible data $\{\rho_{\beta\alpha}\},\{f_\alpha\},
\{K_\alpha\},\{\epsilon_\alpha\}$. Call this topology the $C^r$-topology.

\begin{prop}
The space of $C^r$-maps $[X;X^\prime]^r$ with the $C^r$-topology
is an orbispace under a canonical \'{e}tale topological groupoid.
\end{prop}

\pf
The inclusions $\O_{\{\rho_{\beta\alpha}\}}^r\subset
\O_{\{\rho_{\beta\alpha}\}}$ are continuous, and each
$\phi_{\bar{\gamma}}$ in Lemma 3.2.2 sends a $C^r$ admissible
homomorphism to a $C^r$ admissible homomorphism. Hence each
$\phi_{\bar{\gamma}}$ is a bijection from an open set in
$\O_{\{\rho_{ba}\}}^r$ to an open set in $\O_{\{\rho_{\beta\alpha}\}}^r$.
On the other hand, $\phi_{\bar{\gamma}}$ is clearly continuous with
respect to the topology generated by the subsets in $(3.3.2)$. Thus
each $\phi_{\bar{\gamma}}$ is a local homeomorphism defined on the
space $U^r=\bigsqcup\O_{\{\rho_{\beta\alpha}\}}^r$.

Consider the groupoid $\Gamma^r=\bigcup_{(\sigma,\tau)}\Gamma_{\sigma\tau}$
where $\sigma,\tau$ are running over the set of $C^r$ admissible
homomorphisms. It is naturally identified with the set of germs of the
elements in the pseudogroup generated by $\{\phi_{\bar{\gamma}}\}$,
which acts on $U^r=\bigsqcup\O_{\{\rho_{\beta\alpha}\}}^r$ by local
homeomorphisms. Given with the \'{e}tale topology, $\Gamma^r$ becomes an
\'{e}tale topological groupoid, giving a canonical orbispace structure on
$[X;X^\prime]^r$, as one argues in the proof of Theorem 1.3.

\hfill $\Box$

We remark that when $X, X^\prime$ are smooth manifolds, $[X;X^\prime]^r$
is simply the space of $C^r$-maps from $X$ to $X^\prime$ given with
the Whitney topology, cf. \cite{Hir}.

\begin{lem}
Let $E\rightarrow X^\prime$ be a smooth orbifold vector bundle over
$X^\prime$.
\begin{itemize}
\item Given any $C^r$ homomorphism $\sigma$, there is a canonical
pull-back $C^r$ orbifold vector bundle $E_\sigma\rightarrow X$,
with a $C^r$ bundle morphism $\overline{\Phi}:E_\sigma\rightarrow E$
covering the equivalence class of $\sigma$, $\Phi:X\rightarrow
X^\prime$.
\item For any $\gamma\in\Gamma_{\sigma\tau}$, there is a bundle
isomorphism $\Xi(\gamma):E_\tau\rightarrow E_\sigma$, which is
compatible with the bundle morphisms $\overline{\Phi}:
E_\sigma\rightarrow E$ and satisfies $\Xi(\gamma_2\circ\gamma_1)=
\Xi(\gamma_2)\circ\Xi(\gamma_1)$. In particular, $G_\sigma$ acts
on $E_\sigma$, and $\Xi(\gamma)$ is $\epsilon(\gamma)$-equivariant.
\end{itemize}
\end{lem}

\pf
A smooth orbifold vector bundle of rank $n$ over $X^\prime$ may be given
by a smooth homomorphism of topological groupoids $\mu:
\Gamma\{U_{\alpha^\prime}^\prime\}\rightarrow {\em GL}_n(\R)$.
Suppose $\sigma=(\{f_\alpha\},\{\rho_{\beta\alpha}\}):\Gamma\{U_\alpha\}
\rightarrow\Gamma\{U_{\alpha^\prime}^\prime\}$. We define the $C^r$ orbifold
vector bundle $E_\sigma\rightarrow X$ to be the one given by the $C^r$
homomorphism $\mu\circ\sigma:\Gamma\{U_\alpha\}\rightarrow {\em GL}_n(\R)$,
where we represent $E\rightarrow X^\prime$ by $\mu:
\Gamma\{U_{\alpha^\prime}^\prime\}\rightarrow {\em GL}_n(\R)$.

More concretely, $\mu:\Gamma\{U_{\alpha^\prime}^\prime\}
\rightarrow {\em GL}_n(\R)$ is given by a collection of smooth
maps $\{\mu_{\xi^\prime}\mid \xi^\prime\in T(U_{\alpha^\prime}^\prime,
U_{\beta^\prime}^\prime)\}$, where $\mu_{\xi^\prime}:\mbox{Domain }
(\phi_{\xi^\prime})\rightarrow {\em GL}_n(\R)$, such that
$$
\mu_{\eta^\prime}(\phi_{\xi^\prime}(x))\mu_{\xi^\prime}(x)=
\mu_{\eta^\prime\circ\xi^\prime(x)}(x), \; \forall x\in
\phi_{\xi^\prime}^{-1}(\mbox{Domain }(\phi_{\eta^\prime})).
\leqno (3.3.3)
$$
The smooth orbifold vector bundle $E\rightarrow X^\prime$ is then given by
a collection of local trivializations $\{\widehat{U_{\alpha^\prime}^\prime}
\times\R^n\rightarrow\widehat{U_{\alpha^\prime}^\prime}\}$, with the
action of $G_{U_{\alpha^\prime}^\prime}$ given by the formula
$$
g\cdot (x,v)=(g\cdot x,\mu_g(x)(v)), \forall g\in
G_{U_{\alpha^\prime}^\prime}, (x,v)\in\widehat{U_{\alpha^\prime}^\prime}
\times\R^n, \leqno (3.3.4)
$$
and the set of transition functions $\{\varphi(\xi^\prime)\mid\xi^\prime
\in T(U_{\alpha^\prime}^\prime,U_{\beta^\prime}^\prime)\}$ given by
$$
\varphi(\xi^\prime):(x,v)\mapsto (\phi_{\xi^\prime}(x),
\mu_{\xi^\prime}(x)(v)), \;\forall x\in\mbox{Domain }(\phi_{\xi^\prime}),
v\in\R^n. \leqno (3.3.5)
$$
Likewise, $\mu^\sigma=\mu\circ\sigma$ is given by the collection of
$C^r$ maps $\{\mu_\xi^\sigma\mid \xi\in T(U_\alpha,U_\beta)\}$,
where $\mu^\sigma_\xi:\mbox{Domain }(\phi_\xi)\rightarrow {\em GL}_n(\R)$
is defined by $\mu^\sigma_\xi=\mu_{\rho_{\beta\alpha}(\xi)}\circ f_\alpha$,
$\forall \xi\in T(U_\alpha,U_\beta)$.
The $C^r$ orbifold vector bundle $E_\sigma\rightarrow X$ is defined by the
collection of local trivializations $\{\widehat{U_\alpha}\times\R^n
\rightarrow\widehat{U_\alpha}\}$, with the action of $G_{U_\alpha}$ given
by the formula
$$
g\cdot (x,v)=(g\cdot x,\mu_{\rho_\alpha(g)}(f_\alpha(x))(v)),
\forall g\in G_{U_\alpha}, (x,v)\in\widehat{U_\alpha}\times\R^n,
\leqno (3.3.6)
$$
and the set of transition functions $\{\varphi(\xi)\mid\xi\in
T(U_\alpha,U_\beta)\}$ given by
$$
\varphi(\xi):(x,v)\mapsto (\phi_{\xi}(x),
\mu_{\rho_{\beta\alpha}(\xi)}(f_\alpha(x))(v)), \;
\forall x\in\mbox{Domain }(\phi_{\xi}), v\in\R^n.
\leqno (3.3.7)
$$
The $C^r$ bundle morphism $\overline{\Phi}:E_\sigma\rightarrow E$ is given
by a collection of $C^r$ maps $\{\bar{f}_\alpha\}$, where
$$
\bar{f}_\alpha: (x,v)\mapsto (f_\alpha(x),v),\;
\forall (x,v)\in\widehat{U_\alpha}\times\R^n, \leqno (3.3.8)
$$
such that $(\phi_{\rho_{\beta\alpha}(\xi)},\mu_{\rho_{\beta\alpha}(\xi)})
\circ\bar{f}_\alpha=\bar{f}_\beta\circ (\phi_\xi,\mu^\sigma_\xi)$
for any $\alpha,\beta$ and $\xi\in T(U_\alpha,U_\beta)$.

As for the bundle isomorphisms $\Xi(\gamma):E_\tau\rightarrow E_\sigma$,
$\gamma\in\Gamma_{\sigma\tau}$, we first consider the special case where
$\tau$ is induced by $\sigma$ with $(3.1.2)$ satisfied. Then given
any $\gamma\in\Gamma_{\sigma\tau}$, we pick a representative $\bar{\gamma}
=(\theta,\{\xi_a\},\{\xi_a^\prime\})$ of $\gamma$, and define
$\Xi(\bar{\gamma}):E_\tau\rightarrow E_\sigma$ by a collection of
maps $\{\Xi(\bar{\gamma})_a\}$, where
$$
\Xi(\bar{\gamma})_a: (x,v)\mapsto (\phi_{\xi_a}(x),\mu_{\xi_a^\prime}
(f_a(x))(v)),\; \forall (x,v)\in\widehat{U_a}\times\R^n.
\leqno (3.3.9)
$$
The equivalence class of $\{\Xi(\bar{\gamma})_a\}$, which is the
bundle isomorphism $\Xi(\bar{\gamma}):E_\tau\rightarrow E_\sigma$,
depends only on the equivalence class of $\bar{\gamma}$ as defined by
$(3.1.4)$. We set $\Xi(\gamma)=\Xi(\bar{\gamma})$. More generally,
for any pair $(\sigma,\tau)$ of equivalent homomorphisms, we define
$\Xi(\gamma)=\Xi(\gamma_1)\circ\Xi(\gamma_2)^{-1}$ where $\gamma=
[\gamma_1,\gamma_2]$. We leave the verification for the claimed
properties of $\Xi(\gamma)$ in the lemma to the reader, which is
straightforward.

\hfill $\Box$

Now we come to the proof of Theorem 1.4.

\begin{thm}
Suppose $X$ is compact.
\begin{itemize}
\item [{(1)}] The space of $C^r$-maps $[X;X^\prime]^r$
is naturally a smooth Banach orbifold. More concretely,
$[X;X^\prime]^r$ is Hausdorff and second countable, and there exists an
open cover $\{\O_i\}$ of $[X;X^\prime]^r$ with the following properties.
\begin{itemize}
\item [{(i)}] Each $\O_i$ is associated with a triple $(\widehat{\O_i},
G_i,\pi_i)$, where $\widehat{\O_i}$ is an open ball in a Banach space,
$G_i$ is a finite group acting linearly on $\widehat{\O_i}$, and
$\pi_i:\widehat{\O_i}\rightarrow\O_i$ is a continuous map inducing
a homeomorphism $\widehat{\O_i}/G_i\cong \O_i$.
\item [{(ii)}] There is an \'{e}tale topological groupoid $\Gamma$
with the space of units $U=\bigsqcup_i\widehat{\O_i}$  such that
{\em (a)} $\Gamma\backslash U=[X;X^\prime]^r$, {\em (b)} the restriction
of $\Gamma$ to each $\widehat{\O_i}$ is the product groupoid
$G_i\times\widehat{\O_i}$, and {\em (c)} for each local section $s$ of
$\alpha:\Gamma\rightarrow U$, i.e., $\alpha\circ s=Id$, the map
$\omega\circ s$ is a local diffeomorphism between Banach spaces.
\end{itemize}
\item [{(2)}] The set of $C^l$-maps is a dense subset of $[X;X^\prime]^r$
for all $l\geq r$.
\end{itemize}
\end{thm}

\pf
(1) We fix an auxiliary Riemannian metric on $X^\prime$, whose existence
is ensured by the paracompactness of $X^\prime$. This is to say
that there is a family of Riemannian metrics $\{g_{i^\prime}\}$,
where each $g_{i^\prime}$ is a $G_{U_{i^\prime}^\prime}$-equivariant
Riemannian metric on the local chart $\widehat{U_{i^\prime}^\prime}$,
with respect to which the local diffeomorphisms $\{\phi_{\xi^\prime}\mid
\xi^\prime\in T(U_{i^\prime}^\prime,U_{j^\prime}^\prime)\}$ are
isometric.

Let $\sigma=\{f_\alpha\}\in\O_{\{\rho_{\beta\alpha}\}}^r$ be any $C^r$
admissible homomorphism. By Lemma 3.3.2, there is a canonical pull-back
$C^r$ orbifold vector bundle $E_\sigma\rightarrow X$, with a bundle
morphism $\overline{\Phi}:E_\sigma\rightarrow TX^\prime$ covering the
equivalence class of $\sigma$, $\Phi:X\rightarrow X^\prime$. Note that
$E_\sigma\rightarrow X$ inherits a natural metric from $TX^\prime$.

Denote by $C^r(E_\sigma)$ the space of $C^r$ sections of $E_\sigma$.
We claim that $C^r(E_\sigma)$ is canonically identified with
the set of admissible representatives of $C^r$ sections which are
defined over an arbitrarily given locally finite, admissible cover
$\{U_\alpha\}$, namely, the set of systems of local $C^r$ equivariant
sections $\{s_\alpha:\widehat{U_\alpha}\rightarrow\widehat{U_\alpha}
\times\R^n\}$, which can be extended over to the compact closure of
$\widehat{U_\alpha}$ and satisfy $s_\beta\circ\phi_{\xi}(x)=\varphi(\xi)
\circ s_\alpha(x),\; \forall x\in\mbox{Domain }(\phi_\xi),\xi\in
T(U_\alpha,U_\beta)$, cf. $(3.3.7)$ for the definition of transition
functions $\varphi(\xi)$. To see this, we first observe that the
mapping $\{s_\alpha\}\mapsto [\{s_\alpha\}]$ is injective by Lemma 3.1.2.
On the other hand, for any representative $\{s_i\}$ of a $C^r$ section
which is defined over a refinement $\{U_i\}$ of $\{U_\alpha\}$, we are
able to find a representative $\{s_\alpha\}$ of the $C^r$ section, which
is defined over $\{U_\alpha\}$, as follows. Consider the corresponding
representatives of the zero section, $\{0_i\}$ and $\{0_\alpha\}$. By
Lemma 3.2.2, there is a local homeomorphism $\phi_{\bar{\gamma}}$ from
an open neighborhood of $\{0_i\}$ to an open neighborhood of $\{0_\alpha\}$,
which is clearly linear in this case. For any $\{s_i\}$, we pick a small
$\epsilon>0$ such that $\{\epsilon s_i\}$ lies in the said open
neighborhood  of $\{0_i\}$ where $\phi_{\bar{\gamma}}$ is defined,
and then we define $\{s_\alpha\}=\epsilon^{-1}\phi_{\bar{\gamma}}
(\{\epsilon s_i\})$. Here the existence of such an $\epsilon>0$ relies
on the assumption that $X$ is compact. It is easy to see that $\{s_i\}$
is indeed equivalent to $\{s_\alpha\}$. Hence $C^r(E_\sigma)$ is
canonically identified with the space of admissible compatible systems
of local $C^r$ equivariant sections over an arbitrarily given locally finite,
admissible cover.

We define a norm on $C^r(E_\sigma)$ as follows, which makes
$C^r(E_\sigma)$ into a Banach space. For any $\{s_\alpha\}\in
C^r(E_\sigma)$, we define its norm by $||\{s_\alpha\}||=
\max_\alpha ||s_\alpha||_{C^r}$, where $||s_\alpha||_{C^r}$ is the
$C^r$-norm of the section $s_\alpha$\footnote{note that $s_\alpha$
can be extended over to the compact closure of $\widehat{U_\alpha}$
by assumption.} defined using the metric on $E_\sigma$. Here the index
set $\{\alpha\}$ is finite because $\{U_\alpha\}$ is locally finite
and $X$ is compact.

Let $\tau$ be any admissible $C^r$ homomorphism which is equivalent to
$\sigma$. According to Lemma 3.3.2, for any $\gamma\in\Gamma_{\sigma\tau}$,
there is a canonical orbifold bundle isomorphism $\Xi(\gamma):E_\tau
\rightarrow E_\sigma$. For simplicity, we still denote by $\Xi(\gamma)$
the induced isomorphism between the Banach spaces $C^r(E_\tau)$
and $C^r(E_\sigma)$, which is clearly norm-preserving. In particular,
the isotropy group $G_\sigma$ defined in $(3.1.1)$, which is
a finite group in this case, acts linearly on $C^r(E_\sigma)$,
preserving the norm of $C^r(E_\sigma)$. Moreover, the map $\Xi(\gamma):
C^r(E_\tau)\rightarrow C^r(E_\sigma)$ is equivariant with
respect to the isomorphism $\epsilon(\gamma):G_\tau\rightarrow
G_\sigma$ constructed in Lemma 3.1.1.

Let $\widehat{\O_\sigma}(\epsilon)$ be the open ball of radius
$\epsilon$ in $C^r(E_\sigma)$. Denote by $\mbox{Exp}$ the exponential
map on each $\widehat{U_\alpha^\prime}$. Then the map $\Theta_\sigma:
\widehat{\O_\sigma}(\epsilon)\rightarrow\O_{\{\rho_{\beta\alpha}\}}^r$
defined by
$$
\Theta_\sigma(s)=\{\mbox{Exp}_{f_\alpha(x)}s_\alpha\} \mbox{ where }
s=\{s_\alpha\} \leqno (3.3.10)
$$
is an open embedding when $\epsilon>0$ is sufficiently small.
Moreover, for any $\gamma\in\Gamma_{\sigma\tau}$ and any representative
$\bar{\gamma}$ of $\gamma$, we have $\Theta_\sigma\circ
\Xi(\gamma)=\phi_{\bar{\gamma}}\circ\Theta_\tau$ on the intersection
of their domains, where $\phi_{\bar{\gamma}}$ is the local homeomorphism
constructed in Lemma 3.2.2. In particular, the open embedding
$\Theta_\sigma$ is equivariant with respect to the natural inclusion
$G_\sigma\hookrightarrow G_{\{\rho_{\beta\alpha}\}}$.

Let $\O_\sigma(\epsilon)$ be the image of $\Theta_\sigma
(\widehat{\O_\sigma}(\epsilon))$ in $[\O_{\{\rho_{\beta\alpha}\}}^r]$,
and let $\pi_\sigma:\widehat{\O_\sigma}(\epsilon)\rightarrow
\O_\sigma(\epsilon)$ be the corresponding orbit map. We shall prove
next that for sufficiently small $\epsilon>0$, the map $\pi_\sigma$
induces a homeomorphism between $\widehat{\O_\sigma}(\epsilon)/G_\sigma$
and $\O_\sigma(\epsilon)$.

Some digression first. (1) Observe that each $\widehat{U_\alpha^
\prime}$ is partitioned into a finite disjoint union $\bigsqcup
\widehat{U_\alpha^\prime}(H)$, where $\widehat{U_\alpha^\prime}(H)=
\{x\in\widehat{U_\alpha^\prime}|G_x=H\}$. Here $H$ is a subgroup
of $G_{U_\alpha^\prime}$, and $G_x$ is the stabilizer at $x$. Note that
$\widehat{U_\alpha^\prime}(H)$ has the following property: if
$g\in G_{U_\alpha^\prime}$ fixes a point in $\widehat{U_\alpha^\prime}(H)$,
then $g$ must be an element of $H$, which implies that $g$ actually
fixes the entire $\widehat{U_\alpha^\prime}(H)$. (2) For any $x\in\widehat
{U_\alpha^\prime}$, there exists a $R_x>0$ with the following
significance. For any $y,z\in\widehat{U_\alpha^\prime}$, if $y,z$ are
joined to $x$ by a geodesic with length less than or equal to $R_x$, and
$y=g\cdot z$ for some $g\in G_{U_\alpha^\prime}$, then $g$ lies in $G_x$.
End of digression.

Now for any $(\alpha,H)$, if $\widehat{U_\alpha^\prime}(H)$ intersects
with $f_\alpha(\widehat{U_\alpha})$, we pick a point $x_{\alpha,H}\in
\widehat{U_\alpha}$ such that $x_{\alpha,H}^\prime=f_\alpha(x_{\alpha,H})$
lies in $\widehat{U_\alpha^\prime}(H)$.
Since there are only finitely many $(\alpha,H)$, the set of positive
numbers $R_{x_{\alpha,H}^\prime}$ has a minimum $R_\sigma>0$.

We require $\epsilon\leq R_\sigma$. Now suppose there are $s_1,s_2\in
\widehat{\O_\sigma}(\epsilon)$ such that $[\Theta_\sigma(s_1)]
=[\Theta_\sigma(s_2)]$ in $[\O_{\{\rho_{\beta\alpha}\}}^r]$. Then there
is a $g=\{g_\alpha\}\in G_{\{\rho_{\beta\alpha}\}}$ such
that $\Theta_\sigma(s_2)=g\cdot\Theta_\sigma(s_1)$. By the assumption
that $\epsilon\leq R_\sigma$, we conclude that for any $\alpha$,
$g_\alpha$ fixes each $x^\prime_{\alpha,H}$, hence fixes the
entire $\widehat{U_\alpha^\prime}(H)$, where $H$ is a subgroup of
$G_{U_\alpha^\prime}$ such that $\widehat{U_\alpha^\prime}(H)$ intersects
with $f_\alpha(\widehat{U_\alpha})$. But $f_\alpha(\widehat{U_\alpha})$
is contained in the union of these $\widehat{U_\alpha^\prime}(H)$'s,
hence $f_\alpha(\widehat{U_\alpha})$ is fixed by $g_\alpha$.
This exactly means that $g=\{g_\alpha\}$ is an element of $G_\sigma$,
which implies that $\pi_\sigma$ induces a homeomorphism between
$\widehat{\O_\sigma}(\epsilon)/G_\sigma$ and $\O_\sigma(\epsilon)$.

We pick a collection of $\sigma_i$ and $\epsilon_i\in (0,R_{\sigma_i}]$
such that $\{\O_{\sigma_i}(\epsilon_i)\}$ is a cover of $[X;X^\prime]^r$.
We let $\O_i=\O_{\sigma_i}(\epsilon_i)$, and $\widehat{\O_i}=\widehat
{\O_{\sigma_i}}(\epsilon_i)$, $G_i=G_{\sigma_i}$ and
$\pi_i=\pi_{\sigma_i}$. Then (1)-(i) of the theorem is proved.

We regard each $\widehat{\O_i}$ as an open subset of
$\O_{\{\rho_{\beta\alpha}\}}^r$ for some data $\{\rho_{\beta\alpha}\}$
via the open embeddings defined in $(3.3.10)$, and define $\Gamma$
to be the \'{e}tale topological groupoid which is the restriction of
the \'{e}tale topological groupoid $\Gamma^r$ in Proposition 3.3.1
to $U=\bigsqcup\widehat{\O_i}$. Then ${\em (a),(b)}$ in (1)-(ii) of the
theorem are immediate from the construction. As for ${\em (c)}$ of
(1)-(ii), a little bit elementary Riemannian geometry shows that it
boils down to the following local problem which can be easily checked: Let
$\Omega$ be a compact domain in a Euclidean space. Any self-diffeomorphism
of $\R^N$ induces a self-diffeomorphism of the Banach space
$C^r(\Omega;\R^N)$ of $C^r$-maps from $\Omega$ into $\R^N$.

\vspace{1.5mm}

It remains to check that $[X;X^\prime]^r$ is Hausdorff and second
countable.

\vspace{1.5mm}

To see that $[X;X^\prime]^r$ is Hausdorff, suppose $\{\Phi_n\}$
is a sequence of $C^r$ maps such that $\lim_{n\rightarrow\infty}
\Phi_n=\Phi$ and $\lim_{n\rightarrow\infty}\Phi_n=\Phi^\prime$. We need
to show that $\Phi=\Phi^\prime$. First of all, at the level of induced
maps between underlying spaces, we have $f=\lim_{n\rightarrow\infty}
f_n$, $f^\prime=\lim_{n\rightarrow\infty}f_n$, where $f_n,f$ and
$f^\prime$ are the induced maps of $\Phi_n$, $\Phi$, and
$\Phi^\prime$ respectively.
Since $X^\prime$ is Hausdorff, we have $f=f^\prime$. This allows us
to find a cover of local charts $\{U_\alpha\}$ on $X$, locally finite
and admissible, and a subset of local charts
$\{U_{\alpha^\prime}^\prime\}$ on $X^\prime$,
with a correspondence $U_\alpha\mapsto U_\alpha^\prime\in
\{U_{\alpha^\prime}^\prime\}$, such that there is an admissible
$C^r$ homomorphism $\sigma=(\{f_\alpha^{(1)}\},\{\rho_{\beta\alpha}
^{(1)}\}):\Gamma\{U_\alpha\}\rightarrow\Gamma\{U_{\alpha^\prime}^\prime\}$
representing $\Phi$, and there is an admissible $C^r$ homomorphism
$\tau=(\{f_\alpha^{(2)}\},\{\rho_{\beta\alpha}^{(2)}\}):
\Gamma\{U_\alpha\}\rightarrow\Gamma\{U_{\alpha^\prime}^\prime\}$
representing $\Phi^\prime$. On the other hand, since both
$G_{\{\rho_{\beta\alpha}^{(1)}\}}$ and $G_{\{\rho_{\beta\alpha}^{(2)}\}}$
are finite groups, there are sequences of representatives
$\{\sigma_n\}$, $\{\tau_n\}$ of $\Phi_n$,
such that (1) for any $n$, $\sigma_n\in\O_{\{\rho_{\beta\alpha}^{(1)}\}}^r$,
$\tau_n\in\O_{\{\rho_{\beta\alpha}^{(2)}\}}^r$, and (2)
$\lim_{n\rightarrow\infty}\sigma_n=\sigma$, $\lim_{n\rightarrow\infty}
\tau_n=\tau$. By Lemma 3.1.2, the fact that both $\sigma_n,\tau_n$
represent the map $\Phi_n$ implies that there are $g_{\alpha,n}\in
G_{U_\alpha^\prime}$
such that $\sigma_n$ is conjugate to $\tau_n$ by $\{g_{\alpha,n}\}$.
But each $G_{U_\alpha^\prime}$ is finite, and
$\max_\alpha |G_{U_\alpha^\prime}|<\infty$. Hence
there is a subsequence $n_i\rightarrow\infty$, with $g_{\alpha,n_i}=
g_\alpha$ independent of $n_i$. Now it is clear that $\sigma$
is conjugate to $\tau$ by $\{g_\alpha\}$, and hence $\Phi=\Phi^\prime$
as desired.

As for the second countability of $[X;X^\prime]^r$, we fix a countable
cover of local charts $\U_0^\prime\subset\U^\prime$ of $X^\prime$,
and a countable base of local charts $\U_0\subset\U$ of $X$. The
existence of $\U_0^\prime$ and $\U_0$ follows from the second countability
of $X^\prime$ and $X$. Now observe that: (a) because $\U_0,\U_0^\prime$
are countable, so is the set $\Lambda$ of data
$(\{U_\alpha\},\{U_\alpha^\prime\},\{\rho_{\beta\alpha}\})$,
where $\{U_\alpha\}$ is a finite subset of $\U_0$, with an assignment
$U_\alpha\mapsto U_\alpha^\prime\in\U_0^\prime$, and $\rho_{\beta\alpha}:
T(U_\alpha,U_\beta)\rightarrow T(U_\alpha^\prime,U_\beta^\prime)$ is a
mapping between finite sets, (b) the set $\{[\O^r_{\{\rho_{\beta\alpha}\}}]
\mid (\{U_\alpha\},\{U_\alpha^\prime\},\{\rho_{\beta\alpha}\})\in\Lambda\}$
is a cover of $[X;X^\prime]^r$, and (c) each space
$\O^r_{\{\rho_{\beta\alpha}\}}$ is second countable --- the validity of
this last assertion boils down to the second countability of the Banach
space $C^r(\Omega;\R^N)$ of $C^r$ maps from a compact domain $\Omega$ in a
Euclidean
space into $\R^N$. The second countability of $[X;X^\prime]^r$ follows
easily
from the observations (a), (b) and (c) above.

\vspace{1.5mm}

(2) In order to prove that the set of $C^l$-maps is dense in $[X;X^\prime]^r$
for any $l\geq r$, it suffices to show that for any $\O_{\{\rho_{ji}\}}^r$,
where we may assume $i,j\in\{1,2,\cdots,N\}$ for some $N>0$ since $X$ is
compact, and for any open neighborhood $\O_{\{\rho_{ji}\}}^r(\{f_i\},\{K_i\},
\{\epsilon_i\})$ of $\sigma=\{f_i\}$, there is a $\tau=\{f_i^\prime\}\in
\O_{\{\rho_{ji}\}}^r(\{f_i\},\{K_i\},\{\epsilon_i\})$ which is of $C^l$ class.
Since the set of indexes $i$ is finite, we may assume $\epsilon_i=\epsilon$,
which is independent of $i$, by taking a smaller neighborhood.
$\tau=\{f_i^\prime\}$
will be defined by induction on $i$. First, we pick any $C^l$-map $h_1$
satisfying $||h_1-f_1||_{C^r}<\epsilon$, where the $C^r$-norm is
taken over the closure of $\widehat{U_1}$. Then we make the $C^l$ map
$\rho_1$-equivariant by setting $f_1^\prime(x)=\frac{1}{|G_{U_1}|}
\sum_{g\in G_{U_1}}\rho_1(g^{-1})h_1(g\cdot x)$. Clearly $f_1^\prime$
is $\rho_1$-equivariant and
satisfies $||f_1^\prime-f_1||_{C^r}<\epsilon$. Now assume $f_i^\prime$
is defined and satisfies $||f_i^\prime-f_i||_{C^r}<\epsilon$ for $1\leq i
\leq n-1$, and $(\{f_i^\prime\},\{\rho_{ji}\})$, where $1\leq i\leq n-1$,
is a homomorphism. We shall construct $f_n^\prime$ as follows. Let
$A\subset\mbox{ closure of }\widehat{U_n}$ be the set of $x$ where there
exists a $\xi\in T(U_n,U_i)$ for some $i<n$, such that
$x\in\overline{\mbox{Domain }
(\phi_\xi)}$. For any $x\in A$, we define
$f_n^\prime(x)=\phi_{\rho_{in}(\xi)}^{-1}
\circ f_i^\prime\circ\phi_\xi(x)$, which is independent of the choice
of $\xi$ because $(\{f_i^\prime\},\{\rho_{ji}\})$, where $1\leq i\leq
n-1$, is a homomorphism. Moreover, $A$ is invariant under the action of
$G_{U_n}$,
and $f_n^\prime$ over $A$ is $\rho_n$-equivariant, of $C^l$ class, and
satisfies
$||f_n^\prime-f_n||_{A}^{C^r}<\epsilon$. Now we simply extend $f_n^\prime$
to a $C^l$-map over the closure of $\widehat{U_n}$, making it
$\rho_n$-equivariant and satisfy $||f_n^\prime-f_n||_{C^r}<\epsilon$.
To see that $(\{f_i^\prime\},\{\rho_{ji}\})$, where $1\leq i\leq n$, is
a homomorphism, it suffices to check that $\rho_{ki}(\eta\circ\xi({\bf a}))=
\rho_{kj}(\eta)\circ\rho_{ji}(\xi)(\underline{\{f_i^\prime\}}({\bf a}))$ for
any $1\leq i,j,k\leq n$ and any ${\bf a}\in\Lambda(\xi,\eta)$. When
$\epsilon>0$ is
sufficiently small, $\underline{\{f_i^\prime\}}({\bf a})
=\underline{\{f_i\}}({\bf a})$, from which the above equation
follows. Now by induction, $\tau=\{f_i^\prime\}\in\O_{\{\rho_{ji}\}}$, where
$1\leq i\leq N$, is defined, which is of $C^l$ class and lies in the given
neighborhood $\O_{\{\rho_{ji}\}}^r(\{f_i\},\{K_i\},\{\epsilon_i\})$. Hence
the set of $C^l$-maps is dense in $[X;X^\prime]^r$.

\vspace{1.5mm}

The proof of Theorem 3.3.3 is completed.

\hfill $\Box$

The rest of this subsection is concerned with two results which
are useful in the theory of pseudoholomorphic curves in symplectic
orbifolds, cf. eg. \cite{C2,C3}.

For the first one, let $X,Y,Z$ be compact smooth orbifolds. Consider
the mapping of composition $(\Phi,\Psi)\mapsto \Psi\circ\Phi$,
where $\Phi\in [X;Y]^r$, $\Psi\in [Y;Z]^r$. Similar argument as in
the proof of Proposition 3.2.3 shows that it defines a smooth
`map' from $[X;Y]^r\times [Y;Z]^r$ to $[X;Z]^r$. But note that the
term `map' is not quite clearly defined here because it is not
known whether $[X;Y]^r\times [Y;Z]^r$, $[X;Z]^r$ as smooth Banach
orbifolds satisfy (C2) or not, except for the case where
they are actually smooth Banach manifolds. However, as far as
applications are concerned, it suffices to only consider finite
dimensional suborbifolds of $[X;Y]^r\times [Y;Z]^r$ and $[X;Z]^r$.

\begin{prop}
Let $M,N$ be finite dimensional suborbifolds
of $[X;Y]^r\times [Y;Z]^r$ and $[X;Z]^r$ such that $\Psi\circ\Phi\in N$
for any $(\Phi,\Psi)\in M$. Then there is a canonically defined smooth
map of orbifolds from $M$ to $N$, which induces
$(\Phi,\Psi)\mapsto \Psi\circ\Phi$ between the underlying spaces.
\end{prop}

\pf
For each $(\Phi,\Psi)\in M$, pick admissible homomorphisms $\sigma,\tau$
representing $\Phi,\Psi$ respectively, such that $\kappa\equiv
\tau\circ\sigma$ is defined. Then as in the proof of Proposition
3.2.3, there are positive real numbers $\epsilon_\sigma,\epsilon_\tau$
and $\epsilon_\kappa$, and open balls of Banach spaces
$\widehat{\O_\sigma}(\epsilon_\sigma),\widehat{\O_\tau}(\epsilon_\tau)$
and $\widehat{\O_\kappa}(\epsilon_\kappa)$ as defined in the proof
of Theorem 3.3.3, such that for any $\sigma^\prime\in
\widehat{\O_\sigma}(\epsilon_\sigma)$, $\tau^\prime\in
\widehat{\O_\tau}(\epsilon_\tau)$, $\kappa^\prime\equiv\tau^\prime\circ
\sigma^\prime$ is defined and lies in
$\widehat{\O_\kappa}(\epsilon_\kappa)$. Let $(\widehat{U_{(\sigma,\tau)}},
G_{(\sigma,\tau)})$, $(\widehat{U_\kappa},G_\kappa)$ be local
uniformizing systems of $M,N$ at $(\Phi,\Psi)$, $\Psi\circ\Phi$, which are
induced from $(\widehat{\O_\sigma}(\epsilon_\sigma)\times
\widehat{\O_\tau}(\epsilon_\tau), G_\sigma\times G_\tau)$,
$(\widehat{\O_\kappa}(\epsilon_\kappa),G_\kappa)$ respectively.
By Proposition 2.1.3, we may choose $\epsilon_\sigma,\epsilon_\tau$
and $\epsilon_\kappa$ sufficiently small so that (C1),(C2) are
satisfied for these local uniformizing systems of $M,N$.

Next we define a homomorphism of topological groupoids from
$\Gamma\{U_{(\sigma,\tau)}\}$ to $\Gamma\{U_\kappa\}$. More
concretely, let $(\sigma_1^\prime,\tau_1^\prime)\in
\widehat{U_{(\sigma_1,\tau_1)}}, (\sigma_2^\prime,\tau_2^\prime)
\in\widehat{U_{(\sigma_2,\tau_2)}}$ such that $\sigma_1^\prime,
\sigma_2^\prime$ and $\tau_1^\prime,\tau_2^\prime$ are equivalent,
and let $\gamma_1\in\Gamma_{\sigma_2^\prime\sigma_1^\prime}$,
$\gamma_2\in\Gamma_{\tau_2^\prime\tau_1^\prime}$. Set $\kappa_1^\prime
\equiv \tau_1^\prime\circ\sigma_1^\prime\in\widehat{U_{\kappa_1}}$,
$\kappa_2^\prime\equiv\tau_2^\prime\circ\sigma_2^\prime\in
\widehat{U_{\kappa_2}}$. We shall define a $\gamma\in
\Gamma_{\kappa_2^\prime\kappa_1^\prime}$ such that
$((\gamma_1,\gamma_2),(\sigma_1^\prime,\tau_1^\prime))\mapsto
(\gamma,\kappa_1^\prime)$ is a groupoid homomorphism from
$\Gamma\{U_{(\sigma,\tau)}\}$ to $\Gamma\{U_\kappa\}$. Note that
it suffices to only consider the case where $\sigma_1^\prime,
\tau_1^\prime$ are induced by $\sigma_2^\prime,\tau_2^\prime$.
To fix the notation, let $\sigma_1^\prime=(\{f_a\},\{\rho_{ba}\}):
\Gamma\{U_a\}\rightarrow\Gamma\{V_a\}$, $\sigma_2^\prime=
(\{f_\alpha\},\{\rho_{\beta\alpha}\}):\Gamma\{U_\alpha\}
\rightarrow\Gamma\{V_\alpha\}$, $\tau_1^\prime=(\{h_s\},\{\zeta_{ts}\}):
\Gamma\{V_s\}\rightarrow\Gamma\{W_s\}$, and $\tau_2^\prime=(\{h_\mu\},
\{\zeta_{\nu\mu}\}):\Gamma\{V_\mu\}\rightarrow\Gamma\{W_\mu\}$.
Let $\imath:a\mapsto s,\jmath:\alpha\mapsto\mu$ be the mappings of
indexes such that $\kappa_1^\prime=(\{h_{\imath(a)}\circ f_a\},
\{\zeta_{\imath(b)\imath(a)}\circ\rho_{ba}\})$ and
$\kappa_2^\prime=(\{h_{\jmath(\alpha)}\circ f_\alpha\},
\{\zeta_{\jmath(\beta)\jmath(\alpha)}\circ\rho_{\beta\alpha}\})$.
Finally, let $\gamma_1,\gamma_2$ be represented by $(\theta_1,\{\xi_a\},
\{\xi_a^\prime\})$, $(\theta_2,\{\eta_s\},\{\eta_s^\prime\})$,
where $\theta_1:a\mapsto\alpha$, $\theta_2:s\mapsto\mu$, $\xi_a\in
T(U_a,U_{\theta_1(a)})$, $\xi_a^\prime\in T(V_a,V_{\theta_1(a)})$,
$\eta_s\in T(V_s,V_{\theta_2(s)})$ and $\eta_s^\prime\in
T(W_s,W_{\theta_2(s)})$. Note that $V_a=V_{\imath(a)}$, $V_{\theta_1(a)}=
V_{\jmath(\theta_1(a))}$, thus for each index $a$, set
$\hat{\xi_a}\equiv\eta_{\imath(a)}\circ (\xi_a^\prime)^{-1}\in
T(V_{\jmath(\theta_1(a))},V_{\theta_2(\imath(a))})$. Then
it is easy to check that $\kappa_1^\prime$ is induced by $\kappa_2^\prime$
via $(\theta_1,\{\xi_a\},\{(\zeta_{\theta_2(\imath(a))\jmath(\theta_1(a))}
(\hat{\xi_a}))^{-1}\circ\eta_{\imath(a)}^\prime\})$. We define $\gamma$
to be the corresponding equivalence class in $\Gamma_{\kappa_2^\prime
\kappa_1^\prime}$. It is a routine exercise to check that
$((\gamma_1,\gamma_2),(\sigma_1^\prime,\tau_1^\prime))\mapsto
(\gamma,\kappa_1^\prime)$ is well-defined, and that it is a groupoid
homomorphism from $\Gamma\{U_{(\sigma,\tau)}\}$ to $\Gamma\{U_\kappa\}$.
We leave the details to the reader.

The equivalence class of the above homomorphism from
$\Gamma\{U_{(\sigma,\tau)}\}$ to $\Gamma\{U_\kappa\}$ defines
a map from $M$ to $N$ (cf. Lemma 2.2.1), which is smooth because
$(\sigma_1^\prime,\tau_1^\prime)\mapsto\kappa_1^\prime$ is smooth, and it
induces $(\Phi,\Psi)\mapsto \Psi\circ\Phi$ between the underlying spaces
of $M,N$. Finally, observe that the map is canonically defined such that it
is actually independent of $\{U_{(\sigma,\tau)}\}$ and
$\{U_\kappa\}$.

\hfill $\Box$

For the second result, let $X,X^\prime$ be compact smooth
orbifolds. Suppose $M$ is a finite dimensional suborbifold of
$[X;X^\prime]^r$, which consists of only smooth maps from $X$
to $X^\prime$.

\begin{prop}
There is a canonically defined smooth map from $M\times X$ to
$X^\prime$ which induces the evaluation map $(\Phi,p)\mapsto
|\Phi|(p)$ between the underlying spaces, where $|\Phi|$ is the map
between the underlying spaces of $X,X^\prime$ induced by $\Phi$.
\end{prop}

\pf
For each $\Phi\in M$, pick an admissible homomorphism $\sigma=(\{f_\alpha\},
\{\rho_{\beta\alpha}\}):\Gamma\{U_\alpha(\sigma)\}\rightarrow\Gamma
\{U_{\alpha^\prime}^\prime(\sigma)\}$. Let $(\widehat{U_\sigma},G_\sigma)$
be a local uniformizing system of $M$ at $\Phi$, which is induced
from $(\widehat{\O_\sigma}(\epsilon_\sigma),G_\sigma)$ for some
sufficiently small $\epsilon_\sigma>0$. Furthermore, we may assume
that (C1), (C2) are satisfied for the local uniformizing systems
$\{(\widehat{U_\sigma},G_\sigma)\}$, cf. Proposition 2.1.3.

Next we define a groupoid homomorphism from $\Gamma\{U_\sigma\times
U_\alpha(\sigma)\}$ to $\Gamma\{U_{\alpha^\prime}^\prime(\sigma)\}$.
To this end, for any $(\tau^\prime,x)\in\widehat{U_\tau}\times
\widehat{U_a(\tau)}$, where $\tau^\prime=(\{f_a^\prime\},\{\rho_{ba}\})$,
and for any $\gamma\in\Gamma_{\sigma^\prime\tau^\prime}$ and $\xi\in
T(U_a(\tau),U_\alpha(\sigma))$ such that $x\in\mbox{Domain }(\phi_\xi)$,
where $\sigma^\prime=(\{f_{\alpha}^\prime\}, \{\rho_{\beta\alpha}\})\in
\widehat{U_\sigma}$ and $\gamma$ is represented by $(\theta,\{\xi_a\},
\{\xi_a^\prime\})$, we let
$\xi^\prime=\rho_{\alpha\theta(a)}(\xi\circ\xi_a^{-1})
\circ\xi_a^\prime(f_a^\prime(x))\in T(U_a^\prime(\tau),
U_\alpha^\prime(\sigma))$.
Then it is a routine exercise to check that $((\gamma,\xi),(\tau^\prime,x))
\mapsto (\xi^\prime, f_a^\prime(x))$ defines a groupoid homomorphism from
$\Gamma\{U_\sigma\times U_\alpha(\sigma)\}$ to
$\Gamma\{U_{\alpha^\prime}^\prime(\sigma)\}$.

The corresponding map from $M\times X$ to $X^\prime$ is smooth
because $(\tau^\prime,x)\mapsto f_a^\prime(x)$ is smooth. It is
clear that it induces the evaluation map $(\Phi,p)\mapsto
|\Phi|(p)$ between the underlying spaces, and that it is canonically
defined such that it is independent of the choices made on
$\{\widehat{U_\sigma}\}$. Details are left to the reader.

\hfill $\Box$

We end here with a discussion concerning group actions on orbifolds.
Let $\G$ be a Lie group and $X$ be a compact smooth orbifold in the
classical sense (ie. all local group actions on the uniformizing systems
are effective). By a smooth, effective left action of $\G$ on $X$, we
mean a smooth injective map $\iota:\G\rightarrow [X;X]^r$ for some $r$
such that the image of $\iota$ consists of only smooth maps, and
$\iota(1)=Id$, $\iota(g_1g_2)=\iota(g_1)\circ\iota(g_2)$ for
any $g_1,g_2\in\G$. Note that the image of $\iota$ is contained in an open
Banach submanifold of $[X;X]^r$, so that $\iota$ being a smooth map is
simply in the usual sense. Then as in Proposition 3.3.5, there
is a canonical smooth map of orbifolds from $\G\times X$ to $X$, which
induces the map $(g,p)\mapsto|\iota(g)|(p)$ between the underlying spaces.
(Note that the induced map defines a continuous left action of $\G$ on the
underlying space of $X$.) On the other hand, given any smooth map of orbifolds
$\psi:\G\times X\rightarrow X$, there is an associated map $\iota:\G\rightarrow
[X;X]^r$, with $\iota(g):X\rightarrow X$ being defined by the restriction
of $\psi$ on the suborbifold $\{g\}\times X$. Clearly $\iota(g)\in [X;X]^r$
is smooth. Suppose furthermore, $\iota(1)=Id$ and $\iota(g_1g_2)
=\iota(g_1)\circ\iota(g_2)$ for any $g_1,g_2\in\G$, and $\iota$ is
injective. Then $\iota:\G\rightarrow [X;X]^r$ defines a smooth, effective
left action of $\G$ on $X$. The only thing remains to be verified is
that $\iota$ is a smooth map. To see this, recall that $X$ is compact,
so that for any $g\in\G$, there is an open neighborhood $U$ of $g$ in
$\G$, such that the map $\psi:\G\times X\rightarrow X$ is represented
by a homomorphism $\sigma=(\{f_a\},\{\rho_{ba}\}):\Gamma\{U_a\}\rightarrow
\Gamma\{U_{a^\prime}^\prime\}$, where a neighborhood of $\{g\}\times X$ in
$\G\times X$ is covered by a subset $\{U\times V_i\}$ of $\{U_a\}$. Let
$\theta:i\mapsto a$ be the mapping of indexes such that
$U\times V_i=U_{\theta(i)}$. Then for any $h\in U$, $\iota(h)\in [X;X]^r$
is represented by the homomorphism $\tau_h=(\{f_{\theta(i)}(h,\cdot)\},
\{\rho_{\theta(j)\theta(i)}\})$. Hence $\iota$ is a smooth map.

Finally, suppose $\iota:\G\rightarrow [X;X]^r$ is a smooth, effective
left action of $\G$ on $X$, and $M\subset [X;X^\prime]^r$ is a
finite dimensional suborbifold such that for any $\Phi\in M$, $g\in\G$,
$\Phi\circ\iota(g)\in M$. Then as in Proposition 3.3.4, there is a
canonically defined smooth map of orbifolds from $\G\times M$ to $M$,
which induces $(g,\Phi)\mapsto\Phi\circ\iota(g^{-1})$ between the underlying
spaces. Clearly, the corresponding map $\jmath:\G\rightarrow [M;M]^r$ satisfies
$\jmath(1)=Id$, $\jmath(g_1g_2)=\jmath(g_1)\circ\jmath(g_2)$ for any
$g_1,g_2\in\G$. If furthermore, when $M$ is a genuine orbifold, $M$ is
compact and is an orbifold in the classical sense, and $\jmath$ is injective,
then $\jmath$ defines a smooth, effective left action of $\G$ on $M$
of reparametrization. (Note that when $M$ is a manifold, the canonical
map $\G\times M\rightarrow M$ defines a smooth action of $\G$ on $M$
in the usual sense.)

\subsection{Proof of Theorem 1.5}

We begin with some basic facts about complete Riemannian
orbifolds, which are straightforward generalizations of the
corresponding results for smooth manifolds.

For any connected Riemannian orbifold $X$, one can define a distance
function $d$ on it as follows. For any $p,q\in X$, the distance $d(p,q)$
between $p$ and $q$ is the infinimum of the lengths of all piecewise
$C^1$ paths in $X$ joining $p$ and $q$. By definition a piecewise $C^1$
path in $X$ joining $p$ and $q$ is a map $[a,b]\rightarrow X$, represented
by a homomorphism $\sigma=(\{f_\alpha\},\{\rho_{\beta\alpha}\})$ where
each $f_\alpha:(a_\alpha,b_\alpha)\rightarrow\widehat{U_\alpha}$ is
piecewise $C^1$, such that for the induced map $f:[a,b]\rightarrow X$
between the underlying spaces, $f(a)=p, f(b)=q$. Its length is
defined to be the integral $\int_{a}^b|\frac{df_\alpha}{dt}|dt$, where
the norm $|\frac{df_\alpha}{dt}|$ is defined with respect to the given
Riemannian
metric on $X$. A $C^1$ path in $X$ is a parametric geodesic if it is
represented by a homomorphism $\sigma=(\{f_\alpha\},\{\rho_{\beta\alpha}\})$
where each $f_\alpha:(a_\alpha,b_\alpha)\rightarrow\widehat{U_\alpha}$ is
a parametric geodesic. A Riemannian metric on an orbifold is called
`complete' if every parametric geodesic in $X$ has $(-\infty,\infty)$
as its maximal domain.

The following basic results concerning complete Riemannian metrics on
smooth manifolds, cf. e.g. pp. 172 in \cite{KN}, are readily extended
to the orbifold case.

\begin{itemize}
\item [{(1)}] For any two points in a connected complete Riemannian
manifold, there is a minimizing geodesic connecting them.
\item [{(2)}] The following are equivalent:
\begin{itemize}
\item The Riemannian metric is complete.
\item The metric space defined by the distance function $d$ is complete.
\item Every bounded subset is precompact.
\end{itemize}
\end{itemize}

Now let $X,X^\prime$ be any complete Riemannian orbifolds where $X$
is also compact. We fix some distance functions $d$, $d^\prime$ on
$X$ and $X^\prime$ respectively. Given any $\Phi\in [X;X^\prime]^r$,
which is represented by a $C^r$ admissible homomorphism $\sigma=
(\{f_\alpha\},\{\rho_{\beta\alpha}\})$, the differential $d\Phi$ is
a $C^{r-1}$ section of the orbifold vector bundle $\mbox{End}(TX,
E_\sigma)\rightarrow X$, where $E_\sigma\rightarrow X$ is the canonical
pull-back $C^r$ orbifold vector bundle constructed in Lemma 3.3.2 with
$E=TX^\prime$. Note that $\mbox{End}(TX,E_\sigma)\rightarrow X$ inherits
a natural metric from both $X$ and $X^\prime$, hence $d\Phi$ as a
$C^{r-1}$ section has a natural $C^{r-1}$-norm, which is denoted
by $||d\Phi||_{C^{r-1}}$, and is simply $\max_\alpha ||df_\alpha||
_{C^{r-1}}$ in terms of $\sigma=(\{f_\alpha\},\{\rho_{\beta\alpha}\})$.
Now pick an arbitrary point $p_0^\prime\in X^\prime$, and denote by
$f:X\rightarrow X^\prime$ the induced map of $\Phi$ between the
underlying spaces. We define the $C^r$-norm of $\Phi$ by
$$
||\Phi||_{C^r}=\max_{p\in X} d^\prime(p_0^\prime,f(p))+
||d\Phi||_{C^{r-1}}. \leqno (3.4.1)
$$

With the proceeding understood, we have

\begin{thm}
Let $\{\Phi_n\mid n\geq 1\}$ be any sequence of $C^r$-maps such that
$||\Phi_n||_{C^r}\leq C$ for some constant $C>0$ independent of $n$.
Then there exists a $C^{r-1}$-map $\Phi_0$, and a subsequence $\{\Phi_{n_i}\}$
of $\{\Phi_n\}$ such that $\Phi_{n_i}$ converges to $\Phi_0$ in the
$C^{r-1}$-topology.
\end{thm}

\pf
Let $f_n:X\rightarrow X^\prime$ be the induced map of $\Phi_n$
between the underlying spaces. The assumption that $||\Phi_n||_{C^r}\leq C$
(with $r\geq 1$) for some constant $C>0$ independent of $n$ implies that
all $f_n(X)$ are contained in a fixed bounded subset of $X^\prime$,
which is precompact by the completeness of the Riemannian metric on $X^\prime$,
and moreover, there exists a constant $C^\prime>0$ independent of $n$ such that
$d^\prime(f_n(p),f_n(q))\leq C^\prime d(p,q)$, $\forall p,q\in X,
\forall n$. By the classical Arzela-Ascoli theorem, there is
a subsequence of $\{f_n\}$, still denoted by $\{f_n\}$ for
simplicity, and there is a continuous map $f_0:X\rightarrow
X^\prime$, such that $f_n$ converges to $f_0$ in the $C^0$-topology.

The image $f_0(X)$ is a compact subset of $X^\prime$. Thus we can
cover $f_0(X)$ by a finite set $\{U_{\alpha^\prime}^\prime\}$ of
local charts on $X^\prime$. Consider the cover of open subsets
$\{f_0^{-1}(U_{\alpha^\prime}^\prime)\}$ of $X$. There exists a
finite cover $\{U_\alpha\}$ of local charts on $X$, such that
each $U_\alpha$ is admissible and the closure of
$\widehat{U_\alpha}$, $\widehat{\overline{U_\alpha}}$,
is a closed ball in the Euclidean space. Moreover,
$\{\overline{U_\alpha}\}$ is a refinement of
$\{f_0^{-1}(U_{\alpha^\prime}^\prime)\}$, which in particular
means that there is a correspondence $U_\alpha\mapsto U_\alpha^\prime\in
\{U_{\alpha^\prime}^\prime\}$ such that $f_0(\overline{U_\alpha})\subset
U_{\alpha}^\prime$. For sufficiently large $n$, we have $f_n(\overline
{U_\alpha})\subset U_\alpha^\prime$ since $f_n$ converges to $f_0$
in $C^0$-topology.

A digression is in order. Suppose $\Phi:X\rightarrow X^\prime$ is
any map of orbispaces where $X=Y/G$, $X^\prime=Y^\prime/G^\prime$
are global quotients. Then by the covering space theory developed
in the sequel \cite{C1}, $\Phi$ is represented by a pair $(f,\rho):(Y,G)
\rightarrow (Y^\prime,G^\prime)$, where $f$ is $\rho$-equivariant,
if and only if under the induced homomorphism $\Phi_\ast:\pi_1(X)
\rightarrow \pi_1(X^\prime)$, the subgroup $\pi_1(Y)\subset\pi_1(X)$
is sent into the subgroup $\pi_1(Y^\prime)\subset\pi_1(X^\prime)$,
cf. \S 2.4 of \cite{C1}. End of digression.

Observe that because of $f_n(\overline{U_\alpha})\subset
U_\alpha^\prime$ for sufficiently large $n$, the restriction of
$\Phi_n$ to the subspace
$\overline{U_\alpha}$ is a map into the subspace $U_\alpha^\prime$
of $X^\prime$. Since each $\widehat{\overline{U_\alpha}}$ is a
closed ball, we have $(\Phi_n)_\ast(\pi_1(\widehat{\overline{U_\alpha}}))
=\{1\}\subset \pi_1(\widehat{U_\alpha^\prime})$. Hence
the restriction of each $\Phi_n$ to $\overline{U_\alpha}$ for $n$
sufficiently large is represented by a pair $(f_\alpha^{(n)},
\rho_\alpha^{(n)}):(\widehat{\overline{U_\alpha}},G_{U_\alpha})
\rightarrow (\widehat{U_\alpha^\prime},G_{U_\alpha^\prime})$,
where $f_\alpha^{(n)}$ is $\rho_\alpha^{(n)}$-equivariant.

By Lemma 3.1.3, there exist mappings
$\rho_{\beta\alpha}^{(n)}:T(U_\alpha,U_\beta)\rightarrow T(U_\alpha^\prime,
U_\beta^\prime)$ with
$\rho_{\alpha\alpha}^{(n)}=\rho_\alpha^{(n)}$, such that $\sigma_n
=(\{f_\alpha^{(n)}\},\{\rho_{\beta\alpha}^{(n)}\})$ is a $C^r$
admissible homomorphism representing the map $\Phi_n$. On the
other hand, there is an infinite sequence $n_i\rightarrow\infty$
such that $\rho_{\beta\alpha}^{(n_i)}=\rho_{\beta\alpha}$ for
some $\rho_{\beta\alpha}:T(U_\alpha,U_\beta)\rightarrow T(U_\alpha^\prime,
U_\beta^\prime)$ independent of $n_i$ for all
indexes $\alpha,\beta$, because there are only finitely many
indexes and each $\rho_{\beta\alpha}^{(n)}$ is a mapping between
two finite sets independent of $n$. Finally, since for each $\alpha$,
$\{f_\alpha^{(n_i)}\}$
has bounded $C^r$-norms, we apply the classical Arzela-Ascoli
theorem to conclude that there is a $C^{r-1}$-map $f_\alpha^{(0)}:
\widehat{\overline{U_\alpha}}\rightarrow\widehat{U_\alpha^\prime}$
such that a subsequence of $f_\alpha^{(n_i)}$, still denoted by
$f_\alpha^{(n_i)}$, converges to $f_\alpha^{(0)}$ in $C^{r-1}$-topology.

To see that $(\{f^{(0)}_\alpha\},\{\rho_{\beta\alpha}\})$ is a
homomorphism, it suffices to check that $\rho_{\gamma\alpha}
(\eta\circ\xi({\bf a}))=\rho_{\gamma\beta}(\eta)\circ\rho_{\beta\alpha}(\xi)
(\underline{\{f_\alpha^{(0)}\}}({\bf a}))$ for
any ${\bf a}\in\Lambda(\xi,\eta)$. When $n_i$ is sufficiently large,
$\underline{\{f_\alpha^{(0)}\}}({\bf a})=\underline{\{f_\alpha^{(n_i)}\}}
({\bf a})$, from which the above equation follows.

Putting everything together, we may conclude that the sequence
$\sigma_{n_i}=(\{f_\alpha^{(n_i)}\},\{\rho_{\beta\alpha}\})$ of
$C^r$ admissible homomorphisms converges to the $C^{r-1}$
admissible homomorphism $\sigma_0=(\{f_\alpha^{(0)}\},\{\rho_{\beta\alpha}\})$
in $\O^{r-1}_{\{\rho_{\beta\alpha}\}}$. This means that the
corresponding maps $\Phi_{n_i}=[\sigma_{n_i}]$ converges to a $C^{r-1}$-map 
$\Phi_0=[\sigma_0]$ in the $C^{r-1}$-topology.

\hfill $\Box$

\vspace{2mm}

{\Small Current Address: Department of Mathematics and Statistics, 
University of Massachusetts at Amherst, Amherst, MA 01003. 
{\it e-mail:} wchen@@math.umass.edu}

\end{document}